\documentclass[11pt]{amsart}
\usepackage{amssymb,mathrsfs}
\newtheorem{theorem}{Theorem}[section]
\newtheorem{proposition}[theorem]{Proposition}
\newtheorem{lemma}[theorem]{Lemma}
\newtheorem{corollary}[theorem]{Corollary}
\newtheorem{conjecture}[theorem]{Conjecture}
\theoremstyle{definition}
\newtheorem{definition}[theorem]{Definition}
\newtheorem{remark}[theorem]{Remark}

\begin{document}

\title[Uniqueness of compact ancient solutions]{Uniqueness of compact ancient solutions to three-dimensional Ricci flow}
\author{Simon Brendle, Panagiota Daskalopoulos, and Natasa Sesum}
\address{Department of Mathematics \\ Columbia University \\ New York NY 10027}
\address{Department of Mathematics \\ Columbia University \\ New York NY 10027}
\address{Department of Mathematics \\ Rutgers University \\ Piscataway NJ 08854}
\begin{abstract}
In this paper, we study the classification of $\kappa$-noncollapsed ancient solutions to three-dimensional Ricci flow on $S^3$. We prove that such a solution is either isometric to a family of shrinking round spheres, or the Type II ancient solution constructed by Perelman. 
\end{abstract}
\thanks{The first author was supported by the National Science Foundation under grant DMS-1806190 and by the Simons Foundation. The second author was supported by the National Science Foundation under grant DMS-1266172. The third author was supported by the National Science Foundation under grants DMS-1056387 and DMS-1811833.} 
\maketitle 

\tableofcontents

\section{Introduction}

Consider a solution to the Ricci flow $\frac{\partial}{\partial t} g(t) = -2 \, \text{\rm Ric}_{g(t)}$ on a compact three-manifold which exists for all $t \in (-\infty,T]$. Such a solution is called an ancient solution. The goal in this work is to provide a classification of such solutions under natural geometric assumptions. 

Ancient compact solutions to the two-dimensional Ricci flow were classified by Daskalopoulos, Hamilton, and \v Se\v sum \cite{Daskalopoulos-Hamilton-Sesum2}. It turns out that in this case, the complete list contains only the shrinking spheres (which are noncollapsed) and the King solution (which is collapsed). The King solution is not self-similar, but it can be written in closed form. It was first discovered by King \cite{King} in the context of the logarithmic fast diffusion equation on $\mathbb{R}^2$ and later independently by Rosenau \cite{Rosenau} in the same context. It also appears as the sausage model in the context of quantum field theory (see \cite{Fateev-Onofri-Zamolodchikov}). Let us remark that the classification work in \cite{Daskalopoulos-Hamilton-Sesum2} classifies both collapsed and noncollapsed solutions.
  
We now turn our attention to the three-dimensional Ricci flow. In \cite{Brendle-Huisken-Sinestrari}, it was shown that any three-dimensional ancient solution on $S^3$ with uniformly pinched curvature is a family of shrinking round spheres. In \cite{Perelman1}, Perelman established the existence of a rotationally symmetric ancient solution on $S^3$ which is $\kappa$-noncollapsed and which is not a soliton. This ancient solution is of Type II backwards in time, namely its scalar curvature satisfies $\limsup_{t \to -\infty} (-t) \, R_{\text{\rm max}}(t) = \infty$. Going forward in time, the solution forms a Type I singularity, and shrinks to a round point. Perelman's ancient solution has backward limits which are either the Bryant soliton or the round cylinder $S^2 \times \mathbb{R}$, depending on how we choose the sequence of points about which we rescale. Perelman's ancient solution can be viewed as the three-dimensional analogue of the King solution. However, unlike the King solution, Perelman's ancient solution is noncollapsed.

The noncollapsing property plays a crucial role in the study of the Ricci flow. In fact, in \cite{Perelman1} Perelman proved that every ancient solution arising as a blow-up limit at a finite-time singularity on a compact manifold is $\kappa$-noncollapsed for some $\kappa>0$. Moreover, in dimension $3$, the well-known Hamilton-Ivey pinching estimate tells us that any such blow-up limit has nonnegative sectional curvature. It follows from Hamilton's Harnack estimate (see in \cite{Hamilton2}) that $R_t \geq 0$, yielding the existence of a uniform constant $C > 0$ so that $R(\cdot,t) \leq C$, for all $t\in (-\infty,T]$.  Since the curvature is positive, one concludes that 
$ \|\text{\rm Rm}\|_{g(t)} \leq C$, for all  $t \leq T$, for a  uniform constant $C$. Following Perelman, we say that $(M,g(t))$ is an ancient $\kappa$-solution if $(M,g(t))$ is defined on $(-\infty,T]$, is  non-flat and $\kappa$-noncollapsed, and has bounded nonnegative curvature. In \cite{Perelman1}, Perelman proposed the following conjecture:

\begin{conjecture}[Perelman \cite{Perelman1}]
\label{perelman.conjecture.noncompact.case}
Let $(M,g(t))$ be a noncompact ancient $\kappa$-solution to the Ricci flow in dimension $3$ with positive curvature. Then $(M,g(t))$ is the Bryant soliton.
\end{conjecture}

This conjecture was proved in \cite{Brendle1} in the class of steady gradient Ricci solitons, and in full generality in \cite{Brendle2}. The proof in \cite{Brendle2} has two main parts. In a first step, it is shown that the Bryant soliton is the only noncompact ancient $\kappa$-solution which has positive curvature and is rotationally symmetric. In a second step, it is shown that every noncompact ancient $\kappa$-solution with positive curvature must be rotationally symmetric.

The following is the analogue of Perelman's conjecture in the compact setting: 

\begin{conjecture}
\label{perelman.conjecture.compact.case}
Let $(S^3,g(t))$ be a compact ancient $\kappa$-solution to the Ricci flow on $S^3$. Then $g(t)$ is either a family of shrinking spheres or Perelman's ancient solution.
\end{conjecture}

As announced in \cite{Brendle2}, the techniques in that paper can also be applied to show that any ancient $\kappa$-solution on $S^3$ is rotationally symmetric. We include the proof of this fact in Section \ref{rotational.symmetry}. 

\begin{theorem}
\label{rotational.symmetry}
Let $(S^3,g(t))$ be an ancient $\kappa$-solution on $S^3$. Then $(S^3,g(t))$ is rotationally symmetric.
\end{theorem}

After this result was announced in \cite{Brendle2}, an alternative approach to Theorem \ref{rotational.symmetry} was proposed in \cite{Bamler-Kleiner}.

Next, we give a complete classification of all ancient $\kappa$-solutions on $S^3$ with rotational symmetry: 

\begin{theorem}
\label{uniqueness.theorem}
Let $(S^3,g_1(t))$ and $(S^3,g_2(t))$ be two ancient $\kappa$-solutions on $S^3$ which are rotationally symmetric. Assume that neither $(S^3,g_1(t))$ nor $(S^3,g_2(t))$ is a family of shrinking round spheres. Then $(S^3,g_1(t))$ and $(S^3,g_2(t))$ coincide up to a reparametrization in space, a translation in time, and a parabolic rescaling.
\end{theorem}

Combining Theorem \ref{rotational.symmetry} and Theorem \ref{uniqueness.theorem}, we can draw the following conclusion:

\begin{theorem} 
Let $(S^3,g(t))$ be an ancient $\kappa$-solution on $S^3$ which is not a family of shrinking round spheres. Then $(S^3,g(t))$ coincides with Perelman's solution up to diffeomorphisms, translations in time, and parabolic rescalings.
\end{theorem}

Let us mention some related work in the mean curvature flow setting. Compact, convex ancient solutions to the curve shortening flow were classified in \cite{Daskalopoulos-Hamilton-Sesum1}. In \cite{Brendle-Choi1},\cite{Brendle-Choi2}, it was shown that the bowl soliton is the only ancient solution to mean curvature flow which is noncompact, noncollapsed, strictly convex, and uniformly two-convex. In \cite{Angenent-Daskalopoulos-Sesum1},\cite{Angenent-Daskalopoulos-Sesum2}, it was shown that every ancient solution to mean curvature flow which is compact, noncollapsed, strictly convex, and uniformly two-convex is either the family of shrinking spheres or the family of ancient ovals constructed by White (cf. \cite{White}) and Haslhofer-Hershkovits (cf. \cite{Haslhofer-Hershkovits}). Collapsed ancient solutions to mean curvature flow were studied in \cite{Bourni-Langford-Tinaglia}.

The outline of the paper is as follows. In Section \ref{structure}, we recall some qualitative properties of ancient $\kappa$-solutions on $S^3$. In particular, an ancient $\kappa$-solution on $S^3$ is either a family of shrinking round spheres, or it has the structure of two caps joined by a tube (in which the solution is nearly cylindrical). In Section \ref{proof.of.rotational.symmetry}, we give the proof of Theorem \ref{rotational.symmetry}. 

In Section \ref{a.priori.estimates}, we establish various a-priori estimates for rotationally symmetric solutions, building on our earlier work \cite{Angenent-Brendle-Daskalopoulos-Sesum}.

In Section \ref{weights}, we introduce two weight functions $\mu_+(\rho,\tau)$ and $\mu_-(\rho,\tau)$ in the tip regions. This will be needed to prove the tip region estimates (see Proposition \ref{estimate.for.difference.in.tip.region}). 
 
In Section \ref{overview}, we give the proof of Theorem \ref{uniqueness.theorem}. The proof is inspired by the argument in \cite{Angenent-Daskalopoulos-Sesum2}. Let us sketch the main ideas. Suppose we are given two ancient $\kappa$-solutions $(S^3,g_1(t))$ and $(S^3,g_2(t))$ which satisfy the assumptions of Theorem \ref{uniqueness.theorem}. These ancient solutions can be described by profile functions $F_1(z,t)$ and $F_2(z,t)$. As in \cite{Angenent-Brendle-Daskalopoulos-Sesum}, the profile function gives the radius of a sphere of symmetry which has signed distance $z$ from some reference point. The results in Section \ref{a.priori.estimates} (and our earlier results in \cite{Angenent-Brendle-Daskalopoulos-Sesum}) give precise asymptotic estimates for the profile functions $F_1(z,t)$ and $F_2(z,t)$. 

We work on a time interval $(-\infty,t_*]$, where $-t_*$ is very large. We introduce a three-parameter family of profile functions $F_2^{\alpha\beta\gamma}(z,t)$. These differ from the original profile function $F_2(z,t)$ by a change of the reference point (represented by the parameter $\alpha$); a translation in time (represented by the parameter $\beta$); and a parabolic dilation (represented by the parameter $\gamma$). Our goal is to show that there exists a time $t_*$ and parameters $\alpha,\beta,\gamma$ such that $F_1(z,t) = F_2^{\alpha\beta\gamma}(z,t)$ for all $t \in (-\infty,t_*]$. 

To prove this, we consider two regions, the tip region and the cylindrical region. Roughly speaking, the tip region consists of points in space-time where the radius of the sphere of symmetry is $\lesssim \theta \sqrt{-2t}$, while the cylindrical region consists of point in space-time where the radius of the sphere of symmetry is $\gtrsim \theta \sqrt{-2t}$. Here, $\theta$ is a small positive constant which will be chosen later.

The first main ingredient is a weighted estimate for the difference of two solutions in the tip region (see Proposition \ref{estimate.for.difference.in.tip.region}). This estimate uses the weight functions $\mu_+(\rho,\tau)$ and $\mu_-(\rho,\tau)$ introduced in Section \ref{weights}. The tip region estimate works as long as we choose the parameter $\theta$ small enough. From this point on, we fix $\theta$ sufficient sufficiently small, so that the tip region estimate holds. We next analyze the difference of the two solutions in the cylindrical region. To that end, it is useful to perform a rescaling. We define
\begin{align*} 
G_1(\xi,\tau) &:= e^{\frac{\tau}{2}} \, F_1(e^{-\frac{\tau}{2}} \xi,-e^{-\tau}) - \sqrt{2}, \\ 
G_2^{\alpha\beta\gamma}(\xi,\tau) &:= e^{\frac{\tau}{2}} \, F_2^{\alpha\beta\gamma}(e^{-\frac{\tau}{2}} \xi,-e^{-\tau}) - \sqrt{2}, 
\end{align*} 
and 
\[H^{\alpha\beta\gamma}(\xi,\tau) := G_1(\xi,\tau) - G_2^{\alpha\beta\gamma}(\xi,\tau).\] 
Moreover, we introduce a cutoff in space which allows us to localize the function $H^{\alpha\beta\gamma}$ to the cylindrical region. We put 
\[H_{\mathcal{C}}^{\alpha\beta\gamma}(\xi,\tau) := \chi_{\mathcal{C}}((-\tau)^{-\frac{1}{2}} \xi) \, H^{\alpha\beta\gamma}(\xi,\tau),\] 
where $\chi_{\mathcal{C}}$ denotes a smooth, even cutoff function satisfying $\chi_{\mathcal{C}} = 1$ on $[0,\sqrt{4-\frac{\theta^2}{2}}]$ and $\chi_{\mathcal{C}} = 0$ on $[\sqrt{4-\frac{\theta^2}{4}},\infty)$. Finally, we write $t_* = -e^{-\tau_*}$. 

Given a time $\tau_*$, we choose the parameters $\alpha,\beta,\gamma$ such that the function $H_{\mathcal{C}}^{\alpha\beta\gamma}(\cdot,\tau_*)$ is orthogonal to the Hermite polynomials of degree $0$, $1$, and $2$. This gives three orthogonality relations for the three parameters $\alpha,\beta,\gamma$. Note that the orthogonality relations depend on $\theta$, but this does not pose a problem, as we have already fixed $\theta$ at this stage. 

For this choice of the parameters $\alpha,\beta,\gamma$, we are able to prove an energy estimate for the function $H_{\mathcal{C}}^{\alpha\beta\gamma}$ (see Proposition \ref{estimate.for.difference.in.cylindrical.region}). The estimate in Proposition \ref{estimate.for.difference.in.cylindrical.region} contains a term 
\[\sup_{\tau \leq \tau_*} \int_{\tau-1}^\tau \int_{\{\sqrt{4-\frac{\theta^2}{2}} \, (-\tau')^{\frac{1}{2}} \leq |\xi| \leq \sqrt{4-\frac{\theta^2}{4}} \, (-\tau')^{\frac{1}{2}}\}} e^{-\frac{\xi^2}{4}} \, H^{\alpha\beta\gamma}(\xi,\tau')^2 \, d\xi \, d\tau'\] 
which arises from the cutoff. Crucially, this term can be controlled using the tip region estimate (cf. Proposition \ref{estimate.for.difference.in.tip.region}). The upshot is that, for this particular choice of $\alpha,\beta,\gamma$, we can control the function $H_{\mathcal{C}}^{\alpha\beta\gamma}$ in terms of a scalar function $a^{\alpha\beta\gamma}$, which represents the orthogonal projection of $H_{\mathcal{C}}^{\alpha\beta\gamma}$ to the Hermite polynomial of degree $2$ (see Proposition \ref{neutral.mode.dominates}). 

The final ingredient is an ODE for the function $a^{\alpha\beta\gamma}$ (see Proposition \ref{ode.for.a}). Using this ODE together with the relation $a^{\alpha\beta\gamma}(\tau_*) = 0$, we can conclude that $a_{\alpha\beta\gamma}(\tau) = 0$ for all $\tau \in (-\infty,\tau_*]$. From this, we deduce that $G_1(\xi,\tau) = G_2^{\alpha\beta\gamma}(\xi,\tau)$ for all $(-\infty,\tau_*]$. This finally gives $F_1(z,t) = F_2^{\alpha\beta\gamma}(z,t)$ for all $t \in (-\infty,t_*]$.

\section{Structure of compact ancient $\kappa$-solutions}

\label{structure}

In this section, we recall some basic facts about the structure of compact ancient $\kappa$-solutions. Throughout this section, we assume that $(M,g(t))$ is a three-dimensional ancient $\kappa$-solution which is compact and simply connected. Moreover, we assume that $(M,g(t))$ is not a family of shrinking round spheres. Note that $M$ is diffeomorphic to $S^3$ by work of Hamilton.

\begin{proposition}
\label{asymptotic.soliton}
The asymptotic shrinking soliton associated with $(M,g(t))$ is isometric to the cylinder $S^2 \times \mathbb{R}$.
\end{proposition}

\textbf{Proof.} 
By Perelman's classification of shrinking gradient Ricci solitons in dimension $3$ (cf. \cite{Perelman2}), the asymptotic shrinking soliton associated with $(M,g(t))$ either has constant sectional curvature, or it locally splits as a product. If the asymptotic shrinking soliton associated with $(M,g(t))$ has constant sectional curvature, then, by Hamilton's curvature pinching estimates, the solution $(M,g(t))$ has constant sectional curvature for each $t$, contrary to our assumption. Therefore, the asymptotic shrinking soliton associated with $(M,g(t))$ must be isometric to either the cylinder $S^2 \times \mathbb{R}$, or a quotient of the cylinder $S^2 \times \mathbb{R}$. The asymptotic shrinking soliton cannot be a compact quotient of $S^2 \times \mathbb{R}$. Furthermore, if the asymptotic shrinking soliton is isometric to a $\mathbb{Z}_2$-quotient of the cylinder $S^2 \times \mathbb{R}$, then it contains an embedded $\mathbb{RP}^2$, but this is impossible since $M$ is diffeomorphic to $S^3$. Therefore, the asymptotic shrinking soliton must be isometric to the cylinder $S^2 \times \mathbb{R}$. \\

\begin{proposition}
\label{classification.of.limit.flows}
Let $(x_k,t_k)$ be an arbitrary sequence of points in space-time satisfying $\lim_{k \to \infty} t_k = -\infty$. Let us perform a parabolic rescaling around the point $(x_k,t_k)$ by the factor $R(x_k,t_k)$. After passing to a subsequence, the rescaled flows converge to a limit which is either a family of shrinking cylinders or the Bryant soliton. 
\end{proposition}

\textbf{Proof.} 
By Perelman's work \cite{Perelman1}, the rescaled manifolds converge to an ancient $\kappa$-solution. If the limiting ancient solution is noncompact, then, by \cite{Brendle2}, it must be either a family of shrinking cylinders or the Bryant soliton, and we are done. Hence, it remains to consider the case when the limiting ancient solution is compact. In this case, we have 
\[\limsup_{k \to \infty} R_{\text{\rm max}}(t_k) \, \text{\rm diam}_{g(t_k)}(M)^2 < \infty.\] 
This implies that $(M,g(t_k))$ cannot contain arbitrarily long necks. On the other hand, since the asymptotic shrinking soliton is a cylinder by Proposition \ref{asymptotic.soliton}, we know that $(M,g(t_k))$ must contain arbitrarily long necks if $k$ is sufficiently large. This is a contradiction. This completes the proof of Proposition \ref{classification.of.limit.flows}. \\

In the next step, we fix a small number $\varepsilon_1>0$. For later purposes, it is important that we choose $\varepsilon_1$ small enough so that the conclusion of the Neck Improvement Theorem in \cite{Brendle2} holds. Moreover, we fix a small number $\theta>0$ with the following property: if $(x,t)$ is a point in space-time satisfying $\lambda_1(x,t) \leq \theta R(x,t)$, then the point $(x,t)$ lies at the center of an evolving $\varepsilon_1$-neck. Here, $\lambda_1(x,t)$ denotes the smallest eigenvalue of the Ricci tensor at $(x,t)$. 

\begin{definition}
\label{definition.of.tips}
We say that $p$ is a tip of $(M,g(t))$ if $\lambda_1(p,t) > \frac{1}{6} \, R(p,t)$ and $\nabla R(p,t) = 0$. 
\end{definition}

By work of Hamilton, every neck admits a canonical foliation by CMC spheres. This will be referred to as Hamilton's CMC foliation.

\begin{proposition}
\label{caps}
Consider a sequence of times $t_k \to -\infty$. If $k$ is sufficiently large, we can find two disjoint compact domains $\Omega_{1,k}$ and $\Omega_{2,k}$ with the following properties: 
\begin{itemize}
\item $\Omega_{1,k}$ and $\Omega_{2,k}$ are diffeomorphic to $B^3$.
\item For each point $x \in M \setminus (\Omega_{1,k} \cup \Omega_{2,k})$, we have $\lambda_1(x,t_k) < \theta R(x,t_k)$. In particular, the point $(x,t_k)$ lies at the center of an evolving $\varepsilon_1$-neck.
\item For each point $x \in \Omega_{1,k} \cup \Omega_{2,k}$, we have $\lambda_1(x,t_k) > \frac{1}{2} \, \theta R(x,t_k)$.
\item $\partial \Omega_{1,k}$ and $\partial \Omega_{2,k}$ are leaves of Hamilton's CMC foliation in $(M,g(t_k))$.
\item For each $k$, there exists a leaf $\Sigma_k$ of the CMC foliation with the property that $\Omega_{1,k}$ and $\Omega_{2,k}$ lie in different connected components of $M \setminus \Sigma_k$, and $\sup_{x \in \Sigma_k} \frac{\lambda_1(x,t_k)}{R(x,t_k)} \to 0$. 
\item The domains $(\Omega_{1,k},g(t_k))$ and $(\Omega_{2,k},g(t_k))$ converge to the corresponding subset of the Bryant soliton after rescaling.
\end{itemize}
\end{proposition}

\textbf{Proof.} 
By Proposition \ref{asymptotic.soliton}, the asymptotic shrinking soliton associated with $(M,g(t))$ is a cylinder. Hence, we can find a sequence of points $q_k \in M$ such that $\frac{\lambda_1(q_k,t_k)}{R(q_k,t_k)} \to 0$. In particular, $q_k$ lies at the center of an $\varepsilon_1$-neck if $k$ is sufficiently large. Let $\Sigma_k$ denote the center sphere of this neck. Since $M$ is diffeomorphic to $S^3$, the complement $M \setminus \Sigma_k$ has two connected components. Let us follow Hamilton's CMC foliation outward to either side of the neck, until we encounter a point where $\lambda_1 \geq \frac{2}{3} \, \theta R$. Therefore, we can find points $q_{1,k}$ and $q_{2,k}$ such that $\lambda_1(q_{1,k},t_k) = \frac{2}{3} \, \theta R(q_{1,k},t_k)$ and $\lambda_1(q_{2,k},t_k) = \frac{2}{3} \, \theta R(q_{2,k},t_k)$. Moreover, $q_{1,k}$ and $q_{2,k}$ lie in different connected components of $M \setminus \Sigma_k$. By our choice of $\theta$, the point $q_{1,k}$ lies at the center of an $\varepsilon_1$-neck, and $q_{2,k}$ also lies at the center of an $\varepsilon_1$-neck. Let $\Sigma_{1,k}$ denote the leaf of Hamilton's CMC foliation passing through $q_{1,k}$, and let $\Sigma_{2,k}$ denote the leaf of Hamilton's CMC foliation passing through $q_{2,k}$. Moreover, let $N_k$ denote the tube bounded by $\Sigma_{1,k}$ and $\Sigma_{2,k}$. Clearly, $\Sigma_k \subset N_k$, and $\lambda_1(x,t_k) < \theta R(x,t_k)$ for all $x \in N_k$.

If we rescale the flow around the point $(q_{1,k},t_k)$, then the rescaled flows must converge to the Bryant soliton by Proposition \ref{classification.of.limit.flows}. Consequently, there exists a compact domain $\Omega_{1,k}$ such that $\partial \Omega_{1,k} = \Sigma_{1,k}$, $\Omega_{1,k}$ is diffeomorphic to $B^3$, and $\lambda_1(x,t_k) > \frac{1}{2} \, \theta R(x,t_k)$ for all $x \in \Omega_{1,k}$. Similarly, there exists a compact domain $\Omega_{2,k}$ such that $\partial \Omega_{2,k} = \Sigma_{2,k}$, $\Omega_{2,k}$ is diffeomorphic to a ball, and $\lambda_1(x,t_k) > \frac{1}{2} \, \theta R(x,t_k)$ for all $x \in \Omega_{2,k}$. Since $\sup_{x \in \Sigma_k} \frac{\lambda_1(x,t_k)}{R(x,t_k)} \to 0$, it follows that $\Omega_{1,k} \cup \Omega_{2,k} \subset M \setminus \Sigma_k$ if $k$ is sufficiently large. Since $q_{1,k}$ and $q_{2,k}$ lie in different connected components of $M \setminus \Sigma_k$, we conclude that $\Omega_{1,k}$ and $\Omega_{2,k}$ are contained in different connected components of $M \setminus \Sigma_k$. In particular, $\Omega_{1,k}$ and $\Omega_{2,k}$ are disjoint. Finally, the complement $M \setminus (\Omega_{1,k} \cup \Omega_{2,k})$ is contained in the tube $N_k$; therefore, $\lambda_1(x,t_k) < \theta R(x,t_k)$ for all $x \in M \setminus (\Omega_{1,k} \cup \Omega_{2,k})$. This completes the proof of Proposition \ref{caps}. \\

\begin{corollary}
\label{number.of.tips}
If $k$ is sufficiently large, then the manifold $(M,g(t_k))$ has exactly two tips. One of these points lies in $\Omega_{1,k}$ and the other lies in $\Omega_{2,k}$. In particular, these points are contained in different connected components of $M \setminus \Sigma_k$.
\end{corollary}

\textbf{Proof.} 
On the Bryant soliton, the scalar curvature has exactly one critical point (namely, the tip), and this critical point is non-degenerate. By Proposition \ref{caps}, the domains $(\Omega_{1,k},g(t_k))$ converge to a domain in the Bryant soliton. Hence, if $k$ is sufficiently large, then the set $\{x \in \Omega_{1,k}: \lambda_1(x,t_k) > \frac{1}{6} \, R(x,t_k), \, \nabla R(x,t_k) = 0\}$ consists of exactly one element. An analogous argument shows that the set $\{x \in \Omega_{2,k}: \lambda_1(x,t_k) > \frac{1}{6} \, R(x,t_k), \, \nabla R(x,t_k) = 0\}$ consists of exactly one element. Since $\{x \in M: \lambda_1(x,t_k) > \frac{1}{6} \, R(x,t_k)\} \subset \Omega_{1,k} \cup \Omega_{2,k}$, we conclude that  $\{x \in M: \lambda_1(x,t_k) > \frac{1}{6} \, R(x,t_k), \, \nabla R(x,t_k) = 0\}$ consists of exactly two elements. \\

\begin{proposition}
\label{asymptotics.around.tips}
Consider a sequence of times $t_k \to -\infty$. Let $p_{1,t_k}$ and $p_{2,t_k}$ denote the tips in $(M,g(t_k))$. If we rescale the flow around $(p_{1,t_k},t_k)$ or $(p_{2,t_k},t_k)$, then the rescaled flows converge to the Bryant soliton in the Cheeger-Gromov sense.
\end{proposition} 

\textbf{Proof.} 
This follows immediately from Proposition \ref{classification.of.limit.flows}. \\

\begin{proposition}
\label{distance.between.tips}
Consider a sequence of times $t_k \to -\infty$. Let $p_{1,t_k}$ and $p_{2,t_k}$ denote the tips in $(M,g(t_k))$. Then $R(p_{1,t_k},t_k) \, d_{g(t_k)}(p_{1,t_k},p_{2,t_k})^2 \to \infty$ and $R(p_{2,k},t_k) \, d_{g(t_k)}(p_{1,t_k},p_{2,t_k})^2 \to \infty$.
\end{proposition}

\textbf{Proof.} 
Suppose that $\limsup_{k \to \infty} R(p_{1,t_k},t_k) \, d_{g(t_k)}(p_{1,t_k},p_{2,t_k})^2 < \infty$. By Corollary \ref{number.of.tips}, $p_{1,t_k}$ and $p_{2,t_k}$ are contained in different connected components of $M \setminus \Sigma_k$. Consequently, we can find a sequence of points $y_k \in \Sigma_k$ such that $\limsup_{k \to \infty} R(p_{1,t_k},t_k) \, d_{g(t_k)}(p_{1,t_k},y_k)^2 < \infty$. By Perelman's longrange curvature estimate, we obtain $\limsup_{k \to \infty} R(p_{1,t_k},t_k)^{-1} \, R(y_k,t_k) < \infty$. Putting these facts together yields   $\limsup_{k \to \infty} R(y_k,t_k) \, d_{g(t_k)}(p_{1,t_k},y_k)^2 < \infty$. On the other hand, $\frac{\lambda_1(y_k,t_k)}{R(y_k,t_k)} \to 0$ since $y_k \in \Sigma_k$. From this, we deduce that $\frac{\lambda_1(p_{1,t_k},t_k)}{R(p_{1,t_k},t_k)} \to 0$. This contradicts the fact that $\lambda_1(p_{1,t_k},t_k) > \frac{1}{6} \, R(p_{1,t_k},t_k)$ for each $k$. This completes the proof of Proposition \ref{distance.between.tips}. \\

\begin{proposition}
\label{points.far.away.from.tips.are.necklike}
Consider a sequence of points $(x_k,t_k)$ in spacetime such that $t_k \to -\infty$. Let $p_{1,t_k}$ and $p_{2,t_k}$ denote the tips of $(M,g(t_k))$. If both 
$R(p_{1,t_k},t_k) \, d_{g(t_k)}(p_{1,t_k},x_k)^2 \to \infty$ and $R(p_{2,t_k},t_k) \, d_{g(t_k)}(p_{2,t_k},x_k)^2 \to \infty$, then $\frac{\lambda_1(x_k,t_k)}{R(x_k,t_k)} \to 0$. 
\end{proposition}

\textbf{Proof.} 
Suppose that $\limsup_{k \to \infty} \frac{\lambda_1(x_k,t_k)}{R(x_k,t_k)} > 0$. Let us rescale the flow around the point $(x_k,t_k)$ by the factor $R(x_k,t_k)$, and pass to the limit as $k \to \infty$. By Proposition \ref{classification.of.limit.flows}, the limit must be the Bryant soliton. Consequently, there exists a sequence of points $y_k \in M$ such that $\lambda_1(y_k,t_k) > \frac{1}{6} \, R(y_k,t_k)$, $\nabla R(y_k,t_k) = 0$, and $\limsup_{k \to \infty} R(x_k,t_k) \, d_{g(t_k)}(x_k,y_k)^2 < \infty.$ Perelman's longrange curvature estimate implies  $\limsup_{k \to \infty} R(x_k,t_k)^{-1} \, R(y_k,t_k) < \infty$.  All the above together yield $\limsup_{k \to \infty} R(y_k,t_k) \, d_{g(t_k)}(x_k,y_k)^2 < \infty$. Since $R(p_{1,t_k},t_k) \, d_{g(t_k)}(p_{1,t_k},x_k)^2 \to \infty$ and $R(p_{2,t_k},t_k) \, d_{g(t_k)}(p_{2,t_k},x_k)^2 \to \infty$, it follows that $y_k \notin \{p_{1,t_k},p_{2,t_k}\}$ if $k$ is sufficiently large. Therefore, the set $\{x \in M: \lambda_1(x,t_k) > \frac{1}{6} \, R(x,t_k), \, \nabla R(x,t_k) = 0\}$ contains at least three elements if $k$ is sufficiently large. This contradicts Corollary \ref{number.of.tips}. This completes the proof of Proposition \ref{points.far.away.from.tips.are.necklike}. \\

By combining Corollary \ref{number.of.tips}, Proposition \ref{asymptotics.around.tips}, Proposition \ref{distance.between.tips}, and Proposition \ref{points.far.away.from.tips.are.necklike}, we can draw the following conclusion:

\begin{corollary}
\label{summary}
(i) If $-t$ is sufficiently large, then the manifold $(M,g(t))$ has exactly two tips $p_{1,t}$ and $p_{2,t}$, and these vary smoothly in $t$. \\ 
(ii) Suppose that a large number $A$ is given. If $-t$ is sufficiently large (depending on $A$), then the balls $B_{g(t)}(p_{1,t},A R(p_{1,t},t)^{-\frac{1}{2}})$ and $B_{g(t)}(p_{2,t},A R(p_{2,t},t)^{-\frac{1}{2}})$ are disjoint. \\
(iii) Suppose that a large number $A$ and a small number $\varepsilon>0$ are given. If $-t$ is sufficiently large (depending on $A$ and $\varepsilon$), then the solution in the ball $B_{g(t)}(p_{1,t},A R(p_{1,t},t)^{-\frac{1}{2}})$ is (after a suitable rescaling) $\varepsilon$-close to the corresponding piece of the Bryant soliton in the Cheeger-Gromov sense. Similarly, the solution in the ball $B_{g(t)}(p_{2,t},A R(p_{2,t},t)^{-\frac{1}{2}})$ is (after a suitable rescaling) $\varepsilon$-close to the corresponding piece of the Bryant soliton in the Cheeger-Gromov sense. \\ 
(iv) Given $\varepsilon>0$, we can find a time $T \in (-\infty,0]$ and a large constant $A$ with the following property. If $t \leq T$ and $x \notin B_{g(t)}(p_{1,t},A R(p_{1,t},t)^{-\frac{1}{2}}) \cup B_{g(t)}(p_{2,t},A R(p_{2,t},t)^{-\frac{1}{2}})$, then $(x,t)$ lies at the center of an evolving $\varepsilon$-neck.
\end{corollary}

\section{Rotational symmetry of compact ancient $\kappa$-solutions and proof of Theorem \ref{rotational.symmetry}}

\label{proof.of.rotational.symmetry}

In this section, we give the proof of rotational symmetry. Throughout this section, we assume that $(M,g(t))$ is a three-dimensional ancient $\kappa$-solution which is compact and simply connected. Moreover, we assume that $(M,g(t))$ is not a family of shrinking round spheres. We claim that $(M,g(t))$ is rotationally symmetric. The proof is by contradiction. \textbf{We will assume throughout this section that $(M,g(t))$ is not rotationally symmetric.} 

As in the previous section, we fix a small number $\varepsilon_1 > 0$ and a large number $L$ so that the conclusion of the Neck Improvement Theorem in \cite{Brendle2} holds. Moreover, we fix a small number $\theta>0$ with the following property: if $(x,t)$ is a point in space-time satisfying $\lambda_1(x,t) \leq \theta R(x,t)$, then the point $(x,t)$ lies at the center of an evolving $\varepsilon_1$-neck. 

We begin with a definition, which is adapted from \cite{Brendle2}:

\begin{definition}
\label{symmetry.of.cap}
We say that the flow is $\varepsilon$-symmetric at time $\bar{t}$ if there exist a compact domain $D \subset M$ and time-independent vector fields $U^{(1)},U^{(2)},U^{(3)}$ which are defined on an open set containing $D$ such that the following statements hold:
\begin{itemize} 
\item The domain $D$ is a disjoint union of two domains $D_1$ and $D_2$, each of which is diffeomorphic to $B^3$. 
\item $\lambda_1(x,\bar{t}) < \theta R(x,\bar{t})$ for all points $x \in M \setminus D$.
\item $\lambda_1(x,\bar{t}) > \frac{1}{2} \, \theta R(x,\bar{t})$ for all points $x \in D$.
\item $\partial D_1$ and $\partial D_2$ are leaves of Hamilton's CMC foliation of $(M,g(\bar{t}))$. 
\item For each $x \in M \setminus D$, the point $(x,\bar{t})$ is $\varepsilon$-symmetric in the sense of Definition 8.2 in \cite{Brendle2}.
\item $\sup_{D_1 \times [\bar{t}-\rho_1^2,\bar{t}]} \sum_{l=0}^2 \sum_{a=1}^3 \rho_1^{2l} \, |D^l(\mathscr{L}_{U^{(a)}}(g(t)))|^2 \leq \varepsilon^2$, where $\rho_1^{-2} := \sup_{x \in D_1} R(x,\bar{t})$.
\item $\sup_{D_2 \times [\bar{t}-\rho_2^2,\bar{t}]} \sum_{l=0}^2 \sum_{a=1}^3 \rho_2^{2l} \, |D^l(\mathscr{L}_{U^{(a)}}(g(t)))|^2 \leq \varepsilon^2$, where $\rho_2^{-2} := \sup_{x \in D_2} R(x,\bar{t})$.
\item If $\Sigma \subset D_1$ is a leaf of the CMC foliation of $(M,g(\bar{t}))$ satisfying $\sup_{x \in \Sigma} d_{g(\bar{t})}(x,\partial D_1) \leq 10 \, \text{\rm area}_{g(\bar{t})}(\partial D_1)^{\frac{1}{2}}$, then $\sup_\Sigma \sum_{a=1}^3 \rho_1^{-2} \, |\langle U^{(a)},\nu \rangle|^2 \leq \varepsilon^2$, where $\nu$ denotes the unit normal vector to $\Sigma$ in $(M,g(\bar{t}))$.
\item If $\Sigma \subset D_2$ is a leaf of the CMC foliation of $(M,g(\bar{t}))$ satisfying $\sup_{x \in \Sigma} d_{g(\bar{t})}(x,\partial D_2) \leq 10 \, \text{\rm area}_{g(\bar{t})}(\partial D_2)^{\frac{1}{2}}$, then $\sup_\Sigma \sum_{a=1}^3 \rho_2^{-2} \, |\langle U^{(a)},\nu \rangle|^2 \leq \varepsilon^2$, where $\nu$ denotes the unit normal vector to $\Sigma$ in $(M,g(\bar{t}))$.
\item If $\Sigma \subset D_1$ is a leaf of the CMC foliation of $(M,g(\bar{t}))$ satisfying $\sup_{x \in \Sigma} d_{g(\bar{t})}(x,\partial D_1) \leq 10 \, \text{\rm area}_{g(\bar{t})}(\partial D_1)^{\frac{1}{2}}$, then 
\[\sum_{a,b=1}^3 \bigg | \delta_{ab} - \text{\rm area}_{g(\bar{t})}(\Sigma)^{-2} \int_\Sigma \langle U^{(a)},U^{(b)} \rangle_{g(\bar{t})} \, d\mu_{g(\bar{t})} \bigg |^2 \leq \varepsilon^2.\]
\item If $\Sigma \subset D_2$ is a leaf of the CMC foliation of $(M,g(\bar{t}))$ satisfying $\sup_{x \in \Sigma} d_{g(\bar{t})}(x,\partial D_2) \leq 10 \, \text{\rm area}_{g(\bar{t})}(\partial D_2)^{\frac{1}{2}}$, then 
\[\sum_{a,b=1}^3 \bigg | \delta_{ab} - \text{\rm area}_{g(\bar{t})}(\Sigma)^{-2} \int_\Sigma \langle U^{(a)},U^{(b)} \rangle_{g(\bar{t})} \, d\mu_{g(\bar{t})} \bigg |^2 \leq \varepsilon^2.\]
\end{itemize}
\end{definition}

\begin{remark}
Each tip in $(M,g(\bar{t}))$ is contained in $D$. Moreover, the two tips lie in different connected components of $\{x \in M: \lambda_1(x,\bar{t}) > \frac{1}{2} \theta R(x,\bar{t})\}$; in particular, the tips lie in different connected components of $D$. Hence, after relabeling $D_1$ and $D_2$ if necessary, we have $p_{1,\bar{t}} \in D_1$ and $p_{2,\bar{t}} \in D_2$. With this understood, we have $\text{\rm diam}_{g(\bar{t})}(D_1) \leq C \, R(p_{1,\bar{t}},\bar{t})^{-\frac{1}{2}}$ and $\text{\rm diam}_{g(\bar{t})}(D_2) \leq C \, R(p_{2,\bar{t}},\bar{t})^{-\frac{1}{2}}$. This gives $\frac{1}{C} \, R(p_{1,\bar{t}},\bar{t}) \leq R(x,\bar{t}) \leq C \, R(p_{1,\bar{t}},\bar{t})$ for all $x \in D_1$, and $\frac{1}{C} \, R(p_{2,\bar{t}},\bar{t}) \leq R(x,\bar{t}) \leq C \, R(p_{2,\bar{t}},\bar{t})$ for all $x \in D_2$.
\end{remark}

\begin{lemma}
\label{openness.property} 
Suppose that the flow is $\varepsilon$-symmetric at time $\bar{t}$. If $\tilde{t}$ is sufficiently close to $\bar{t}$, then the flow is $2\varepsilon$-symmetric at time $\bar{t}$.
\end{lemma}

\textbf{Proof.} 
The proof is analogous to the proof of Lemma 9.5 in \cite{Brendle2}. \\

\begin{proposition}
\label{solution.becomes.arbitrarily.symmetric.as.t.approaches.minus.infinity}
Let $\varepsilon>0$ be given. If $-t$ is sufficiently large (depending on $\varepsilon$), then the flow is $\varepsilon$-symmetric at time $t$.
\end{proposition}

\textbf{Proof.} 
This follows from Corollary \ref{summary}. \\

We next consider an arbitrary sequence $\varepsilon_k \to 0$. For $k$ large, we define 
\[t_k = \inf \{t \in (-\infty,0]: \text{\rm The flow is not $\varepsilon_k$-symmetric at time $t$}\}.\]  
If $\limsup_{k \to \infty} t_k > -\infty$, it follows that $(M,g(t))$ is rotationally symmetric for $-t$ sufficiently large, and this contradicts our assumption. Therefore, $\limsup_{k \to \infty} t_k = -\infty$. 

For $-t$ sufficiently large, we denote by $p_{1,t}$ and $p_{2,t}$ the tips of $(M,g(t))$. Since $t_k \to -\infty$, Proposition \ref{asymptotics.around.tips} implies that, if we rescale the solution around $(p_{1,t_k},t_k)$ by the factor $r_{1,k}^{-2} := R(p_{1,t_k},t_k)$, then the rescaled flows converge to the Bryant soliton. Similarly, if we rescale the solution around $(p_{1,t_k},t_k)$ by the factor $r_{2,k}^{-2} := R(p_{2,t_k},t_k)$, then the rescaled flows converge to the Bryant soliton. 

Hence, we can draw the following conclusion:

\begin{proposition} 
\label{perturbation.of.Bryant.soliton}
There exists a sequence $\delta_k \to 0$ such that the following statements hold when $k$ is sufficiently large: 
\begin{itemize}
\item For each $t \in [t_k-\delta_k^{-1} r_{1,k}^2,t_k]$, we have $d_{g(t)}(p_{1,t_k},p_{1,t}) \leq \delta_k r_{1,k}$ and $1-\delta_k \leq r_{1,k}^2 \, R(p_{1,t},t) \leq 1+\delta_k$.
\item For each $t \in [t_k-\delta_k^{-1} r_{2,k}^2,t_k]$, we have $d_{g(t)}(p_{2,t_k},p_{2,t}) \leq \delta_k r_{2,k}$ and $1-\delta_k \leq r_{2,k}^2 \, R(p_{2,t},t) \leq 1+\delta_k$.
\item The scalar curvature satisfies $r_{1,k}^2 \, R(x,t) \leq 4$ and $\frac{1}{2K} \, (r_{1,k}^{-1} \, d_{g(t)}(p_{1,t_k},x)+1)^{-1} \leq r_{1,k}^2 \, R(x,t) \leq 2K \, (r_{1,k}^{-1} \, d_{g(t)}(p_{1,t_k},x)+1)^{-1}$ for all points $(x,t) \in B_{g(t_k)}(p_{1,t_k},\delta_k^{-1} r_{1,k}) \times [t_k-\delta_k^{-1} r_{1,k}^2,t_k]$.
\item The scalar curvature satisfies $r_{2,k}^2 \, R(x,t) \leq 4$ and $\frac{1}{2K} \, (r_{2,k}^{-1} \, d_{g(t)}(p_{2,t_k},x)+1)^{-1} \leq r_{2,k}^2 \, R(x,t) \leq 2K \, (r_{2,k}^{-1} \, d_{g({2,t})}(p_{2,t_k},x)+1)^{-1}$ for all points $(x,t) \in B_{g(t_k)}(p_{2,t_k},\delta_k^{-1} r_{2,k}) \times [t_k-\delta_k^{-1} r_{2,k}^2,t_k]$.
\item There exists a nonnegative function $f_1: B_{g(t_k)}(p_{1,t_k},\delta_k^{-1} r_{1,k}) \times [t_k-\delta_k^{-1} r_{1,k}^2,t_k] \to \mathbb{R}$ such that $|\text{\rm Ric}-D^2 f_1| \leq \delta_k r_{1,k}^{-2}$, $|\Delta f_1 + |\nabla f_1|^2 - r_{1,k}^{-2}| \leq \delta_k r_{1,k}^{-2}$, and $|\frac{\partial}{\partial t} f_1 + |\nabla f_1|^2| \leq \delta_k r_{1,k}^{-2}$. Moreover, the function $f_1$ satisfies $\frac{1}{2K} \, (r_{1,k}^{-1} \, d_{g(t)}(p_{1,t_k},x)+1) \leq f_1(x,t)+1 \leq 2K \, (r_{1,k}^{-1} \, d_{g(t)}(p_{1,t_k},x)+1)$ for all points $(x,t) \in B_{g(t_k)}(p_{1,t_k},\delta_k^{-1} r_{1,k}) \times [t_k-\delta_k^{-1} r_{1,k}^2,t_k]$.
\item There exists a nonnegative function $f_2: B_{g(t_k)}(p_{2,t_k},\delta_k^{-1} r_{2,k}) \times [t_k-\delta_k^{-1} r_{2,k}^2,t_k] \to \mathbb{R}$ such that $|\text{\rm Ric}-D^2 f_2| \leq \delta_k r_{2,k}^{-2}$, $|\Delta f_2 + |\nabla f_2|^2 - r_{2,k}^{-2}| \leq \delta_k r_{2,k}^{-2}$, and $|\frac{\partial}{\partial t} f_2 + |\nabla f_2|^2| \leq \delta_k r_{2,k}^{-2}$. Moreover, the function $f_2$ satisfies $\frac{1}{2K} \, (r_{2,k}^{-1} \, d_{g(t)}(p_{2,t_k},x)+1) \leq f_2(x,t)+1 \leq 2K \, (r_{2,k}^{-1} \, d_{g(t)}(p_{2,t_k},x)+1)$ for all points $(x,t) \in B_{g(t_k)}(p_{2,t_k},\delta_k^{-1} r_{2,k}) \times [t_k-\delta_k^{-1} r_{2,k}^2,t_k]$.
\end{itemize}
Here, $K$ is a universal constant. 
\end{proposition}

\textbf{Proof.} 
By Proposition \ref{asymptotics.around.tips}, the solution looks like the Bryant soliton near each tip. From this, the assertion follows. \\

\begin{lemma}
\label{derivative.of.distance.function}
By a suitable choice of $\delta_k$, we can arrange that the following holds. If $t \in [t_k-\delta_k^{-1} r_{1,k}^2,t_k]$ and $d_{g(t)}(p_{1,t_k},x) \leq \delta_k^{-1} \, r_{1,k}$, then $0 \leq -\frac{d}{dt} d_{g(t)}(p_{1,t_k},x) \leq 80 \, r_{1,k}^{-1}$. Similarly, if $t \in [t_k-\delta_k^{-1} r_{2,k}^2,t_k]$ and $d_{g(t)}(p_{2,t_k},x) \leq \delta_k^{-1} \, r_{2,k}$, then $0 \leq -\frac{d}{dt} d_{g(t)}(p_{2,t_k},x) \leq 80 \, r_{2,k}^{-1}$.
\end{lemma} 

\textbf{Proof.} 
By Proposition \ref{asymptotics.around.tips}, the solution looks like the Bryant soliton near each tip. From this, the assertion follows. \\

\begin{lemma} 
\label{balls.around.tips.are.disjoint}
By a suitable choice of $\delta_k$, we can arrange that the balls $B_{g(t)}(p_{1,t},\delta_k^{-2} R(p_{1,t},t)^{-\frac{1}{2}})$ and $B_{g(t)}(p_{2,t},\delta_k^{-2} R(p_{2,t},t)^{-\frac{1}{2}})$ are disjoint for each $t \in (-\infty,t_k]$.
\end{lemma}

\textbf{Proof.} 
Since $t_k \to -\infty$, the assertion follows from Corollary \ref{summary}. \\

\begin{lemma}
\label{varepsilon_k.symmetry.at.earlier.times}
If $t \in (-\infty,t_k)$, then the flow is $\varepsilon_k$-symmetric at time $t$. In particular, if $(x,t) \in M \times (-\infty,t_k)$ is a point in spacetime satisfying $\lambda_1(x,t) < \frac{1}{2} \theta R(x,t)$, then the point $(x,t)$ is $\varepsilon_k$-symmetric in the sense of Definition 8.2 in \cite{Brendle2}.
\end{lemma}

\textbf{Proof.} 
The first statement follows directly from the definition of $t_k$. The second statement follows from Definition \ref{symmetry.of.cap}. \\

Recall that $L$ has been defined as the constant in the Neck Improvement Theorem in \cite{Brendle2}. By Corollary \ref{summary}, we can find a time $T \in (-\infty,0]$ and a large constant $\Lambda$ with the following properties: 
\begin{itemize} 
\item $L \sqrt{\frac{4K}{\Lambda}} \leq 10^{-6}$.
\item If $(\bar{x},\bar{t}) \in M \times (-\infty,T]$ satisfies $d_{g(\bar{t})}(p_{1,\bar{t}},\bar{x}) \geq \frac{\Lambda}{2} \, R(p_{1,\bar{t}},\bar{t})^{-\frac{1}{2}}$ and $d_{g(\bar{t})}(p_{2,\bar{t}},\bar{x}) \geq \frac{\Lambda}{2} \, R(p_{2,\bar{t}},\bar{t})^{-\frac{1}{2}}$, then $\lambda_1(x,t) < \frac{1}{2} \theta R(x,t)$ for all points $(x,t) \in B_{g(\bar{t})}(\bar{x},L \, R(\bar{x},\bar{t})^{-\frac{1}{2}}) \times [\bar{t} - L \, R(\bar{x},\bar{t})^{-1},\bar{t}]$.
\end{itemize}

\begin{lemma}
\label{improved.symmetry.for.far.away.points}
If $(\bar{x},\bar{t}) \in M \times (-\infty,t_k]$ satisfies $d_{g(\bar{t})}(p_{1,\bar{t}},\bar{x}) \geq \frac{\Lambda}{2} \, R(p_{1,\bar{t}},\bar{t})^{-\frac{1}{2}}$ and $d_{g(\bar{t})}(p_{2,\bar{t}},\bar{x}) \geq \frac{\Lambda}{2} \, R(p_{2,\bar{t}},\bar{t})^{-\frac{1}{2}}$, then $(\bar{x},\bar{t})$ is $\frac{\varepsilon_k}{2}$-symmetric.
\end{lemma}

\textbf{Proof.} 
By our choice of $\Lambda$, every point in the parabolic neighborhood $B_{g(\bar{t})}(\bar{x},L \, R(\bar{x},\bar{t})^{-\frac{1}{2}}) \times [\bar{t} - L \, R(\bar{x},\bar{t})^{-1},\bar{t}]$ satisfies $\lambda_1(x,t) < \frac{1}{2} \theta R(x,t)$. By our choice of $\theta$, every point in the parabolic neighborhood $B_{g(\bar{t})}(\bar{x},L \, R(\bar{x},\bar{t})^{-\frac{1}{2}}) \times [\bar{t} - L \, R(\bar{x},\bar{t})^{-1},\bar{t}]$ lies at the center of an evolving $\varepsilon_1$-neck. Moreover, by Lemma \ref{varepsilon_k.symmetry.at.earlier.times}, every point in the parabolic neighborhood $B_{g(\bar{t})}(\bar{x},L \, R(\bar{x},\bar{t})^{-\frac{1}{2}}) \times [\bar{t} - L \, R(\bar{x},\bar{t})^{-1},\bar{t})$ is $\varepsilon_k$-symmetric. Hence, the Neck Improvement Theorem in \cite{Brendle2} implies that $(\bar{x},\bar{t})$ is $\frac{\varepsilon_k}{2}$-symmetric. \\

\begin{lemma}
\label{technical}
If $(\bar{x},\bar{t}) \in M \times [t_k-\delta_k^{-1} r_{1,k}^2,t_k]$ satisfies $\Lambda r_{1,k} \leq d_{g(\bar{t})}(p_{1,t_k},\bar{x}) \leq \delta_k^{-1} \, r_{1,k}$, then $d_{g(\bar{t})}(p_{1,\bar{t}},\bar{x}) \geq \frac{\Lambda}{2} \, R(p_{1,\bar{t}},\bar{t})^{-\frac{1}{2}}$ and $d_{g(\bar{t})}(p_{2,\bar{t}},\bar{x}) \geq \frac{\Lambda}{2} \, R(p_{2,\bar{t}},\bar{t})^{-\frac{1}{2}}$. Similarly, if $(\bar{x},\bar{t}) \in M \times [t_k-\delta_k^{-1} r_{2,k}^2,t_k]$ satisfies $\Lambda r_{2,k} \leq d_{g(\bar{t})}(p_{2,t_k},\bar{x}) \leq \delta_k^{-1} \, r_{2,k}$, then $d_{g(\bar{t})}(p_{1,\bar{t}},\bar{x}) \geq \frac{\Lambda}{2} \, R(p_{1,\bar{t}},\bar{t})^{-\frac{1}{2}}$ and $d_{g(\bar{t})}(p_{2,\bar{t}},\bar{x}) \geq \frac{\Lambda}{2} \, R(p_{2,\bar{t}},\bar{t})^{-\frac{1}{2}}$.
\end{lemma}

\textbf{Proof.}
Suppose that $(\bar{x},\bar{t}) \in M \times [t_k-\delta_k^{-1} r_{1,k}^2,t_k]$ satisfies $\Lambda r_{1,k} \leq d_{g(\bar{t})}(p_{1,t_k},\bar{x}) \leq \delta_k^{-1} \, r_{1,k}$. Using Proposition \ref{perturbation.of.Bryant.soliton}, we obtain $(\Lambda-\delta_k) \, r_{1,k} \leq d_{g(\bar{t})}(p_{1,\bar{t}},\bar{x}) \leq (\delta_k^{-1}+\delta_k) \, r_{1,k}$ and $1-\delta_k \leq r_{1,k}^2 \, R(p_{1,\bar{t}},\bar{t}) \leq 1+\delta_k$. Putting these facts together, we obtain $\frac{\Lambda}{2} \, R(p_{1,\bar{t}},\bar{t})^{-\frac{1}{2}} \leq d_{g(\bar{t})}(p_{1,\bar{t}},\bar{x}) \leq 2\delta_k^{-1} \, R(p_{1,\bar{t}},\bar{t})^{-\frac{1}{2}}$ if $k$ is sufficiently large. By Lemma \ref{balls.around.tips.are.disjoint}, the balls $B_{g(\bar{t})}(p_{1,\bar{t}},\delta_k^{-2} R(p_{1,\bar{t}},\bar{t})^{-\frac{1}{2}})$ and $B_{g(\bar{t})}(p_{2,\bar{t}},\delta_k^{-2} R(p_{2,\bar{t}},\bar{t})^{-\frac{1}{2}})$ are disjoint. Since $d_{g(\bar{t})}(p_{1,\bar{t}},\bar{x}) < \delta_k^{-2} \, R(p_{1,\bar{t}},\bar{t})^{-\frac{1}{2}}$, we conclude that $d_{g(\bar{t})}(p_{2,\bar{t}},\bar{x}) \geq \delta_k^{-2} \, R(p_{2,\bar{t}},\bar{t})^{-\frac{1}{2}} \geq \frac{\Lambda}{2} \, R(p_{2,\bar{t}},\bar{t})^{-\frac{1}{2}}$. This proves the assertion. \\

\begin{proposition}
\label{iteration}
If $(\bar{x},\bar{t}) \in M \times [t_k-2^{-j} \delta_k^{-1} r_{1,k}^2,t_k]$ satisfies $2^{\frac{j}{400}} \, \Lambda r_{1,k} \leq d_{g(\bar{t})}(p_{1,t_k},\bar{x}) \leq (400KL)^{-j} \, \delta_k^{-1} \, r_{1,k}$, then $(\bar{x},\bar{t})$ is $2^{-j-1} \varepsilon_k$-symmetric. Similarly, if $(\bar{x},\bar{t}) \in M \times [t_k-2^{-j} \delta_k^{-1} r_{2,k}^2,t_k]$ satisfies $2^{\frac{j}{400}} \, \Lambda r_{2,k} \leq d_{g(\bar{t})}(p_{2,t_k},\bar{x}) \leq (400KL)^{-j} \, \delta_k^{-1} \, r_{2,k}$, then $(\bar{x},\bar{t})$ is $2^{-j-1} \varepsilon_k$-symmetric. 
\end{proposition}

\textbf{Proof.} 
The proof is by induction on $j$. We first verify the assertion for $j=0$. Suppose that $(\bar{x},\bar{t}) \in M \times [t_k-\delta_k^{-1} r_{1,k}^2,t_k]$ satisfies $\Lambda r_{1,k} \leq d_{g(\bar{t})}(p_{1,t_k},\bar{x}) \leq \delta_k^{-1} \, r_{1,k}$. By Lemma \ref{technical}, we know that $d_{g(\bar{t})}(p_{1,\bar{t}},\bar{x}) \geq \frac{\Lambda}{2} \, R(p_{1,\bar{t}},\bar{t})^{-\frac{1}{2}}$ and $d_{g(\bar{t})}(p_{2,\bar{t}},\bar{x}) \geq \frac{\Lambda}{2} \, R(p_{2,\bar{t}},\bar{t})^{-\frac{1}{2}}$. Hence, Lemma \ref{improved.symmetry.for.far.away.points} implies that $(\bar{x},\bar{t})$ is $\frac{\varepsilon_k}{2}$-symmetric. This proves the assertion for $j=0$. 

The inductive step is analogous to the proof of Proposition 9.16 in \cite{Brendle2}. Suppose that $j \geq 1$ and the assertion holds for $j-1$. We claim that the assertion holds for $j$. To that end, we consider a point $(\bar{x},\bar{t}) \in M \times [t_k-2^{-j} \delta_k^{-1} r_{1,k}^2,t_k]$ such that $2^{\frac{j}{400}} \, \Lambda r_{1,k} \leq d_{g(\bar{t})}(p_{1,t_k},\bar{x}) \leq (400KL)^{-j} \, \delta_k^{-1} \, r_{1,k}$. Lemma \ref{technical} implies that $d_{g(\bar{t})}(p_{1,\bar{t}},\bar{x}) \geq \frac{\Lambda}{2} \, R(p_{1,\bar{t}},\bar{t})^{-\frac{1}{2}}$ and $d_{g(\bar{t})}(p_{2,\bar{t}},\bar{x}) \geq \frac{\Lambda}{2} \, R(p_{2,\bar{t}},\bar{t})^{-\frac{1}{2}}$. In view of our definition of $\Lambda$, we obtain $\lambda_1(\bar{x},\bar{t}) < \frac{1}{2} \theta R(\bar{x},\bar{t})$. Hence, by our choice of $\theta$, $(\bar{x},\bar{t})$ lies at the center of an evolving $\varepsilon_1$-neck. For abbreviation, we put $R(\bar{x},\bar{t}) = r^{-2}$. By Proposition \ref{perturbation.of.Bryant.soliton}, $\frac{1}{4} \, r_{1,k}^2 \leq r^2 \leq 4K \, r_{1,k} \, d_{g(\bar{t})}(p_{1,t_k},\bar{x})$. This implies 
\begin{align*} 
\bar{t}-Lr^2 
&\geq \bar{t} - 4KL \, r_{1,k} \, d_{g(\bar{t})}(p_{1,t_k},\bar{x}) \\ 
&\geq \bar{t} - 4KL \, (400KL)^{-j} \, \delta_k^{-1} \, r_{1,k}^2 \\ 
&\geq \bar{t} - 2^{-j} \, \delta_k^{-1} \, r_{1,k}^2 \\ 
&\geq t_k - 2^{-j+1} \, \delta_k^{-1} \, r_{1,k}^2. 
\end{align*} 
Furthermore, $r^2 \leq 4K \, r_{1,k} \, d_{g(\bar{t})}(p_{1,t_k},\bar{x}) \leq \frac{4K}{\Lambda} \, d_{g(\bar{t})}(p_{1,t_k},\bar{x})^2$. Since $L \sqrt{\frac{4K}{\Lambda}} \leq 10^{-6}$, we obtain $r \leq \sqrt{\frac{4K}{\Lambda}} \, d_{g(\bar{t})}(p_{1,t_k},\bar{x}) \leq 10^{-6} \, L^{-1} \, d_{g(\bar{t})}(p_{1,t_k},\bar{x})$. Consequently, 
\begin{align*} 
d_{g(\bar{t})}(p_{1,t_k},x) 
&\geq d_{g(\bar{t})}(p_{1,t_k},\bar{x}) - Lr \\ 
&\geq (1-10^{-6}) \, d_{g(\bar{t})}(p_{1,t_k},\bar{x}) \\ 
&\geq (1-10^{-6}) \, 2^{\frac{j}{400}} \, \Lambda r_{1,k} \\ 
&\geq 2^{\frac{j-1}{400}} \, \Lambda r_{1,k} 
\end{align*}
for all $x \in B_{g(\bar{t})}(\bar{x},Lr)$. Moreover, using the inequalities $\frac{1}{2} \, r_{1,k} \leq r$ and $r^2 \leq 4K \, r_{1,k} \, d_{g(\bar{t})}(p_{1,t_k},\bar{x})$, we obtain 
\begin{align*} 
d_{g(\bar{t})}(p_{1,t_k},x) + 80Lr^2 r_{1,k}^{-1} 
&\leq d_{g(\bar{t})}(p_{1,t_k},\bar{x}) + Lr + 80Lr^2 r_{1,k}^{-1} \\ 
&\leq d_{g(\bar{t})}(p_{1,t_k},\bar{x}) + 82Lr^2 r_{1,k}^{-1} \\ 
&\leq 400KL \, d_{g(\bar{t})}(p_{1,t_k},\bar{x}) \\ 
&\leq (400KL)^{-j+1} \, \delta_k^{-1} \, r_{1,k} 
\end{align*} 
for all $x \in B_{g(\bar{t})}(\bar{x},Lr)$. Lemma \ref{derivative.of.distance.function} implies 
\[d_{g(\bar{t})}(p_{1,t_k},x) \leq d_{g(t)}(p_{1,t_k},x) \leq d_{g(\bar{t})}(p_{1,t_k},x) + 80Lr^2 r_{1,k}^{-1},\]
hence 
\[2^{\frac{j-1}{400}} \, \Lambda r_{1,k} \leq d_{g(t)}(p_{1,t_k},x) \leq (400KL)^{-j+1} \, \delta_k^{-1} \, r_{1,k}\] 
for all $(x,t) \in B_{g(\bar{t})}(\bar{x},Lr) \times [\bar{t}-Lr^2,\bar{t}]$. Therefore, the induction hypothesis guarantees that every point in $B_{g(\bar{t})}(\bar{x},Lr) \times [\bar{t}-Lr^2,\bar{t}]$ is $2^{-j} \varepsilon_k$-symmetric. By the Neck Improvement Theorem in \cite{Brendle2}, the point $(\bar{x},\bar{t})$ must be $2^{-j-1} \varepsilon_k$-symmetric. This completes the proof. \\

\begin{proposition}
\label{improvement.of.symmetry.of.cap}
If $k$ is sufficiently large, then the flow is $\frac{\varepsilon_k}{2}$-symmetric at time $t_k$. 
\end{proposition}

\textbf{Proof.} 
The arguments in Section 9 of \cite{Brendle2} go through unchanged. \\

Proposition \ref{improvement.of.symmetry.of.cap} contradicts the definition of $t_k$. This completes the proof of Theorem \ref{rotational.symmetry}.

\section{A priori estimates for compact ancient $\kappa$-solutions with rotational symmetry}

\label{a.priori.estimates}

We first recall some basic facts about the Bryant soliton.

\begin{proposition}[R.~Bryant \cite{Bryant}]
\label{asymptotics.of.bryant.soliton}
Consider the Bryant soliton, normalized so that the scalar curvature at the tip is equal to $1$. Then the metric can be written in the form $\Phi(r)^{-1} \, dr \otimes dr + r^2 \, g_{S^2}$, where $\Phi(r)$ satisfies the ODE 
\[\Phi(r) \Phi''(r) - \frac{1}{2} \, \Phi'(r)^2 + r^{-2} \, (1-\Phi(r)) \, (r \Phi'(r) + 2\Phi(r)) = 0.\] 
Moreover, $\Phi(r) = 1-\frac{r^2}{6} + O(r^4)$ as $r \to 0$ and $\Phi(r) = r^{-2} + 2r^{-4} + O(r^{-6})$ as $r \to \infty$. 
\end{proposition}

\textbf{Proof.} 
See \cite{Bryant}, Theorem 1 on p.~17. \\

\begin{proposition}
\label{estimate.for.normalized.Bryant.profile}
Let $\eta>0$ be given. If $|s|$ is sufficiently small (depending on $\eta$), then 
\[\big|\Phi((1+s)r)^{-1} - \Phi(r)^{-1}\big| \leq \eta \, \big(\Phi(r)^{-1} - 1\big)\] 
for all $r \geq 0$.
\end{proposition}

\textbf{Proof.}
We define $\chi(r) = r^{-2} \, (\Phi(r)^{-1} - 1)$. Note that $\chi(r)$ is a positive smooth function. Moreover, $\chi(r)$ satisfies the asymptotic expansions $\chi(r) = \frac{1}{6} + O(r^2)$ as $r \to 0$ and $\chi(r) = 1 + O(r^{-2})$ as $r \to \infty$. In particular, $\chi(r)$ is uniformly bounded above and below by positive constants. Hence, if $|s|$ is sufficiently small (depending on $\eta$), then 
\[|\chi((1+s)r) - \chi(r)| \leq \frac{\eta}{2} \, \chi(r)\] 
for all $r>0$. Therefore, if $|s|$ is sufficiently small (depending on $\eta$), then 
\begin{align*} 
|(1+s)^2 \, \chi((1+s)r) - \chi(r)| 
&\leq (1+s)^2 \, |\chi((1+s)r) - \chi(r)| + |(1+s)^2 -1| \, \chi(r) \\ 
&\leq \eta \, \chi(r) 
\end{align*} 
for all $r>0$. This gives 
\[\big | \Phi((1+s)r)^{-1} - \Phi(r)^{-1} \big | \leq \eta \, \big ( \Phi(r)^{-1} - 1 \big )\] 
for all $r>0$. \\

\begin{corollary} 
\label{concavity.of.F^2.on.bryant.soliton}
Consider the Bryant soliton, normalized so that the scalar curvature at the tip is equal to $1$. Let us write the metric in the form $dz \otimes dz + B(z)^2 \, g_{S^2}$. Then there exists a large constant $L_0$ such that $\frac{d^2}{dz^2} B(z)^2 < 0$ if $B(z)^2 \geq \frac{L_0^2}{4}$. 
\end{corollary} 

\textbf{Proof.} 
Since $r \Phi'(r) + 2\Phi(r) = -4r^{-4} + O(r^{-6})$ as $r \to \infty$, we conclude that $r \Phi'(r) + 2\Phi(r) < 0$ for $r$ sufficiently large. We next observe that $\big ( \frac{d}{dz} B(z) \big )^2 = \Phi(B(z))$. Differentiating this identity with respect to $z$ gives $2 \, \frac{d^2}{dz^2} B(z) = \Phi'(B(z))$. Thus, we conclude that $\frac{d^2}{dz^2} B(z)^2 = B(z) \, \Phi'(B(z)) + 2 \Phi(B(z)) < 0$ if $B(z)$ is sufficiently large. \\

\begin{corollary}
\label{asymptotics.of.bryant.soliton.2}
Consider the Bryant soliton, normalized so that the scalar curvature at the tip is equal to $1$. Let us write the metric in the form $dz \otimes dz + B(z)^2 \, g_{S^2}$. Then $B(z) \, \frac{d}{dz} B(z) \to 1$ as $z \to \infty$.
\end{corollary}

\textbf{Proof.} 
Note that $r \, \Phi(r)^{\frac{1}{2}} \to 1$ as $r \to \infty$. Using the identity $\big ( \frac{d}{dz} B(z) \big )^2 = \Phi(B(z))$, we obtain $B(z) \, \frac{d}{dz} B(z) = B(z) \, \Phi(B(z))^{\frac{1}{2}} \to 1$ as $z \to \infty$. \\
 
We now assume that $(S^3,g(t))$ is an ancient $\kappa$-solution which is not a family of shrinking round spheres. Let $q\in S^3$ be a reference point chosen as in \cite{Angenent-Brendle-Daskalopoulos-Sesum}. Recall that $q$ is chosen so that $\limsup_{t\to-\infty} (-t) R(q,t) \leq 100$ (see Proposition 3.1 in \cite{Angenent-Brendle-Daskalopoulos-Sesum}). In the same paper we showed that if $t_j \to -\infty$ and if we dilate the flow around the point $(q,t_j)$ by the factor $(-t_j)^{-\frac{1}{2}}$, then the rescaled manifolds converge to a cylinder of radius $\sqrt{2}$.  Let $F(z,t)$ denote  the radius of a sphere of symmetry in $(S^3, g(t))$ which has signed distance $z$ from point $q$. The function $F(z,t)$ satisfies the PDE 
\[F_t(z,t) - F_{zz}(z,t) = -F(z,t)^{-1} \, (1-F_z(z,t)^2) - 2 \, F_z(z,t) \int_0^z \frac{F_{zz}(z',t)}{F(z',t)} \, dz'.\] 
For abbreviation, let $H(z,t) := \frac{1}{2} \, F(z,t)^2 + t$. \\

\begin{lemma}
\label{boundary.points}
Let $L_0$ be the constant in Corollary \ref{concavity.of.F^2.on.bryant.soliton}. There exists a time $T_0 < 0$ with the following property. If $t \leq T_0$ and $F(z,t)^2 = L_0^2 \, \frac{(-t)}{\log (-t)}$, then $H_{zz}(z,t) < 0$.
\end{lemma}

\textbf{Proof.} 
Suppose this is false. Then there exists a sequence of times $t_j \to -\infty$ and a sequence of points $z_j$ such that $F(z_j,t_j)^2 = L_0^2 \, \frac{(-t_j)}{\log (-t_j)}$ and $H_{zz}(z_j,t_j) \geq 0$. By the result in \cite{Angenent-Brendle-Daskalopoulos-Sesum}, the curvature in the tip region behaves like $(2+o(1)) \, \frac{\log(-t_j)}{(-t_j)}$. Since $F(z_j,t_j)^2 = L_0^2 \, \frac{(-t_j)}{\log (-t_j)}$, it follows that the point $(z_j,t_j)$ has distance at most $C(L_0) \, \sqrt{\frac{(-t_j)}{\log(-t_j)}}$ from one of the tips. Hence, if we rescale around the point $(z_j,t_j)$, the rescaled metrics converge to the Bryant soliton. Passing to the limit, we find a point $z_\infty$ on the Bryant soliton such that $B(z_\infty)^2 = L_0^2$ and $\frac{d^2}{dz^2} B(z)^2 \big |_{z=z_\infty} \geq 0$. This contradicts Corollary \ref{concavity.of.F^2.on.bryant.soliton}. \\

\begin{lemma}
\label{evolution.of.H_zz}
The function $H_{zz}(z,t)$ satisfies the evolution equation 
\begin{align*} 
&H_{zzt}(z,t) - H_{zzzz}(z,t) \\ 
&\leq -2 \, H_{zzz}(z,t) \int_0^z \frac{F_{zz}(z',t)}{F(z',t)} \, dz' - 2 \, F(z,t)^{-1} \, F_z(z,t) \, H_{zzz}(z,t). 
\end{align*}
\end{lemma}

\textbf{Proof.} 
The function $H(z,t)$ satisfies 
\[H_t(z,t) - H_{zz}(z,t) = -2 \, H_z(z,t) \int_0^z \frac{F_{zz}(z',t)}{F(z',t)} \, dz'.\]
This implies  
\begin{align*} 
&H_{zzt}(z,t) - H_{zzzz}(z,t) \\ 
&= -2 \, H_{zzz}(z,t) \int_0^z \frac{F_{zz}(z',t)}{F(z',t)} \, dz' - 2 \, F(z,t)^{-1} \, F_z(z,t) \, H_{zzz}(z,t) \\ 
&- 4 \, F_{zz}(z,t)^2 + 4 \, F(z,t)^{-1} \, F_z(z,t)^2 \, F_{zz}(z,t).
\end{align*}
Since $F_{zz}(z,t) \leq 0$, the assertion follows. \\

\begin{proposition} 
\label{concavity.of.F2}
Let $L_0$ be chosen as in Corollary \ref{concavity.of.F^2.on.bryant.soliton} and let $T_0$ be chosen as in Lemma \ref{boundary.points}. If $t \leq T_0$ and $F(z,t)^2 \geq L_0^2 \, \frac{(-t)}{\log (-t)}$, then $H_{zz}(z,t) \leq 0$.
\end{proposition}

\textbf{Proof.} 
Suppose this is false. Then we can find a point $(z_0,t_0)$ such that $t_0 \leq T_0$, $F(z_0,t_0)^2 \geq L_0^2 \, \frac{(-t_0)}{\log (-t_0)}$, and $H_{zz}(z_0,t_0) > 0$. In view of Lemma \ref{boundary.points} and Lemma \ref{evolution.of.H_zz}, the maximum principle gives 
\[\sup_{F(z,t)^2 \geq L_0^2 \, \frac{(-t)}{\log (-t)}} H_{zz}(z,t) \geq H_{zz}(z_0,t_0) > 0\] 
for each $t \leq t_0$. Let us consider a sequence $t_j \to -\infty$. For $j$ large, we can find a point $z_j$ such that $F(z_j,t_j)^2 \geq L_0^2 \, \frac{(-t_j)}{\log (-t_j)}$ and $H_{zz}(z_j,t_j) \geq H_{zz}(z_0,t_0) > 0$. Using the inequality $F_{zz} \leq 0$, we obtain $F_z(z_j,t_j)^2 \geq H_{zz}(z_j,t_j) \geq H_{zz}(z_0,t_0) > 0$ for $j$ large. Hence, if we rescale around the points $(z_j,t_j)$ and pass to the limit, then the limit cannot be a cylinder. Consequently, the limit of these rescalings must be the Bryant soliton. Hence, after passing to the limit, we obtain a point $z_\infty$ on the Bryant soliton such that $B(z_\infty)^2 \geq L_0^2$ and $\frac{d^2}{dz^2} B(z)^2 \big |_{z=z_\infty} \geq 0$. This contradicts Corollary \ref{concavity.of.F^2.on.bryant.soliton}. \\

We next recall a crucial estimate from \cite{Angenent-Brendle-Daskalopoulos-Sesum}. \\

\begin{proposition}[cf. \cite{Angenent-Brendle-Daskalopoulos-Sesum}]
\label{precise.estimate.for.F}
Fix a small number $\theta > 0$ and a small number $\eta > 0$. Then 
\[\Big | \frac{1}{2} \, F(z,t)^2 + t + \frac{z^2+2t}{4 \log(-t)} \Big | \leq \eta \, \frac{z^2-t}{\log(-t)}\] 
if $F(z,t) \geq \frac{\theta}{400} \sqrt{-t}$ and $-t$ is sufficiently large (depending on $\eta$ and $\theta$). 
\end{proposition}

\textbf{Proof.} 
By Proposition 6.3 and Proposition 6.4 in \cite{Angenent-Brendle-Daskalopoulos-Sesum}, we can find a large number $M$ (depending on $\eta$ and $\theta$) with the property that 
\[\Big | \frac{1}{2} \, F(z,t)^2 + t + \frac{z^2+2t}{4 \log(-t)} \Big | \leq \eta \, \frac{z^2}{\log(-t)}\] 
whenever $|z| \geq M \sqrt{-t}$, $F(z,t) \geq \frac{\theta}{400} \sqrt{-t}$, and $-t$ is sufficiently large. Having fixed $M$, Propositon 5.10 in \cite{Angenent-Brendle-Daskalopoulos-Sesum} implies that 
\[\Big | \frac{1}{2} \, F(z,t)^2 + t + \frac{z^2+2t}{4  \log(-t)} \Big | \leq \eta \, \frac{(-t)}{\log(-t)}\] 
whenever $|z| \leq M \sqrt{-t}$ and $-t$ is sufficiently large. Putting these facts together, we conclude that 
\[\Big | \frac{1}{2} \, F(z,t)^2 + t + \frac{z^2+2t}{4 \log(-t)} \Big | \leq \eta \, \frac{z^2-t}{\log(-t)}\] 
whenever $F(z,t) \geq \frac{\theta}{400} \sqrt{-t}$ and $-t$ is sufficiently large. \\

\begin{proposition}
\label{precise.estimate.for.F_z}
Let us fix a small number $\theta > 0$ and a small number $\eta > 0$. Then 
\[\Big | F(z,t) \, F_z(z,t) + \frac{z}{2 \log(-t)} \Big | \leq \eta \, \frac{|z|+\sqrt{-t}}{\log(-t)}\] 
if $F(z,t) \geq \frac{\theta}{200} \sqrt{-t}$ and $-t$ is sufficiently large (depending on $\eta$ and $\theta$).
\end{proposition}

\textbf{Proof.}
Let $\theta \in (0,\frac{1}{2})$ and $\eta \in (0,\frac{1}{2})$ be given. We can find a small positive number $\mu \in (0,\eta)$ and time $T_0$ with the property that $F((1+\mu)z,t) \geq \frac{\theta}{400} \sqrt{-t}$ whenever $F(z,t) \geq \frac{\theta}{200} \sqrt{-t}$ and $t \leq T_0$. Moreover, by Proposition \ref{precise.estimate.for.F}, we can find a time $T \leq T_0$ such that 
\[\Big | \frac{1}{2} \, F(z,t)^2 + t + \frac{z^2+2t}{4 \log(-t)} \Big | \leq \eta \, \mu \, \frac{z^2}{16 \log(-t)}\] 
whenever $z \geq 2\sqrt{-t_0}$, $F(z,t) \geq \frac{\theta}{400} \sqrt{-t}$, and $t \leq T$. 

Suppose now that $(z_0,t_0)$ is a point in spacetime satisfying $z_0 \geq 4\sqrt{-t_0}$, $F(z_0,t_0) \geq \frac{\theta}{200} \sqrt{-t_0}$, and $t_0 \leq T$. Then $F(z,t_0) \geq \frac{\theta}{400} \sqrt{-t_0}$ for all $z \in [(1-\mu)z_0,(1+\mu)z_0]$. Consequently, 
\[\Big | \frac{1}{2} \, F(z,t_0)^2 + t_0 + \frac{z^2+2t_0}{4 \log(-t_0)} \Big | \leq \eta \, \mu \, \frac{z_0^2}{4 \log(-t_0)}\] 
for all $z \in [(1-\mu)z_0,(1+\mu)z_0]$. This implies 
\[\inf_{z \in [(1-\mu)z_0,z_0]} \Big ( F(z,t_0) \, F_z(z,t_0) + \frac{z}{2 \log(-t_0)} \Big ) \leq \eta \, \frac{z_0}{2 \log(-t_0)}\] 
and 
\[\sup_{z \in [z_0,(1+\mu)z_0]} \Big ( F(z,t_0) \, F_z(z,t_0) + \frac{z}{2 \log(-t_0)} \Big ) \geq -\eta \, \frac{z_0}{2 \log(-t_0)}\] 
By Proposition \ref{concavity.of.F2}, the function $z \mapsto FF_z$ is monotone decreasing in the relevant region. This gives 
\[F(z_0,t_0) \, F_z(z_0,t_0) + \frac{(1-\mu)z_0}{2 \log(-t_0)} \leq \eta \, \frac{z_0}{2 \log(-t_0)}\] 
and 
\[F(z_0,t_0) \, F_z(z_0,t_0) + \frac{(1+\mu)z_0}{2 \log(-t_0)} \geq -\eta \, \frac{z_0}{2 \log(-t_0)}.\] 
Since $\mu \in (0,\eta)$, it follows that 
\[\Big | F(z_0,t_0) \, F_z(z_0,t_0) + \frac{z_0}{2 \log(-t_0)} \Big | \leq (\eta+\mu) \, \frac{z_0}{2 \log(-t_0)} \leq \eta \, \frac{z_0}{\log(-t_0)}.\] 
To summarize, we have verified the assertion for $z \geq 4\sqrt{-t}$. An analogous argument shows that the assertion holds for $z \leq -4\sqrt{-t}$. Finally, if $|z| \leq 4\sqrt{-t}$, then the assertion follows from Proposition 5.10 in \cite{Angenent-Brendle-Daskalopoulos-Sesum}. This completes the proof of Proposition \ref{precise.estimate.for.F_z}. \\

\begin{corollary}
\label{bound.for.F_z} 
Let us fix a small number $\theta > 0$. Then 
\[|F_z(z,t)| \leq \frac{C(\theta)}{\sqrt{\log(-t)}}\] 
if $F(z,t) \geq \frac{\theta}{200} \sqrt{-t}$ and $-t$ is sufficiently large (depending on $\theta$).
\end{corollary}

\textbf{Proof.} 
The asymptotic estimates in \cite{Angenent-Brendle-Daskalopoulos-Sesum} imply that $|z| \leq (2+o(1)) \sqrt{(-t) \log(-t)}$. Hence, the assertion follows from Proposition \ref{precise.estimate.for.F_z}. \\

\begin{proposition} 
\label{higher.derivative.bounds.in.cylindrical.region}
Let us fix a small number $\theta > 0$. Then 
\[F(z,t) \, |F_{zz}(z,t)| + F(z,t)^2 \, |F_{zzz}(z,t)| \leq \frac{C(\theta)}{\sqrt{\log (-t)}}\] 
if $F(z,t) \geq \frac{\theta}{100} \sqrt{-t}$ and $-t$ is sufficiently large (depending on $\theta$).
\end{proposition}

\textbf{Proof.} 
Let us fix a small number $\varepsilon > 0$. Moreover, we consider a point $(p_0,t_0)$ in space-time with the property that the sphere of symmetry passing through $(p_0,t_0)$ has radius $r_0 \geq \frac{\theta}{100} \sqrt{-t_0}$ at $(p_0,t_0)$. If $-t_0$ is sufficiently large (depending on $\theta$ and $\varepsilon$), then the point $(p_0,t_0)$ lies at the center of an evolving $\varepsilon$-neck. Let $\tilde{F}(z,t)$ denote the radius of the sphere of symmetry which has signed distance $z$ from the point $p_0$. By assumption, $\tilde{F}(0,t_0) = r_0 \geq \frac{\theta}{100} \sqrt{-t_0}$. Since the point $(p_0,t_0)$ lies on a neck, we have $\frac{1}{2} \, r_0 \leq \tilde{F}(z,t) \leq 100 \, r_0$ and $|\tilde{F}_z(z,t)| + r_0 \, |\tilde{F}_{zz}(z,t)| + r_0^2 \, |\tilde{F}_{zzz}(z,t)| \leq 1$ for all $(z,t) \in [-r_0,r_0] \times [t_0-r_0^2,t_0]$. Moreover, since $(p_0,t_0)$ lies on a a neck, we obtain 
\[\tilde{F}(z,t)^2 \geq \sqrt{\frac{1}{4} \, \tilde{F}(0,t_0)^2 + (t_0-t)} \geq \sqrt{\Big ( \frac{\theta}{200} \Big )^2 (-t_0) + (t_0-t)} \geq \frac{\theta}{200} \sqrt{-t}\] 
for all $(z,t) \in [-r_0,r_0] \times [t_0-r_0^2,t_0]$. Hence, Corollary \ref{bound.for.F_z} implies $|\tilde{F}_z(z,t)| \leq \frac{C(\theta)}{\sqrt{\log (-t)}}$ for all $(z,t) \in [-r_0,r_0] \times [t_0-r_0^2,t_0]$. 

The function $\tilde{F}$ satisfies the same PDE as the original function $F$. In other words, 
\begin{align*} 
\tilde{F}_t(z,t) - \tilde{F}_{zz}(z,t) 
&= -\tilde{F}(z,t)^{-1} \, (1+\tilde{F}_z(z,t)^2) \\ 
&+ 2 \, \tilde{F}_z(z,t) \, \bigg [ \tilde{F}(0,t)^{-1} \, \tilde{F}_z(0,t) - \int_0^z \frac{\tilde{F}_z(z',t)^2}{\tilde{F}(z',t)^2} \, dz' \bigg ]. 
\end{align*} 
Differentiating this equation with respect to $z$ gives 
\begin{align*} 
\tilde{F}_{zt}(z,t) - \tilde{F}_{zzz}(z,t) 
&= \tilde{F}(z,t)^{-2} \, \tilde{F}_z(z,t) \, (1-\tilde{F}_z(z,t)^2) \\
&- 2 \, \tilde{F}(z,t)^{-1} \, \tilde{F}_z(z,t) \, \tilde{F}_{zz}(z,t) \\ 
&+ 2 \, \tilde{F}_{zz}(z,t) \, \bigg [ \tilde{F}(0,t)^{-1} \, \tilde{F}_z(0,t) - \int_0^z \frac{\tilde{F}_z(z',t)^2}{\tilde{F}(z',t)^2} \, dz' \bigg ]
\end{align*} 
and 
\begin{align*} 
\tilde{F}_{zzt}(z,t) - \tilde{F}_{zzzz}(z,t) 
&= -2 \, \tilde{F}(z,t)^{-3} \, \tilde{F}_z(z,t)^2 \, (1-\tilde{F}_z(z,t)^2) \\ 
&+ \tilde{F}(z,t)^{-2} \, \tilde{F}_{zz}(z,t) \, (1-3 \, \tilde{F}_z(z,t)^2) \\
&- 2 \, \tilde{F}(z,t)^{-1} \, \tilde{F}_{zz}(z,t)^2 - 2 \, \tilde{F}(z,t)^{-1} \, \tilde{F}_z(z,t) \, \tilde{F}_{zzz}(z,t) \\ 
&+ 2 \, \tilde{F}_{zzz}(z,t) \, \bigg [ \tilde{F}(0,t)^{-1} \, \tilde{F}_z(0,t) - \int_0^z \frac{\tilde{F}_z(z',t)^2}{\tilde{F}(z',t)^2} \, dz' \bigg ]. 
\end{align*} 
For $(z,t) \in [-r_0,r_0] \times [t_0-r_0^2,t_0]$, we have $|\tilde{F}_z(z,t)| \leq \frac{C(\theta)}{\sqrt{\log (-t_0)}}$, hence $r_0^2 \, |\tilde{F}_{zt}(z,t) - \tilde{F}_{zzz}(z,t)| \leq \frac{C(\theta)}{\sqrt{\log (-t_0)}}$. Using standard interior estimates for parabolic equations, we obtain $r_0 \, |\tilde{F}_{zz}(z,t)| \leq \frac{C(\theta)}{\sqrt{\log (-t_0)}}$ for all $(z,t) \in [-\frac{r_0}{2},\frac{r_0}{2}] \times [t_0-\frac{r_0^2}{4},t_0]$. This implies $r_0^3 \, |\tilde{F}_{zzt}(z,t) - \tilde{F}_{zzzz}(z,t)| \leq \frac{C(\theta)}{\sqrt{\log (-t_0)}}$ for all $(z,t) \in [-\frac{r_0}{2},\frac{r_0}{2}] \times [t_0-\frac{r_0^2}{4},t_0]$. Hence, standard interior estimates for parabolic equations give $r_0^2 \, |\tilde{F}_{zzz}(z,t)| \leq \frac{C(\theta)}{\sqrt{\log (-t_0)}}$ for all $(z,t) \in [-\frac{r_0}{4},\frac{r_0}{4}] \times [t_0-\frac{r_0^2}{16},t_0]$. 

Thus, $r_0 \, |\tilde{F}_{zz}(0,t_0)| + r_0^2 \, |\tilde{F}_{zzz}(0,t_0)| \leq \frac{C(\theta)}{\sqrt{\log (-t_0)}}$. This finally implies $F \, |F_{zz}| + F^2 \, |F_{zzz}| \leq \frac{C(\theta)}{\sqrt{\log (-t_0)}}$ at the point $(p_0,t_0)$. This completes the proof of Proposition \ref{higher.derivative.bounds.in.cylindrical.region}. \\

\begin{proposition}
\label{derivative.of.F.wrt.t}
Let us fix a small number $\theta > 0$. Then 
\[|1 + FF_t| \leq \frac{C(\theta)}{\sqrt{\log(-t)}}\] 
whenever $F \geq \frac{\theta}{100} \sqrt{-t}$, and $-t$ is sufficiently large (depending on $\theta$).
\end{proposition}

\textbf{Proof.} 
Using the evolution equation for $F$, we obtain 
\begin{align*} 
&1 + F(z,t) \, F_t(z,t) \\ 
&= F(z,t) \, F_{zz}(z,t) - F_z(z,t)^2 \\ 
&+ 2 \, F(z,t) \, F_z(z,t) \, \bigg [ F(0,t)^{-1} \, F_z(0,t) - \int_0^z \frac{F_z(z',t)^2}{F(z',t)^2} \, dz' \bigg ]. 
\end{align*} 
The asymptotic estimates in \cite{Angenent-Brendle-Daskalopoulos-Sesum} imply that, for $-t$ sufficiently large, the domain of definition of the function $z \mapsto F(z,t)$ is contained in the interval $[-4 \sqrt{(-t) \log(-t)},4\sqrt{(-t) \log(-t)}]$. Moreover, if $F(z,t) \geq \frac{\theta}{100} \sqrt{-t}$, then $F(z',t) \geq \frac{\theta}{100} \sqrt{-t}$ for all $z'$ between $0$ and $z$. Using Corollary \ref{bound.for.F_z}, we obtain 
\[\bigg | \int_0^z \frac{F_z(z',t)^2}{F(z',t)^2} \, dz' \bigg | \leq \frac{C(\theta)}{(-t) \log(-t)} \, |z| \leq \frac{C(\theta)}{\sqrt{(-t) \log(-t)}}\] 
whenever $F(z,t) \geq \frac{\theta}{100} \sqrt{-t}$, and $-t$ is sufficiently large. Using Corollary \ref{bound.for.F_z} and Proposition \ref{higher.derivative.bounds.in.cylindrical.region}, we conclude that 
\[|1 + F(z,t) \, F_t(z,t)| \leq \frac{C(\theta)}{\sqrt{\log(-t)}}\]  
whenever $F \geq \frac{\theta}{100} \sqrt{-t}$, and $-t$ is sufficiently large. This completes the proof. \\

\begin{proposition}
\label{neck}
Let $\varepsilon>0$ be given. Then there exists a large number $L$ (depending on $\varepsilon$) and a time $T$ such that the following holds. If $F \geq L \, \sqrt{\frac{(-t)}{\log(-t)}}$ and $t \leq T$ at some point in space-time, then that point lies at the center of an evolving $\varepsilon$-neck. 
\end{proposition}

\textbf{Proof.} This follows from the fact, established in \cite{Angenent-Brendle-Daskalopoulos-Sesum}, that the scalar curvature at each tip is $(1+o(1)) \, \frac{\log(-t)}{(-t)}$. \\

\begin{corollary}
\label{higher.derivative.bounds.in.collar.region}
Let $\eta>0$ be given. Then there exists a large number $L$ (depending on $\eta$) and a time $T$ such that $|F_z| + F \, |F_{zz}| + F^2 \, |F_{zzz}| \leq \eta$ whenever $F \geq L \, \sqrt{\frac{(-t)}{\log(-t)}}$ and $t \leq T$.
\end{corollary}

\textbf{Proof.} 
This follows directly from Proposition \ref{neck}. \\

\begin{proposition}
\label{estimate.for.F_z.in.collar.region}
Let $\eta>0$ be given. Then there exist a large number $L \in (\eta^{-1},\infty)$ and a small number $\theta \in (0,\eta)$ (depending on $\eta$), and a time $T$ with the property that 
\[\bigg | 1 - \sqrt{\frac{\log(-t)}{(-t)}} \, F \, |F_z| \bigg | \leq \eta\] 
whenever $L \, \sqrt{\frac{(-t)}{\log(-t)}} \leq F \leq 100\theta \sqrt{-2t}$ and $t \leq T$.
\end{proposition}

\textbf{Proof.} 
By Corollary \ref{asymptotics.of.bryant.soliton.2}, we can find a large number $L \in (\eta^{-1},\infty)$ such that $\big | 1 - B(z) \, \frac{d}{dz} B(z) \big | \leq \frac{\eta}{2}$ for $z \geq \frac{L}{2}$. Recall that the solution looks like the Bryant soliton near each tip, and the scalar curvature at each tip equals $(1+o(1)) \, \frac{\log(-t)}{(-t)}$. Consequently, 
\[\bigg | 1 - \sqrt{\frac{\log(-t)}{(-t)}} \, F \, |F_z| \bigg | \leq \eta\] 
if $F = L \, \sqrt{\frac{(-t)}{\log(-t)}}$ and $-t$ is sufficiently large. On the other hand, for each $\theta \in (0,\frac{1}{1000})$, Proposition \ref{precise.estimate.for.F} implies
\[z^2 = (4+o(1)) \, (1-(100\theta)^2) \, (-t) \log(-t)\] 
if $F = 100\theta \sqrt{-2t}$. Using Proposition \ref{precise.estimate.for.F_z}, we obtain 
\[F \, |F_z| = (1+o(1)) \, \frac{|z|}{2 \log(-t)}\] 
if $F = 100\theta \sqrt{-2t}$. Consequently, 
\[1 - \sqrt{\frac{\log(-t)}{(-t)}} \, F \, |F_z| = 1-\sqrt{1-(100\theta)^2}+o(1)\]
if $F = 100\theta \sqrt{-2t}$. Therefore, if we choose $\theta$ sufficiently small (depending on $\eta$), then we obtain 
\[\bigg | 1 - \sqrt{\frac{\log(-t)}{(-t)}} \, F \, |F_z| \bigg | \leq \eta\]
if $F = 100\theta \sqrt{-2t}$ and $-t$ is sufficiently large. Hence, the assertion follows from the fact that the function $z \mapsto FF_z$ is monotone decreasing in the relevant region (see Proposition \ref{concavity.of.F2}). This completes the proof of Proposition \ref{estimate.for.F_z.in.collar.region}. \\

In the remainder of this section, we define functions $U_+(r,t)$ and $U_-(r,t)$ so that 
\[U_+(r,t) = \Big ( \frac{\partial}{\partial z} F(z,t) \Big )^2\] 
for $r = F(z,t)$ and $z \geq 2\sqrt{-t}$ and 
\[U_-(r,t) = \Big ( \frac{\partial}{\partial z} F(z,t) \Big )^2\] 
for $r = F(z,t)$ and $z \leq -2\sqrt{-t}$. Let us consider the rescaled functions 
\begin{align*}
V_+(\rho,\tau) &:= \sqrt{U_+(e^{-\frac{\tau}{2}} \rho,-e^{-\tau})}, \\ 
V_-(\rho,\tau) &:= \sqrt{U_-(e^{-\frac{\tau}{2}} \rho,-e^{-\tau})}. 
\end{align*}
For each $\rho \in (0,1)$, we denote by $\xi_+(\rho,\tau)$ the unique positive solution of the equation $F(e^{-\frac{\tau}{2}} \xi,-e^{-\tau}) = e^{-\frac{\tau}{2}} \rho$; moreover, we denote by $\xi_-(\rho,\tau)$ the unique negative solution of the equation $F(e^{-\frac{\tau}{2}} \xi,-e^{-\tau}) = e^{-\frac{\tau}{2}} \rho$. \\

\begin{proposition}
\label{bound.for.V.in.transition.region}
Let us fix a small number $\theta > 0$. If $-\tau$ is sufficiently large (depending on $\theta$), then $\frac{1}{C(\theta)} \, (-\tau)^{-\frac{1}{2}} \leq V_+(\rho,\tau) \leq C(\theta) \, (-\tau)^{-\frac{1}{2}}$ and $\big | \frac{\partial}{\partial \rho} V_+(\rho,\tau) \big | \leq C(\theta)$ for every $\rho \in [\frac{\theta}{100},100\theta]$.
\end{proposition}

\textbf{Proof.} 
Proposition \ref{precise.estimate.for.F_z} implies that $\frac{1}{C(\theta) \, \sqrt{\log(-t)}} \leq |F_z(z,t)| \leq \frac{C(\theta)}{\sqrt{\log(-t)}}$ whenever $\frac{\theta}{100} \, \sqrt{-t} \leq F(z,t) \leq 100\theta \, \sqrt{-t}$. Moreover, Proposition \ref{higher.derivative.bounds.in.cylindrical.region} gives $|F_{zz}(z,t)| \leq \frac{C(\theta)}{\sqrt{(-t) \log(-t)}}$ whenever $\frac{\theta}{100} \, \sqrt{-t} \leq F(z,t) \leq 100\theta \, \sqrt{-t}$. From this, the assertion follows easily. \\

\begin{proposition}
\label{comparison.of.V.with.bryant.soliton.profile}
Fix a small number $\eta>0$. Then we can find a small number $\theta \in (0,\eta)$ (depending on $\eta$) such that, for $-\tau$ sufficiently large, we have 
\[|V_+(\rho,\tau)^{-2} - \Phi((-\tau)^{\frac{1}{2}} \rho)^{-1}| \leq \eta \, (V_+(\rho,\tau)^{-2}-1)\] 
in the region $\{\rho \leq 100\theta\}$. Here, $\Phi$ denotes the profile of the Bryant soliton.
\end{proposition}

\textbf{Proof.} 
By Proposition \ref{estimate.for.F_z.in.collar.region}, we can find a large number $L \in (\eta^{-1},\infty)$ and a small number $\theta \in (0,\eta)$ with the property that 
\[\bigg | 1 - \sqrt{\frac{\log(-t)}{(-t)}} \, F \, |F_z| \bigg | \leq \frac{\eta}{32}\] 
whenever $L \, \sqrt{\frac{(-t)}{\log(-t)}} \leq F \leq 100\theta \sqrt{-t}$ and $-t$ is sufficiently large. This implies 
\[\big | 1 - (-\tau)^{\frac{1}{2}} \, \rho \, V_+(\rho,\tau) \big | \leq \frac{\eta}{32}\]
whenever $L \, (-\tau)^{-\frac{1}{2}} \leq \rho \leq 100\theta$ and $-\tau$ is sufficiently large. On the other hand, we can find a number $\tilde{L} \geq L$ such that 
\[\big | (-\tau)^{\frac{1}{2}} \, \rho \, \Phi((-\tau)^{\frac{1}{2}} \rho)^{\frac{1}{2}} - 1 \big | \leq \frac{\eta}{32}\] 
whenever $\rho \geq \tilde{L} \, (-\tau)^{-\frac{1}{2}}$ and $-\tau$ is sufficiently large. Putting these facts together, we conclude that 
\[\big | \Phi((-\tau)^{\frac{1}{2}} \rho)^{\frac{1}{2}} - V_+(\rho,\tau) \big | \leq \frac{\eta}{16} \, (-\tau)^{-\frac{1}{2}} \, \rho^{-1} \leq \frac{\eta}{8} \, V_+(\rho,\tau)\] 
whenever $\tilde{L} \, (-\tau)^{-\frac{1}{2}} \leq \rho \leq 100\theta$ and $-\tau$ is sufficiently large. This gives 
\[\big | V_+(\rho,\tau)^{-2} - \Phi((-\tau)^{\frac{1}{2}} \rho)^{-1} \big | \leq \frac{\eta}{2} \, V_+(\rho,\tau)^{-2} \leq \eta \, (V_+(\rho,\tau)^{-2}-1)\] 
whenever $\tilde{L} \, (-\tau)^{-\frac{1}{2}} \leq \rho \leq 100\theta$ and $-\tau$ is sufficiently large. 

On the other hand, since the solution looks like the Bryant soliton near each tip, we know that 
\[\big | V_+(\rho,\tau)^{-2} - \Phi((-\tau)^{\frac{1}{2}} \rho)^{-1} \big | \leq \eta \, (V_+(\rho,\tau)^{-2}-1)\] 
whenever $\rho \leq \tilde{L} \, (-\tau)^{-\frac{1}{2}}$ and $-\tau$ is sufficiently large. Putting these facts together, the assertion follows. \\

\begin{proposition}
\label{higher.derivatives.of.V.in.collar.region}
Fix a small number $\eta>0$. Then we can find a large number $L$ (depending on $\eta$) such that, for $-\tau$ sufficiently large, we have 
\[V_+(\rho,\tau) \leq \eta,\] 
\[\Big | \frac{\partial}{\partial \rho} V_+(\rho,\tau) \Big | \leq \eta \, \rho^{-1} \, V_+(\rho,\tau)^{-1}\] 
and 
\[\Big | \frac{\partial^2}{\partial \rho^2} V_+(\rho,\tau) \Big | \leq \eta \, \rho^{-2} \, V_+(\rho,\tau)^{-3}\] 
in the region $\{L \, (-\tau)^{-\frac{1}{2}} \leq \rho \leq \frac{1}{4}\}$.
\end{proposition}

\textbf{Proof.} 
By Corollary \ref{higher.derivative.bounds.in.collar.region}, we can find a large number $L$ (depending on $\eta$) such that $|F_z| + F \, |F_{zz}| + F^2 \, |F_{zzz}| \leq \eta$ whenever $F \geq L \, \sqrt{\frac{(-t)}{\log(-t)}}$. This implies 
\[V_+(\rho,\tau) \leq \eta,\] 
\[\Big | \rho \, V_+(\rho,\tau) \, \frac{\partial}{\partial \rho} V_+(\rho,\tau) \Big | \leq \eta,\] 
and 
\[\Big | \rho^2 \, V_+(\rho,\tau)^2 \, \frac{\partial^2}{\partial \rho^2} V_+(\rho,\tau) + \rho^2 \, V_+(\rho,\tau) \, \Big ( \frac{\partial}{\partial \rho} V_+(\rho,\tau) \Big )^2 \Big | \leq \eta\] 
in the region $\{L \, (-\tau)^{-\frac{1}{2}} \leq \rho \leq \frac{1}{4}\}$. From this, the assertion follows. \\

\begin{corollary}
\label{derivative.of.V.wrt.tau}
Fix a small number $\eta>0$. Then, for $-\tau$ sufficiently large, we have 
\[\Big | \frac{\partial}{\partial \tau} V_+(\rho,\tau) \Big | \leq \eta \, \rho^{-2} \, (V_+(\rho,\tau)^{-1}-1)\] 
in the region $\{\rho \leq \frac{1}{4}\}$. 
\end{corollary}

\textbf{Proof.} 
By Proposition \ref{higher.derivatives.of.V.in.collar.region}, we can find a large number $L$ (depending on $\eta$) such that, for $-\tau$ sufficiently large, we have 
\[V_+^2 \, \Big | \frac{\partial^2 V_+}{\partial \rho^2} + \rho^{-2} \, (V_+^{-2}-1) \, \Big ( \rho \, \frac{\partial V_+}{\partial \rho}+V_+ \Big ) \Big | \leq \frac{\eta}{2} \, \rho^{-2} \, (V_+^{-1}-1)\] 
and 
\[\Big | \rho \, \frac{\partial V_+}{\partial \rho} \Big | \leq \eta \, (V_+^{-1}-1)\] 
in the region $\{L \, (-\tau)^{-\frac{1}{2}} \leq \rho \leq \frac{1}{4}\}$. On the other hand, since the solution looks like the Bryant soliton near each tip, we know that 
\[V_+^2 \, \Big | \frac{\partial^2 V_+}{\partial \rho^2} + \rho^{-2} \, (V_+^{-2}-1) \, \Big ( \rho \, \frac{\partial V_+}{\partial \rho}+V_+ \Big ) \Big | \leq \frac{\eta}{2} \, \rho^{-2} \, (V_+^{-1}-1)\] 
and 
\[\Big | \rho \, \frac{\partial V_+}{\partial \rho} \Big | \leq C \, (V_+^{-1}-1)\] 
whenever $\rho \leq L \, (-\tau)^{-\frac{1}{2}}$ and $-\tau$ is sufficiently large. Putting these facts together, we conclude that 
\[V_+^2 \, \Big | \frac{\partial^2 V_+}{\partial \rho^2} + \rho^{-2} \, (V_+^{-2}-1) \, \Big ( \rho \, \frac{\partial V_+}{\partial \rho}+V_+ \Big ) \Big | \leq \frac{\eta}{2} \, \rho^{-2} \, (V_+^{-1}-1)\] 
and 
\[\Big | \rho \, \frac{\partial V_+}{\partial \rho} \Big | \leq \eta \, \rho^{-2} \, (V_+^{-1}-1)\] 
whenever $\rho \leq \frac{1}{4}$ and $-\tau$ is sufficiently large. Using the equation 
\[\frac{\partial V_+}{\partial \tau} + \frac{\rho}{2} \, \frac{\partial V_+}{\partial \rho} = V_+^2 \, \Big ( \frac{\partial^2 V_+}{\partial \rho^2} + \rho^{-2} \, (V_+^{-2}-1) \, \Big ( \rho \, \frac{\partial V_+}{\partial \rho}+V_+ \Big ) \Big ),\] 
we conclude that 
\[\Big | \frac{\partial V_+}{\partial \tau} \Big | \leq \eta \, \rho^{-2} \, (V_+^{-1}-1)\] 
ihenever $\rho \leq \frac{1}{4}$ and $-\tau$ is sufficiently large. This completes the proof of Corollary \ref{derivative.of.V.wrt.tau}. \\

\begin{proposition}
\label{derivative.of.xi.wrt.rho}
Fix a small number $\eta>0$. Then we can find a small number $\theta \in (0,\eta)$ (depending on $\eta$) such that, for $-\tau$ sufficiently large, we have 
\[\Big | \frac{\partial}{\partial \rho} \Big ( \frac{\xi_+(\rho,\tau)^2}{4} \Big ) + \rho^{-1} \, (V_+(\rho,\tau)^{-2}-1) \Big | \leq \eta \, \rho^{-1} \, (V_+(\rho,\tau)^{-2}-1)\] 
in the region $\{\frac{\theta}{8} \leq \rho \leq 2\theta\}$. 
\end{proposition}

\textbf{Proof.} 
In the following, $\theta>0$ will denote a small positive number which will be specified later. Proposition \ref{precise.estimate.for.F} and Proposition \ref{precise.estimate.for.F_z} imply 
\[\xi_+(\rho,\tau)^2 = (2+o(1)) \, (2-\rho^2) \, (-\tau)\] 
and 
\[\rho \, F_z(e^{-\frac{\tau}{2}} \, \xi_+(\rho,\tau),-e^{-\tau}) = -(1+o(1)) \, \frac{\xi_+(\rho,\tau)}{(-2\tau)}\]
in the region $\{\frac{\theta}{8} \leq \rho \leq 2\theta\}$. On the other hand, differentiating the relation $F(e^{-\frac{\tau}{2}} \, \xi_+(\rho,\tau),-e^{-\tau}) = e^{-\frac{\tau}{2}} \rho$ with respect to $\rho$ gives 
\[\frac{\partial}{\partial \rho} \xi_+(\rho,\tau) = F_z(e^{-\frac{\tau}{2}} \, \xi_+(\rho,\tau),-e^{-\tau})^{-1}.\] 
This implies 
\begin{align*} 
\frac{\partial}{\partial \rho} \Big ( \frac{\xi_+(\rho,\tau)^2}{4} \Big ) 
&= \frac{1}{2} \, \xi_+(\rho,\tau) \, F_z(e^{-\frac{\tau}{2}} \, \xi_+(\rho,\tau),-e^{-\tau})^{-1} \\ 
&= -(1+o(1)) \, \Big ( 1-\frac{\rho^2}{2} \Big ) \, \rho^{-1} \, F_z(e^{-\frac{\tau}{2}} \, \xi_+(\rho,\tau),-e^{-\tau})^{-2} \\ 
&= -(1+o(1)) \, \Big ( 1-\frac{\rho^2}{2} \Big ) \, \rho^{-1} \, V_+(\rho,\tau)^{-2} 
\end{align*}
in the region $\{\frac{\theta}{8} \leq \rho \leq 2\theta\}$. Hence, if we choose $\theta$ sufficiently small (depending on $\eta$), then 
\[\Big | \frac{\partial}{\partial \rho} \Big ( \frac{\xi_+(\rho,\tau)^2}{4} \Big ) + \rho^{-1} \, (V_+(\rho,\tau)^{-2}-1) \Big | \leq \eta \, \rho^{-1} \, (V_+(\rho,\tau)^{-2}-1)\] 
in the region $\{\frac{\theta}{8} \leq \rho \leq 2\theta\}$. This completes the proof of Proposition \ref{derivative.of.xi.wrt.rho}. \\

\begin{proposition}
\label{second.derivative.of.xi.wrt.rho}
Fix a small number $\theta>0$. Then, for $-\tau$ large, we have 
\[\Big | \frac{\partial}{\partial \rho} \Big ( \frac{\xi_+(\rho,\tau)^2}{4} \Big ) \Big | \leq C(\theta) \, (-\tau)\] 
and 
\[\Big | \frac{\partial^2}{\partial \rho^2} \Big ( \frac{\xi_+(\rho,\tau)^2}{4} \Big ) \Big | \leq C(\theta) \, (-\tau)^{\frac{3}{2}}\] 
in the region $\{\frac{\theta}{8} \leq \rho \leq 2\theta\}$. 
\end{proposition}

\textbf{Proof.}
Differentiating the relation $F(e^{-\frac{\tau}{2}} \, \xi_+(\rho,\tau),-e^{-\tau}) = e^{-\frac{\tau}{2}} \rho$ with respect to $\rho$ gives 
\[\frac{\partial}{\partial \rho} \xi_+(\rho,\tau) = F_z(e^{-\frac{\tau}{2}} \, \xi_+(\rho,\tau),-e^{-\tau})^{-1}\] 
and 
\[\frac{\partial^2}{\partial \rho^2} \xi_+(\rho,\tau) = -e^{-\frac{\tau}{2}} \, F_z(e^{-\frac{\tau}{2}} \, \xi_+(\rho,\tau),-e^{-\tau})^{-3} \, F_{zz}(e^{-\frac{\tau}{2}} \, \xi_+(\rho,\tau),-e^{-\tau}).\] 
Using Proposition \ref{higher.derivative.bounds.in.cylindrical.region}, we obtain 
\[\Big | \frac{\partial}{\partial \rho} \xi_+(\rho,\tau) \Big | \leq C(\theta) \, (-\tau)^{\frac{1}{2}}\] 
and 
\[\Big | \frac{\partial^2}{\partial \rho^2} \xi_+(\rho,\tau) \Big | \leq C(\theta) \, (-\tau)\] 
in the region $\{\frac{\theta}{8} \leq \rho \leq 2\theta\}$. This finally implies 
\[\Big | \frac{\partial}{\partial \rho} \Big ( \frac{\xi_+(\rho,\tau)^2}{4} \Big ) \Big | \leq C(\theta) \, (-\tau)\] 
and 
\[\Big | \frac{\partial^2}{\partial \rho^2} \Big ( \frac{\xi_+(\rho,\tau)^2}{4} \Big ) \Big | \leq C(\theta) \, (-\tau)^{\frac{3}{2}}\] 
in the region $\{\frac{\theta}{8} \leq \rho \leq 2\theta\}$. This proves the assertion. \\

\begin{proposition}
\label{derivative.of.xi.wrt.tau}
Fix a small number $\theta>0$. Then, for $-\tau$ large, we have 
\[\Big | \frac{\partial}{\partial \tau} \Big ( \frac{\xi_+(\rho,\tau)^2}{4} \Big ) \Big | \leq o(1) \, (-\tau)\] 
in the region $\{\frac{\theta}{8} \leq \rho \leq 2\theta\}$.
\end{proposition}

\textbf{Proof.} 
Let us fix a small number $\theta>0$. Differentiating the relation $F(e^{-\frac{\tau}{2}} \, \xi_+(\rho,\tau),-e^{-\tau}) = e^{-\frac{\tau}{2}} \rho$ with respect to $\tau$ gives  
\begin{align*} 
-\frac{1}{2} \, \rho 
&= e^{-\frac{\tau}{2}} \, F_t(e^{-\frac{\tau}{2}} \, \xi_+(\rho,\tau),-e^{-\tau}) \\ 
&+ \Big ( \frac{\partial}{\partial \tau} \xi_+(\rho,\tau) - \frac{1}{2} \, \xi_+(\rho,\tau) \Big ) \, F_z(e^{-\frac{\tau}{2}} \, \xi_+(\rho,\tau),-e^{-\tau}). 
\end{align*} 
Using Proposition \ref{derivative.of.F.wrt.t}, we obtain 
\begin{align*} 
e^{-\frac{\tau}{2}} \, F_t(e^{-\frac{\tau}{2}} \, \xi_+(\rho,\tau),-e^{-\tau}) 
&= -(1+o(1)) \, e^{-\frac{\tau}{2}} \, F(e^{-\frac{\tau}{2}} \, \xi_+(\rho,\tau),-e^{-\tau})^{-1} \\ 
&= -(1+o(1)) \, \rho^{-1} 
\end{align*} 
for $\frac{\theta}{8} \leq \rho \leq 2\theta$. Moreover, 
\[F_z(e^{-\frac{\tau}{2}} \, \xi_+(\rho,\tau),-e^{-\tau}) = -(1+o(1)) \, \Big ( 1-\frac{\rho^2}{2} \Big )^{\frac{1}{2}} \, \rho^{-1} \, (-\tau)^{-\frac{1}{2}}\] 
for $\frac{\theta}{8} \leq \rho \leq 2\theta$. Putting these facts together, we obtain 
\[\frac{\partial}{\partial \tau} \xi_+(\rho,\tau) - \frac{1}{2} \, \xi_+(\rho,\tau) = -(1+o(1)) \, \Big ( 1-\frac{\rho^2}{2} \Big )^{\frac{1}{2}} \, (-\tau)^{\frac{1}{2}}\] 
for $\frac{\theta}{8} \leq \rho \leq 2\theta$. Moreover, 
\[\xi_+(\rho,\tau) = (2+o(1)) \, \Big ( 1-\frac{\rho^2}{2} \Big )^{\frac{1}{2}} \, (-\tau)^{\frac{1}{2}}\] 
for $\frac{\theta}{8} \leq \rho \leq 2\theta$. Thus, we conclude that 
\[\frac{\partial}{\partial \tau} \xi_+(\rho,\tau) = o(1) \, (-\tau)^{\frac{1}{2}}\] 
for $\frac{\theta}{8} \leq \rho \leq 2\theta$. This finally implies 
\[\frac{\partial}{\partial \tau} \Big ( \frac{\xi_+(\rho,\tau)^2}{4} \Big ) = o(1) \, (-\tau)\] 
for $\frac{\theta}{8} \leq \rho \leq 2\theta$. \\

\section{The tip region weights $\mu_+(\rho,\tau)$ and $\mu_-(\rho,\tau)$}

\label{weights}

In this section, we define weights $\mu_+(\rho,\tau)$ and $\mu_-(\rho,\tau)$ which will be needed in the analysis of the linearized equation in the tip region. Let $\theta>0$ be a small positive number, and let $\zeta: \mathbb{R} \to [0,1]$ be a smooth, monotone increasing cutoff function satisfying $\zeta(\rho) = 0$ for $\rho \leq \frac{\theta}{8}$ and $\zeta(\rho) = 1$ for $\rho \geq \frac{\theta}{4}$. We define the weight $\mu_+(\rho,\tau)$ by 
\begin{align*}
\mu_+(\rho,\tau) &= -\zeta(\rho) \, \frac{\xi_+(\rho,\tau)^2}{4} - \int_\rho^\theta \zeta'(\tilde{\rho}) \, \frac{\xi_+(\tilde{\rho},\tau)^2}{4} \, d\tilde{\rho} \\ 
&- \int_\rho^\theta (1-\zeta(\tilde{\rho})) \, \tilde{\rho}^{-1} \, \big ( \Phi((-\tau)^{\frac{1}{2}} \, \tilde{\rho})^{-1}-1 \big ) \, d\tilde{\rho},
\end{align*}
where $\Phi$ denotes the profile of the Bryant soliton. We can define a weight $\mu_-(\rho,\tau)$ in analogous fashion. Of course, the cutoff function $\zeta$ and the weights $\mu_+(\rho,\tau)$ and $\mu_-(\rho,\tau)$ depend on the choice of the parameter $\theta$, but we suppress that dependence in our notation. \\

\begin{lemma} 
\label{weight}
The weight $\mu_+(\rho,\tau)$ satisfies $\mu_+(\rho,\tau) = -\frac{\xi_+(\rho,\tau)^2}{4}$ for $\rho \geq \frac{\theta}{4}$. Moreover, $\mu_+(\rho,\tau) \leq 0$ for all $\rho \leq \frac{\theta}{4}$.
\end{lemma}

\textbf{Proof.} 
This follows immediately from the definition of $\mu_+(\rho,\tau)$. \\

\begin{lemma}
\label{derivative.of.mu.wrt.rho}
Fix a small number $\eta>0$. Then we can find a small number $\theta \in (0,\eta)$ (depending on $\eta$) such that, for $-\tau$ sufficiently large, we have 
\[\Big | \frac{\partial \mu_+}{\partial \rho}(\rho,\tau) - \rho^{-1} \, (V_+(\rho,\tau)^{-2} - 1) \Big | \leq \eta \, \rho^{-1} \, (V_+(\rho,\tau)^{-2} - 1)\] 
in the tip region $\{\rho \leq 2\theta\}$. 
\end{lemma} 

\textbf{Proof.} 
We compute 
\[\frac{\partial \mu_+}{\partial \rho}(\rho,\tau) = -\zeta(\rho) \, \frac{\partial}{\partial \rho} \Big ( \frac{\xi_+(\rho,\tau)^2}{4} \Big ) + (1-\zeta(\rho)) \, \rho^{-1} \, \big ( \Phi((-\tau)^{\frac{1}{2}} \, \rho)^{-1}-1 \big ).\] 
This gives 
\begin{align*} 
&\frac{\partial \mu_+}{\partial \rho}(\rho,\tau) - \rho^{-1} \, (V_+(\rho,\tau)^{-2} - 1) \\ 
&= -\zeta(\rho) \, \Big ( \frac{\partial}{\partial \rho} \Big ( \frac{\xi_+(\rho,\tau)^2}{4} \Big ) + \rho^{-1} \, (V_+(\rho,\tau)^{-2}-1) \Big ) \\ 
&- (1-\zeta(\rho)) \, \rho^{-1} \, \big ( V_+(\rho,\tau)^{-2} - \Phi((-\tau)^{\frac{1}{2}} \, \rho)^{-1} \big ). 
\end{align*}
Hence, the assertion follows from Proposition \ref{comparison.of.V.with.bryant.soliton.profile} and Proposition \ref{derivative.of.xi.wrt.rho}. \\

\begin{lemma}
\label{second.derivative.of.mu.wrt.rho}
If we choose $\theta>0$ sufficiently small, then the following holds. If $-\tau$ is sufficiently large (depending on $\theta$), then 
\[\frac{\partial^2 \mu_+}{\partial \rho^2}(\rho,\tau) \leq \frac{1}{4} \, \Big ( \frac{\partial \mu_+}{\partial \rho}(\rho,\tau) \Big )^2 + \frac{K_*}{4} \, \rho^{-2}\] 
in the tip region $\{\rho \leq 2\theta\}$. Here, $K_*$ is a universal constant which is independent of $\theta$.
\end{lemma}

\textbf{Proof.} 
We compute 
\begin{align*} 
\frac{\partial^2 \mu_+}{\partial \rho^2}(\rho,\tau) 
&= -\zeta(\rho) \, \frac{\partial^2}{\partial \rho^2} \Big ( \frac{\xi_+(\rho,\tau)^2}{4} \Big ) - \zeta'(\rho) \, \frac{\partial}{\partial \rho} \Big ( \frac{\xi_+(\rho,\tau)^2}{4} \Big ) \\ 
&- [1-\zeta(\rho) + \rho \, \zeta'(\rho)] \, \rho^{-2} \, \big ( \Phi((-\tau)^{\frac{1}{2}} \, \rho)^{-1}-1 \big ) \\ 
&- (1-\zeta(\rho)) \, (-\tau)^{\frac{1}{2}} \, \rho^{-1} \, \Phi((-\tau)^{\frac{1}{2}} \, \rho)^{-2} \, \Phi'((-\tau)^{\frac{1}{2}} \, \rho).  
\end{align*}
Recall that $0 \leq \zeta \leq 1$ and $\zeta' \geq 0$. Moreover, we have $\Phi(r)^{-1}-1 \geq \frac{1}{K} \, r^2$ and $|\Phi(r)^{-2} \, \Phi'(r)| \leq Kr$ for all $r \in [0,\infty)$, where $K$ is a universal constant. This implies 
\begin{align*} 
\frac{\partial^2 \mu_+}{\partial \rho^2}(\rho,\tau) 
&\leq -\zeta(\rho) \, \frac{\partial^2}{\partial \rho^2} \Big ( \frac{\xi_+(\rho,\tau)^2}{4} \Big ) - \zeta'(\rho) \, \frac{\partial}{\partial \rho} \Big ( \frac{\xi_+(\rho,\tau)^2}{4} \Big ) \\ 
&+ K \, (1-\zeta(\rho)) \, (-\tau), 
\end{align*}
where $K$ is a universal constant which is independent of $\theta$. Using Proposition \ref{second.derivative.of.xi.wrt.rho}, we obtain 
\[\frac{\partial^2 \mu_+}{\partial \rho^2}(\rho,\tau) \leq o(1) \, (-\tau)^2\] 
in the region $\{\frac{\theta}{8} \leq \rho \leq 2\theta\}$, and 
\[\frac{\partial^2 \mu_+}{\partial \rho^2}(\rho,\tau) \leq K \, (-\tau)\] 
in the region $\{\rho \leq \frac{\theta}{8}\}$. In the next step, we apply Proposition \ref{comparison.of.V.with.bryant.soliton.profile} and Lemma \ref{derivative.of.mu.wrt.rho} with $\eta = \frac{1}{2}$. If we choose $\theta>0$ sufficiently small, then Proposition \ref{comparison.of.V.with.bryant.soliton.profile} and Lemma \ref{derivative.of.mu.wrt.rho} imply 
\begin{align*} 
\frac{\partial \mu_+}{\partial \rho}(\rho,\tau) 
&\geq \frac{1}{2} \, \rho^{-1} \, (V_+(\rho,\tau)^{-2}-1) \\ 
&\geq \frac{1}{4} \, \rho^{-1} \, (\Phi((-\tau)^{\frac{1}{2}} \rho)^{-1}-1) \\ 
&\geq \frac{1}{4K} \, (-\tau) \, \rho 
\end{align*}
in the region $\{\rho \leq 2\theta\}$, where again $K$ is a universal constant independent of $\theta$. Hence, if $-\tau$ is sufficiently large (depending on $\theta$), then we have 
\[\frac{\partial^2 \mu_+}{\partial \rho^2}(\rho,\tau) \leq \frac{1}{4} \, \Big ( \frac{\partial \mu_+}{\partial \rho}(\rho,\tau) \Big )^2 + 16K^4 \, \rho^{-2}\] 
in the region $\{\rho \leq 2\theta\}$. This completes the proof of Lemma \ref{second.derivative.of.mu.wrt.rho}. \\

\begin{lemma} 
\label{derivative.of.mu.wrt.tau}
Let us fix a small number $\theta>0$. Then, for $-\tau$ large, we have 
\[\Big | \frac{\partial \mu_+}{\partial \tau}(\rho,\tau) \Big | \leq o(1) \, (-\tau)\] 
in the tip region $\{\rho \leq 2\theta\}$.
\end{lemma} 

\textbf{Proof.} 
We compute 
\begin{align*} 
\frac{\partial \mu_+}{\partial \tau}(\rho,\tau) &= -\zeta(\rho) \, \frac{\partial}{\partial \tau} \Big ( \frac{\xi_+(\rho,\tau)^2}{4} \Big ) - \int_\rho^\theta \zeta'(\tilde{\rho}) \, \frac{\partial}{\partial \tau} \Big ( \frac{\xi_+(\tilde{\rho},\tau)^2}{4} \Big ) \, d\tilde{\rho} \\ 
&- \frac{1}{2} \, (-\tau)^{-\frac{1}{2}} \int_\rho^\theta (1-\zeta(\tilde{\rho})) \, \Phi((-\tau)^{\frac{1}{2}} \, \tilde{\rho})^{-2} \, \Phi'((-\tau)^{\frac{1}{2}} \, \tilde{\rho}) \, d\tilde{\rho}.
\end{align*}
Note that $|\Phi(r)^{-2} \, \Phi'(r)| \leq Kr$ for all $r \in [0,\infty)$. This gives 
\begin{align*} 
\Big | \frac{\partial \mu_+}{\partial \tau}(\rho,\tau) \Big | 
&\leq \Big | \zeta(\rho) \, \frac{\partial}{\partial \tau} \Big ( \frac{\xi_+(\rho,\tau)^2}{4} \Big ) \Big | + \bigg | \int_{\frac{\theta}{8}}^\theta \zeta'(\tilde{\rho}) \, \frac{\partial}{\partial \tau} \Big ( \frac{\xi_+(\tilde{\rho},\tau)^2}{4} \Big ) \, d\tilde{\rho} \bigg | \\ 
&+ K \int_{\frac{\theta}{8}}^\theta (1-\zeta(\tilde{\rho})) \, \tilde{\rho} \, d\tilde{\rho}
\end{align*} 
for $\rho \leq 2\theta$. Here, $K$ is constant which is independent of $\theta$. Using Proposition \ref{derivative.of.xi.wrt.tau}, we obtain 
\[\Big | \frac{\partial \mu_+}{\partial \tau}(\rho,\tau) \Big | \leq o(1) \, (-\tau)\] 
for $\rho \leq 2\theta$. This completes the proof. \\

In the remainder of this section, we establish a weighted Poincar\'e inequality. 

\begin{proposition}
\label{Poincare.inequality}
If we choose $\theta>0$ sufficiently small, then the following holds. If $-\tau$ is sufficiently large (depending on $\theta$), then 
\[\int_0^{2\theta} \Big ( \frac{\partial \mu_+}{\partial \rho} \Big )^2 \, f^2 \, e^{-\mu_+} \, d\rho \leq 8 \int_0^{2\theta} \Big ( \frac{\partial f}{\partial \rho} \Big )^2 \, e^{-\mu_+} \, d\rho + K_* \int_0^{2\theta} \rho^{-2} \, f^2 \, e^{-\mu_+} \, d\rho\] 
for every smooth function $f$ which is supported in the region $\{\rho \leq 2\theta\}$. Here, $K_*$ is the constant in Lemma \ref{second.derivative.of.mu.wrt.rho}; in particular, $K_*$ is a universal constant which is independent of $\theta$. Note that the right hand side is infinite unless $f(0)=0$.
\end{proposition}

\textbf{Proof.} 
We compute 
\[\frac{\partial}{\partial \rho} \Big ( \frac{\partial \mu_+}{\partial \rho} \, f^2 \, e^{-\mu_+} \Big ) = \frac{\partial^2 \mu_+}{\partial \rho^2} \, f^2 \, e^{-\mu_+} + 2 \, \frac{\partial \mu_+}{\partial \rho} \, f \, \frac{\partial f}{\partial \rho} \, e^{-\mu_+} - \Big ( \frac{\partial \mu_+}{\partial \rho} \Big )^2 \, f^2 \, e^{-\mu_+}.\] 
Using Young's inequality, we obtain 
\[\frac{\partial}{\partial \rho} \Big ( \frac{\partial \mu_+}{\partial \rho} \, f^2 \, e^{-\mu_+} \Big ) \leq \frac{\partial^2 \mu_+}{\partial \rho^2} \, f^2 \, e^{-\mu_+} + 2 \, \Big ( \frac{\partial f}{\partial \rho} \Big )^2 \, e^{-\mu_+} - \frac{1}{2} \, \Big ( \frac{\partial \mu_+}{\partial \rho} \Big )^2 \, f^2 \, e^{-\mu_+}.\] 
Hence, Lemma \ref{second.derivative.of.mu.wrt.rho} gives 
\[\frac{\partial}{\partial \rho} \Big ( \frac{\partial \mu_+}{\partial \rho} \, f^2 \, e^{-\mu_+} \Big ) \leq 2 \, \Big ( \frac{\partial f}{\partial \rho} \Big )^2 \, e^{-\mu_+} - \frac{1}{4} \, \Big ( \frac{\partial \mu_+}{\partial \rho} \Big )^2 \, f^2 \, e^{-\mu_+} + \frac{K_*}{4} \, \rho^{-2} \, f^2 \, e^{-\mu_+}.\] 
From this, the assertion follows. \\

\section{Overview of the proof of Theorem \ref{uniqueness.theorem}}

\label{overview}

We now consider two ancient $\kappa$-solutions $(S^3,g_1(t))$ and $(S^3,g_2(t))$. We assume throughout that neither $(S^3,g_1(t))$ nor $(S^3,g_2(t))$ is a family of shrinking round spheres. We know that both solutions are rotationally symmetric. Let us choose reference points $q_1,q_2 \in S^3$ such that 
\[\limsup_{t \to -\infty} (-t) \, R_{g_1(t)}(q_1) \leq 100 \qquad  \mbox{and} \qquad  \limsup_{t \to -\infty} (-t) \, R_{g_2(t)}(q_2) \leq 100.\] 
Let $F_1(z,t)$ denote the radius of a sphere of symmetry in $(S^3,g_1(t))$ which has signed distance $z$ from the reference point $q_1$. Similarly, let $F_2(z,t)$ denote the radius of a sphere of symmetry in $(S^3,g_1(t))$ which has signed distance $z$ from the reference point $q_2$. 

The functions $F_1(z,t)$ and $F_2(z,t)$ satisfy the PDE 
\begin{align*} 
F_t(z,t)  
&= F_{zz}(z,t) - F(z,t)^{-1} \, (1-F_z(z,t)^2) \\ 
&- 2 \, F_z(z,t) \int_0^z \frac{F_{zz}(z',t)}{F(z',t)} \, dz'. 
\end{align*} 
In the next step, we replace the function $F_2(z,t)$ by a new function $F_2^{\alpha\beta\gamma}(z,t)$. Here, $(\alpha,\beta,\gamma)$ is a triplet of real numbers satisfying the following admissibility condition: \\

\begin{definition}
\label{admissibility}
Given a real number $\varepsilon \in (0,\frac{1}{2})$, we say that the triplet $(\alpha,\beta,\gamma)$ is $\varepsilon$-admissible with respect to time $t_*$ if 
\[|\alpha| \leq \varepsilon \sqrt{-t_*}, \qquad |\beta| \leq \varepsilon \, \frac{(-t_*)}{\log(-t_*)}, \qquad |\gamma| \leq \varepsilon \log (-t_*).\]
\end{definition}

In the following, we consider a time $t_* < 0$, where $-t_*$ is very large. Suppose that $(\alpha,\beta,\gamma)$ is a triplet of real numbers which is $\varepsilon$-admissible with respect to time $t_*$ for some $\varepsilon \in (0,\frac{1}{2})$. For each $t \leq t_*$, we consider the rescaled metrics
\[g_2^{\beta\gamma}(t) := e^\gamma \, g_2(e^{-\gamma} (t-\beta)).\] 
Moreover, let $F_2^{\beta\gamma}(z,t)$ denote the profile function associated with the metric $g_2^{\beta\gamma}(t)$. In other words, $F_2^{\beta\gamma}(z,t)$ denotes the radius of a sphere of symmetry in $(S^3,g_2^{\beta\gamma}(t))$ which has signed distance $z$ from the reference point $q_2$. Clearly, 
\[F_2^{\beta\gamma}(z,t) = e^{\frac{\gamma}{2}} \, F_2(e^{-\frac{\gamma}{2}} z,e^{-\gamma} (t-\beta)).\] 
Note that the metrics $g_2^{\beta\gamma}(t)$ form a solution to the Ricci flow, and the profile function $F_2^{\beta\gamma}(z,t)$ satisfies the same PDE as the original profile function $F_2(z,t)$. In the next step, we choose a new base point $q_2^{\alpha\beta\gamma}$ with the property that the sphere of symmetry passing through $q_2^{\alpha\beta\gamma}$ has signed distance $\alpha$ from the point $q_2$ with respect to the metric $g_2^{\beta\gamma}(t_*)$. For each time $t \leq t_*$, we denote by $s^{\alpha\beta\gamma}(t)$ the signed distance of the sphere of symmetry passing through $q_2^{\alpha\beta\gamma}$ from the point $q_2$ with respect to the metric $g_2^{\beta\gamma}(t)$. The function $s^{\alpha\beta\gamma}(t)$ can be characterized as the solution of the ODE 
\[\frac{d}{dt} s^{\alpha\beta\gamma}(t) = 2 \int_0^{s^{\alpha\beta\gamma}(t)} \frac{F_{2,zz}^{\beta\gamma}(z',t)}{F_2^{\beta\gamma}(z',t)} \, dz', \qquad s^{\alpha\beta\gamma}(t_*) = \alpha.\] 
For each time $t \leq t_*$, we denote by $F_2^{\alpha\beta\gamma}(z,t)$ the radius of the sphere of symmetry in $(S^3,g_2^{\beta\gamma}(t))$ which has signed distance $z$ from the point $q_2^{\alpha\beta\gamma}$. Clearly, the function $F_2^{\alpha\beta\gamma}(z,t)$ satisfies the same PDE as the function $F_2(z,t)$. Moreover, the function $F_2^{\alpha\beta\gamma}(z,t)$ is related to the function $F_2^{\beta\gamma}(z,t)$ by the formula 
\[F_2^{\alpha\beta\gamma}(z,t) = F_2^{\beta\gamma}(z+s^{\alpha\beta\gamma}(t),t) = e^{\frac{\gamma}{2}} \, F_2 \big ( e^{-\frac{\gamma}{2}} (z + s^{\alpha\beta\gamma}(t)),e^{-\gamma} (t-\beta) \big ).\]
In other words, the modified solution $F_2^{\alpha\beta\gamma}(z,t)$ differs from $F_2^{\beta\gamma}(z,t)$ by a translation in space. In particular, for $t=t_*$, we obtain 
\[F_2^{\alpha\beta\gamma}(z,t_*) = F_2^{\beta\gamma}(z+\alpha,t_*) = e^{\frac{\gamma}{2}} \, F_2(e^{-\frac{\gamma}{2}} (z+\alpha),e^{-\gamma} (t_*-\beta)).\]

\begin{lemma}
\label{bound.for.s^alpha}
If $-t_*$ is sufficiently large, then the following holds. Suppose that the triplet $(\alpha,\beta,\gamma)$ is $\varepsilon$-admissible with respect to time $t_*$, where $\varepsilon \in (0,\frac{1}{2})$. Let $s^{\alpha\beta\gamma}(t)$ denote the solution of the ODE 
\[\frac{d}{dt} s^{\alpha\beta\gamma}(t) = 2 \int_0^{s^{\alpha\beta\gamma}(t)} \frac{F_{2,zz}^{\beta\gamma}(z',t)}{F_2^{\beta\gamma}(z',t)} \, dz'\] 
with terminal condition $s^{\alpha\beta\gamma}(t_*) = \alpha$. Then $|s^{\alpha\beta\gamma}(t)| \leq \varepsilon \sqrt{-t}$ for all $t \leq t_*$. 
\end{lemma}

\textbf{Proof.} 
Recall that the reference point $q_2$ has been chosen such that the blow-down limit of $(S^3,g_2(t))$ around $q_2$ is a cylinder. Hence, if we choose $-t_*$ is sufficiently large, then 
\[0 \leq -\frac{F_{2,zz}(z,t)}{F_2(z,t)} \leq \frac{1}{(-8t)}\] 
whenever $t \leq -\frac{1}{2} \, \sqrt{-t_*}$ and $|z| \leq \sqrt{-2t}$. In the next step, we replace $t$ by $e^{-\gamma} (t-\beta)$, and we replace $z$ by $e^{-\frac{\gamma}{2}} z$. This gives 
\[0 \leq -\frac{F_{2,zz}^{\beta\gamma}(z,t)}{F_2^{\beta\gamma}(z,t)} \leq \frac{1}{(-8(t-\beta))}\] 
whenever $t-\beta \leq -\frac{1}{2} \, e^\gamma \sqrt{-t_*}$ and $|z| \leq \sqrt{-2(t-\beta)}$. The condition $|\gamma| \leq \varepsilon \log (-t_*) \leq \frac{1}{2} \log(-t_*)$ implies $t_* \leq -e^\gamma \sqrt{-t_*}$. Moreover, the condition $|\beta| \leq \varepsilon \, \frac{(-t_*)}{\log(-t_*)} \leq \frac{1}{2} \, \frac{(-t_*)}{\log(-t_*)}$ ensures that $t-\beta \leq \frac{1}{2} \, t$ for all $t \leq t_*$. Consequently, 
\[0 \leq -\frac{F_{2,zz}^{\beta\gamma}(z,t)}{F_2^{\beta\gamma}(z,t)} \leq \frac{1}{(-4t)}\] 
whenever $t \leq t_*$ and $|z| \leq \sqrt{-t}$. Hence, if $s^{\alpha\beta\gamma}(t)$ is a solution of the ODE above, then 
\[\Big | \frac{d}{dt} s^{\alpha\beta\gamma}(t) \Big | \leq \frac{1}{(-2t)} \, |s^{\alpha\beta\gamma}(t)|\] 
whenever $t \leq t_*$ and $|s^{\alpha\beta\gamma}(t)| \leq \sqrt{-t}$. From this, we deduce that 
\[\frac{d}{dt} \big ( (-t)^{-1} \, (s^{\alpha\beta\gamma}(t))^2 \big ) = (-t)^{-2} \, (s^{\alpha\beta\gamma}(t))^2 + 2 \, (-t)^{-1} \, s^{\alpha\beta\gamma}(t) \, \frac{d}{dt} s^{\alpha\beta\gamma}(t) \geq 0\] 
whenever $t \leq t_*$ and $|s^{\alpha\beta\gamma}(t)| \leq \sqrt{-t}$. By assumption, $(-t_*)^{-1} \, (s^{\alpha\beta\gamma}(t_*))^2 = (-t_*)^{-1} \, \alpha^2 \leq \varepsilon^2$. Since $\varepsilon \in (0,\frac{1}{2})$, we conclude that $(-t)^{-1} \, (s^{\alpha\beta\gamma}(t))^2 \leq \varepsilon^2$ for all $t \leq t_*$. This completes the proof of Lemma \ref{bound.for.s^alpha}. \\

Using the admissibility conditions in Definition \ref{admissibility} and Lemma \ref{bound.for.s^alpha}, we can estimate the modified profile function $F_2^{\alpha\beta\gamma}$: \\

\begin{proposition}
\label{estimate.for.modified.profile}
Fix a small number $\theta > 0$ and a small number $\eta > 0$. Then there exists a small number $\varepsilon>0$ (depending on $\theta$ and $\eta$) with the following property. If the triplet $(\alpha,\beta,\gamma)$ is $\varepsilon$-admissible with respect to time $t_*$ and $-t_*$ is sufficiently large, then  
\[\Big | \frac{1}{2} \, F_2^{\alpha\beta\gamma}(z,t)^2 + t + \frac{z^2+2t}{4 \log(-t)} \Big | \leq \eta \, \frac{z^2-t}{\log(-t)}\] 
and 
\[\Big | F_2^{\alpha\beta\gamma}(z,t) \, F_{2z}^{\alpha\beta\gamma}(z,t) + \frac{z}{2 \log(-t)} \Big | \leq \eta \, \frac{|z|+\sqrt{-t}}{\log(-t)}\] 
whenever $F_2^{\alpha\beta\gamma}(z,t) \geq \frac{\theta}{10} \sqrt{-t}$ and $t \leq t_*$.
\end{proposition}

\textbf{Proof.} 
Using Proposition \ref{precise.estimate.for.F} and Proposition \ref{precise.estimate.for.F_z}, we obtain 
\[\Big | \frac{1}{2} \, F_2(z,t)^2 + t + \frac{z^2+2t}{4 \log(-t)} \Big | \leq \frac{\eta}{4} \, \frac{z^2-t}{\log(-t)}\] 
and 
\[\Big | F_2(z,t) \, F_{2z}(z,t) + \frac{z}{2 \log(-t)} \Big | \leq \frac{\eta}{4} \, \frac{|z|+\sqrt{-t}}{\log(-t)}\] 
whenever $F_2(z,t) \geq \frac{\theta}{20} \sqrt{-t}$ and $-t$ is sufficiently large. We now replace $t$ by $e^{-\gamma} (t-\beta)$ and $z$ by $e^{-\frac{\gamma}{2}} z$. This gives 
\[\Big | \frac{1}{2} \, F_2^{\beta\gamma}(z,t)^2 + (t-\beta) + \frac{z^2+2(t-\beta)}{4 \log(-(t-\beta)) - 4\gamma} \Big | \leq \frac{\eta}{4} \, \frac{z^2-(t-\beta)}{\log(-(t-\beta))-\gamma}\] 
and 
\[\Big | F_2^{\beta\gamma}(z,t) \, F_{2z}^{\beta\gamma}(z,t) + \frac{z}{2 \log(-(t-\beta)) - 2\gamma} \Big | \leq \frac{\eta}{4} \, \frac{|z|+\sqrt{-(t-\beta)}}{\log(-(t-\beta))-\gamma}\] 
whenever $F_2^{\beta\gamma}(z,t) \geq \frac{\theta}{20} \sqrt{-(t-\beta)}$ and $-e^{-\gamma} (t-\beta)$ is sufficiently large. By assumption, the triplet $(\alpha,\beta,\gamma)$ is $\varepsilon$-admissible with respect to time $t_*$. If $\varepsilon$ is sufficiently small (depending on $\theta$ and $\eta$) and $-t_*$ is sufficiently large (depending on $\theta$ and $\eta$), then we obtain 
\[\Big | \frac{1}{2} \, F_2^{\beta\gamma}(z,t)^2 + t + \frac{z^2+2t}{4 \log(-t)} \Big | \leq \frac{\eta}{2} \, \frac{z^2-t}{\log(-t)}\] 
and 
\[\Big | F_2^{\beta\gamma}(z,t) \, F_{2z}^{\beta\gamma}(z,t) + \frac{z}{2 \log(-t)} \Big | \leq \frac{\eta}{2} \, \frac{|z|+\sqrt{-t}}{\log(-t)}\] 
whenever $F_2^{\beta\gamma}(z,t) \geq \frac{\theta}{10} \sqrt{-t}$ and $t \leq t_*$. By Lemma \ref{bound.for.s^alpha}, $|s^{\alpha\beta\gamma}(t)| \leq \varepsilon \sqrt{-t}$ for $t \leq t_*$. Hence, we obtain 
\[\Big | \frac{1}{2} \, F_2^{\alpha\beta\gamma}(z,t)^2 + t + \frac{z^2+2t}{4 \log(-t)} \Big | \leq \eta \, \frac{z^2-t}{\log(-t)}\] 
and 
\[\Big | F_2^{\alpha\beta\gamma}(z,t) \, F_{2z}^{\alpha\beta\gamma}(z,t) + \frac{z}{2 \log(-t)} \Big | \leq \eta \, \frac{|z|+\sqrt{-t}}{\log(-t)}\] 
whenever $F_2^{\alpha\beta\gamma}(z,t) \geq \frac{\theta}{10} \sqrt{-t}$ and $t \leq t_*$. This completes the proof of Proposition \ref{estimate.for.modified.profile}. \\

We define functions $U_{1+}(r,t)$ and $U_{1-}(r,t)$ by 
\[U_{1+}(r,t) = \Big ( \frac{\partial}{\partial z} F_1(z,t) \Big )^2\] 
for $r = F_1(z,t)$ and $z \geq 2\sqrt{-t}$ and 
\[U_{1-}(r,t) = \Big ( \frac{\partial}{\partial z} F_1(z,t) \Big )^2\] 
for $r = F_1(z,t)$ and $z \leq -2\sqrt{-t}$. Similarly, we define functions $U_{2+}(r,t)$ and $U_{2-}(r,t)$ by 
\[U_{2+}(r,t) = \Big ( \frac{\partial}{\partial z} F_2(z,t) \Big )^2\] 
for $r = F_2(z,t)$ and $z \geq 2\sqrt{-t}$ and 
\[U_{2-}(r,t) = \Big ( \frac{\partial}{\partial z} F_2(z,t) \Big )^2\] 
for $r = F_2(z,t)$ and $z \leq -2\sqrt{-t}$. Moreover, we define 
\begin{align*} 
U_{2+}^{\beta\gamma}(r,t) &:= U_{2+}(e^{-\frac{\gamma}{2}} r,e^{-\gamma} (t-\beta)), \\ 
U_{2-}^{\beta\gamma}(r,t) &:= U_{2-}(e^{-\frac{\gamma}{2}} r,e^{-\gamma} (t-\beta)). 
\end{align*} 
With this understood, we have 
\[U_{2+}^{\beta\gamma}(r,t) = \Big ( \frac{\partial}{\partial z} F_2^{\alpha\beta\gamma}(z,t) \Big )^2\] 
for $r = F_2^{\alpha\beta\gamma}(z,t)$, $z \geq 4\sqrt{-t}$, and $t \leq t_*$, and 
\[U_{2-}^{\beta\gamma}(r,t) = \Big ( \frac{\partial}{\partial z} F_2^{\alpha\beta\gamma}(z,t) \Big )^2\] 
for $r = F_2^{\alpha\beta\gamma}(z,t)$, $z \leq -4\sqrt{-t}$, and $t \leq t_*$. 

In the next step, we perform the usual rescaling. We put $t = -e^{-\tau}$ and $r = e^{-\frac{\tau}{2}} \rho$. This gives functions $V_{1+}(\rho,\tau)$, $V_{1-}(\rho,\tau)$, $V_{2+}(\rho,\tau)$, $V_{2-}(\rho,\tau)$, $V_{2+}^{\beta\gamma}(\rho,\tau)$, and $V_{2-}^{\beta\gamma}(\rho,\tau)$, where 
\begin{align*} 
V_{1+}(\rho,\tau) &:= \sqrt{U_{1+}(e^{-\frac{\tau}{2}} \rho,-e^{-\tau})}, \\ 
V_{1-}(\rho,\tau) &:= \sqrt{U_{1-}(e^{-\frac{\tau}{2}} \rho,-e^{-\tau})}, \\ 
V_{2+}(\rho,\tau) &:= \sqrt{U_{2+}(e^{-\frac{\tau}{2}} \rho,-e^{-\tau})}, \\ 
V_{2-}(\rho,\tau) &:= \sqrt{U_{2-}(e^{-\frac{\tau}{2}} \rho,-e^{-\tau})}, \\ 
V_{2+}^{\beta\gamma}(\rho,\tau) &:= \sqrt{U_{2+}^{\beta\gamma}(e^{-\frac{\tau}{2}} \rho,-e^{-\tau})}, \\ 
V_{2-}^{\beta\gamma}(\rho,\tau) &:= \sqrt{U_{2-}^{\beta\gamma}(e^{-\frac{\tau}{2}} \rho,-e^{-\tau})}. 
\end{align*} 
A straightforward calculation gives 
\begin{align*} 
V_{2+}^{\beta\gamma}(\rho,\tau) &= V_{2+} \Big ( \frac{\rho}{\sqrt{1+\beta e^\tau}}, \tau + \gamma - \log(1+\beta e^\tau) \Big ), \\ 
V_{2-}^{\beta\gamma}(\rho,\tau) &= V_{2-} \Big ( \frac{\rho}{\sqrt{1+\beta e^\tau}}, \tau + \gamma - \log(1+\beta e^\tau) \Big ). 
\end{align*}

\begin{proposition}
\label{comparison.of.modified.V.with.bryant.soliton.profile}
Fix a small number $\eta>0$. Then we can find a small number $\theta \in (0,\eta)$ (depending on $\eta$) and a small number $\varepsilon>0$ (depending on $\theta$ and $\eta$) with the following property. If the triplet $(\alpha,\beta,\gamma)$ is $\varepsilon$-admissible with respect to time $t_* = -e^{-\tau_*}$ and $-\tau_*$ is sufficiently large, then 
\[|V_{2+}^{\beta\gamma}(\rho,\tau)^{-2} - \Phi((-\tau)^{\frac{1}{2}} \rho)^{-1}| \leq \eta \, (V_{2+}^{\beta\gamma}(\rho,\tau)^{-2}-1)\] 
for $\rho \leq 10\theta$ and $\tau \leq \tau_*$, and 
\[\Big | \frac{\partial}{\partial \tau} V_{2+}^{\beta\gamma}(\rho,\tau)^{-2} \Big | \leq \eta \, \rho^{-2} \, (V_{2+}^{\beta\gamma}(\rho,\tau)^{-1}-1)\] 
for $\rho \leq \frac{1}{8}$ and $\tau \leq \tau_*$. Here, $\Phi$ denotes the profile of the Bryant soliton.
\end{proposition}

\textbf{Proof.} 
By Proposition \ref{comparison.of.V.with.bryant.soliton.profile}, we can choose $\theta$ sufficiently small (depending on $\eta$) so that 
\[|V_{2+}(\rho,\tau)^{-2} - \Phi((-\tau)^{\frac{1}{2}} \rho)^{-1}| \leq \frac{\eta}{4} \, (V_{2+}(\rho,\tau)^{-2}-1)\] 
whenever $\rho \leq 20\theta$ and $-\tau$ is sufficiently large. We now replace $\tau$ by $\tau + \gamma - \log(1+\beta e^\tau)$ and $\rho$ by $\frac{\rho}{\sqrt{1+\beta e^\tau}}$. This gives 
\begin{align*} 
&\bigg | V_{2+}^{\beta\gamma}(\rho,\tau)^{-2} - \Phi \bigg ( \Big ( \frac{1}{1+\beta e^\tau}+\frac{\gamma-\log(1+\beta e^\tau)}{\tau (1+\beta e^\tau)} \Big )^{\frac{1}{2}} \, (-\tau)^{\frac{1}{2}} \rho \bigg )^{-1} \bigg | \\ 
&\leq \frac{\eta}{4} \, (V_{2+}^{\beta\gamma}(\rho,\tau)^{-2}-1) 
\end{align*}
whenever $\rho \leq 20\theta \sqrt{1+\beta e^\tau}$ and $-\tau$ is sufficiently large. 

By assumption, the triplet $(\alpha,\beta,\gamma)$ is $\varepsilon$-admissible with respect to time $t_* = -e^{-\tau_*}$. If we choose $\varepsilon$ sufficiently small (depending on $\eta$) and $-\tau_*$ sufficiently large (depending on $\eta$), then Proposition \ref{estimate.for.normalized.Bryant.profile} implies that 
\begin{align*} 
&\bigg | \Phi \bigg ( \Big ( \frac{1}{1+\beta e^\tau}+\frac{\gamma-\log(1+\beta e^\tau)}{\tau (1+\beta e^\tau)} \Big )^{\frac{1}{2}} \, (-\tau)^{\frac{1}{2}} \rho \bigg )^{-1} - \Phi((-\tau)^{\frac{1}{2}} \rho)^{-1} \bigg | \\ 
&\leq \frac{\eta}{4} \, (\Phi((-\tau)^{\frac{1}{2}} \rho)^{-1}-1)
\end{align*} 
for all $\rho$ and all $\tau \leq \tau_*$. Hence, if we choose $-\tau_*$ sufficiently large (depending on $\theta$ and $\eta$), then we obtain 
\begin{align*} 
&|V_{2+}^{\beta\gamma}(\rho,\tau)^{-2} - \Phi((-\tau)^{\frac{1}{2}} \rho)^{-1}| \\ 
&\leq \frac{\eta}{4} \, (V_{2+}^{\beta\gamma}(\rho,\tau)^{-2}-1) + \frac{\eta}{4} \, (\Phi((-\tau)^{\frac{1}{2}} \rho)^{-1}-1) \\ 
&\leq \frac{\eta}{2} \, (V_{2+}^{\beta\gamma}(\rho,\tau)^{-2}-1) + \frac{\eta}{4} \, |V_{2+}^{\beta\gamma}(\rho,\tau)^{-2} - \Phi((-\tau)^{\frac{1}{2}} \rho)^{-1}|
\end{align*}
whenever $\rho \leq 10\theta$ and $\tau \leq \tau_*$. The last term on the right hand side can be absorbed into the left hand side. This gives 
\[|V_{2+}^{\beta\gamma}(\rho,\tau)^{-2} - \Phi((-\tau)^{\frac{1}{2}} \rho)^{-1}| \leq \eta \, (V_{2+}^{\beta\gamma}(\rho,\tau)^{-2}-1)\] 
whenever $\rho \leq 10\theta$ and $\tau \leq \tau_*$. This proves the first statement. 

We now turn to the second statement. Using Proposition \ref{higher.derivatives.of.V.in.collar.region} and Corollary \ref{derivative.of.V.wrt.tau}, we obtain 
\[\Big | \rho \, \frac{\partial}{\partial \rho} V_{2+}(\rho,\tau) \Big | \leq C \, (V_{2+}(\rho,\tau)^{-1}-1)\] 
and 
\[\Big | \frac{\partial}{\partial \tau} V_{2+}(\rho,\tau) \Big | \leq \frac{\eta}{2} \, \rho^{-2} \, (V_{2+}(\rho,\tau)^{-1}-1)\] 
whenever $\rho \leq \frac{1}{4}$ and $-\tau$ is sufficiently large. Using the identity 
\[V_{2+}^{\beta\gamma}(\rho,\tau) = V_{2+} \Big ( \frac{\rho}{\sqrt{1+\beta e^\tau}}, \tau + \gamma - \log(1+\beta e^\tau) \Big )\] 
and the chain rule, we conclude that 
\[\Big | \frac{\partial}{\partial \tau} V_{2+}^{\beta\gamma}(\rho,\tau) \Big | \leq \frac{\eta}{2} \, \rho^{-2} \, (V_{2+}^{\beta\gamma}(\rho,\tau)^{-1}-1) + C \, \Big | \frac{\beta e^\tau}{1+\beta e^\tau} \Big | \, (V_{2+}^{\beta\gamma}(\rho,\tau)^{-1}-1)\] 
whenever $\rho \leq \frac{1}{8}$ and $\tau \leq \tau_*$. Since the triplet $(\alpha,\beta,\gamma)$ is $\varepsilon$-admissible with respect to time $t_* = -e^{-\tau_*}$, the second statement follows. This completes the proof of Proposition \ref{comparison.of.modified.V.with.bryant.soliton.profile}. \\

We next consider the difference between the two solutions: 
\begin{align*} 
W_+^{\beta\gamma}(\rho,\tau) &:= V_{1+}(\rho,\tau) - V_{2+}^{\beta\gamma}(\rho,\tau), \\ 
W_-^{\beta\gamma}(\rho,\tau) &:= V_{1-}(\rho,\tau) - V_{2-}^{\beta\gamma}(\rho,\tau). 
\end{align*} 
For each $\tau$, we know that $1-V_{1+}(\rho,\tau) = O(\rho^2)$, $1-V_{1-}(\rho,\tau) = O(\rho^2)$, $1-V_{2+}^{\beta\gamma}(\rho,\tau) = O(\rho^2)$, $1-V_{2-}^{\beta\gamma}(\rho,\tau) = O(\rho^2)$ as $\rho \to 0$. Hence, for each $\tau$, we have $W_+^{\beta\gamma}(\rho,\tau) = O(\rho^2)$ and $W_-^{\beta\gamma}(\rho,\tau) = O(\rho^2)$ as $\rho \to 0$. 

Let $\mu_+(\rho,\tau)$ and $\mu_-(\rho,\tau)$ denote the weights associated with the solution $(S^3,g_1(t))$. \\

\begin{proposition}
\label{estimate.for.difference.in.tip.region}
We can choose $\theta>0$ and $\varepsilon>0$ sufficiently small so that the following holds. If $-\tau_*$ is sufficiently large (depending on $\theta$) and the triplet $(\alpha,\beta,\gamma)$ is $\varepsilon$-admissible with respect to time $t_* = -e^{-\tau_*}$, then 
\begin{align*} 
&\sup_{\tau \leq \tau_*} (-\tau)^{-\frac{1}{2}} \int_{\tau-1}^\tau \int_0^\theta V_{1+}^{-2} \, (W_+^{\beta\gamma})^2 \, e^{\mu_+} \\ 
&\leq C(\theta) \, (-\tau_*)^{-1} \sup_{\tau \leq \tau_*} (-\tau)^{-\frac{1}{2}} \int_{\tau-1}^\tau \int_\theta^{2\theta} V_{1+}^{-2} \, (W_+^{\beta\gamma})^2 \, e^{\mu_+}. 
\end{align*} 
An analogous estimate holds for $W_-^{\beta\gamma}$.
\end{proposition}

We will give the proof of Proposition \ref{estimate.for.difference.in.tip.region} in Section \ref{difference.tip.region}. \\

From this point on, we fix $\theta$ small enough so that the conclusion of Proposition \ref{estimate.for.difference.in.tip.region} holds. Let $\chi_{\mathcal{C}}$ denote a smooth, even cutoff function satisfying $\chi_{\mathcal{C}} = 1$ on $[0,\sqrt{4-\frac{\theta^2}{2}}]$ and $\chi_{\mathcal{C}} = 0$ on $[\sqrt{4-\frac{\theta^2}{4}},\infty)$. Moreover, we may assume that $\chi_{\mathcal{C}}$ is monotone decreasing on $[0,\infty)$. We define 
\begin{align*} 
G_1(\xi,\tau) &:= e^{\frac{\tau}{2}} \, F_1(e^{-\frac{\tau}{2}} \xi,-e^{-\tau}) - \sqrt{2}, \\ 
G_2(\xi,\tau) &:= e^{\frac{\tau}{2}} \, F_2(e^{-\frac{\tau}{2}} \xi,-e^{-\tau}) - \sqrt{2}, \\ 
G_2^{\alpha\beta\gamma}(\xi,\tau) &:= e^{\frac{\tau}{2}} \, F_2^{\alpha\beta\gamma}(e^{-\frac{\tau}{2}} \xi,-e^{-\tau}) - \sqrt{2}. 
\end{align*} 
Let
\[H^{\alpha\beta\gamma}(\xi,\tau) := G_1(\xi,\tau) - G_2^{\alpha\beta\gamma}(\xi,\tau)\] 
and 
\[H_{\mathcal{C}}^{\alpha\beta\gamma}(\xi,\tau) := \chi_{\mathcal{C}}((-\tau)^{-\frac{1}{2}} \xi) \, H^{\alpha\beta\gamma}(\xi,\tau).\] 
Using the PDEs for $G_1$ and $G_2^{\alpha\beta\gamma}$, we can derive a PDE for the function $H^{\alpha\beta\gamma}$. The leading term in that PDE is given by the operator 
\[\mathcal{L} f := f_{\xi\xi} - \frac{1}{2} \, \xi \, f_\xi + f.\] 
To analyze this operator, we perform a spectral decomposition. As in \cite{Angenent-Daskalopoulos-Sesum2}, we consider the Hilbert space $\mathcal{H} = L^2(\mathbb{R},e^{-\frac{\xi^2}{4}} \, d\xi)$. The Hilbert space $\mathcal{H}$ has a natural direct sum decomposition $\mathcal{H} = \mathcal{H}_+ \oplus \mathcal{H}_0 \oplus \mathcal{H}_-$. Here, $\mathcal{H}_+$ is a two-dimensional subspace spanned by the functions $1$ and $\xi$; $\mathcal{H}_0$ is a one-dimensional subspace spanned by the function $\xi^2-2$; and $\mathcal{H}_-$ is the orthogonal complement of $\mathcal{H}_+ \oplus \mathcal{H}_0$. Finally, let $P_+$, $P_0$, and $P_-$ denote the projection operators associated to the direct sum decomposition $\mathcal{H} = \mathcal{H}_+ \oplus \mathcal{H}_0 \oplus \mathcal{H}_-$.

With this understood, we write 
\[P_0 H_{\mathcal{C}}^{\alpha\beta\gamma}(\xi,\tau) = \sqrt{2} \, a^{\alpha\beta\gamma}(\tau) \, (\xi^2-2),\] 
where 
\[a^{\alpha\beta\gamma}(\tau) := \frac{1}{16\sqrt{2\pi}} \int_{\mathbb{R}} e^{-\frac{\xi^2}{4}} \, (\xi^2-2) \, H_{\mathcal{C}}^{\alpha\beta\gamma}(\xi,\tau) \, d\xi.\] 
Moreover, we put $\hat{H}_{\mathcal{C}}^{\alpha\beta\gamma} = P_+ H_{\mathcal{C}}^{\alpha\beta\gamma} + P_- H_{\mathcal{C}}^{\alpha\beta\gamma}$. \\

\begin{proposition}
\label{choice.of.parameters}
Fix $\theta>0$ and $\varepsilon>0$ small enough so that the conclusion of Proposition \ref{estimate.for.difference.in.tip.region} holds. Let $\delta \in (0,\varepsilon)$ be given. If $-\tau_*$ is sufficiently large (depending on $\delta$), then we can find a triplet $(\alpha,\beta,\gamma)$ (depending on $\tau_*$) such that $P_+ H_{\mathcal{C}}^{\alpha\beta\gamma} = 0$ and $P_0 H_{\mathcal{C}}^{\alpha\beta\gamma} = 0$ at time $\tau_*$. Moreover, if $-\tau_*$ is sufficiently large (depending on $\delta$), then the triplet $(\alpha,\beta,\gamma)$ is $\delta$-admissible with respect to time $t_* = -e^{-\tau_*}$.
\end{proposition}

\textbf{Proof.} 
Using the identity $s^{\alpha\beta\gamma}(t_*) = \alpha$, we obtain  
\[F_2^{\alpha\beta\gamma}(z,t_*) = e^{\frac{\gamma}{2}} \, F_2(e^{-\frac{\gamma}{2}} (z+\alpha),e^{-\gamma} (t_*-\beta)).\] 
Consequently, 
\begin{align*} 
G_2^{\alpha\beta\gamma}(\xi,\tau_*) 
&= \sqrt{1+\beta e^{\tau_*}} \, G_2 \Big ( \frac{\xi + \alpha e^{\frac{\tau_*}{2}}}{\sqrt{1+\beta e^{\tau_*}}},\tau_* + \gamma - \log(1+\beta e^{\tau_*}) \Big ) \\ 
&+ \sqrt{2} \, (\sqrt{1+\beta e^{\tau_*}} - 1). 
\end{align*} 
The proof of Proposition \ref{choice.of.parameters} now proceeds as in \cite{Angenent-Daskalopoulos-Sesum2}. This argument relies only on the asymptotics of our solution in the cylindrical region. Since the asymptotics of our ancient solutions to Ricci flow in the cylindrical region are very similar to the cylindrical region asymptotics of ancient solutions to mean curvature flow, the proof of Proposition \ref{choice.of.parameters} is identical to the proof of the corresponding Proposition 4.1 in \cite{Angenent-Daskalopoulos-Sesum2}. \\

From this point on, we assume that the triplet $(\alpha,\beta,\gamma)$ is chosen as in Proposition \ref{choice.of.parameters}. In particular, this will ensure that $a^{\alpha\beta\gamma}(\tau_*) = 0$. Note that the triplet $(\alpha,\beta,\gamma)$ depends on $\tau_*$ (which we have not yet fixed). \\

\begin{proposition}
\label{estimate.for.difference.in.cylindrical.region}
Fix $\theta>0$ small enough so that the conclusion of Proposition \ref{estimate.for.difference.in.tip.region} holds. Suppose that $-\tau_*$ is sufficiently large, and that the triplet $(\alpha,\beta,\gamma)$ is chosen as in Proposition \ref{choice.of.parameters}. Then 
\begin{align*}
&(-\tau_*) \sup_{\tau \leq \tau_*} \int_{\tau-1}^\tau \int_{\mathbb{R}} e^{-\frac{\xi^2}{4}} \, (\hat{H}_{\mathcal{C},\xi}^{\alpha\beta\gamma}(\xi,\tau')^2 + \hat{H}_{\mathcal{C}}^{\alpha\beta\gamma}(\xi,\tau')^2) \, d\xi \, d\tau' \\ 
&\leq C(\theta) \sup_{\tau \leq \tau_*} \int_{\tau-1}^\tau a^{\alpha\beta\gamma}(\tau')^2 \, d\tau' \\ 
&+ C(\theta) \sup_{\tau \leq \tau_*} \int_{\tau-1}^\tau \int_{\{\sqrt{4-\frac{\theta^2}{2}} \, (-\tau')^{\frac{1}{2}} \leq |\xi| \leq \sqrt{4-\frac{\theta^2}{4}} \, (-\tau')^{\frac{1}{2}}\}} e^{-\frac{\xi^2}{4}} \, H^{\alpha\beta\gamma}(\xi,\tau')^2 \, d\xi \, d\tau'. 
\end{align*}
\end{proposition}

We will give the proof of Proposition \ref{estimate.for.difference.in.cylindrical.region} in Section \ref{difference.cylindrical.region}. \\

By combining Proposition \ref{estimate.for.difference.in.tip.region} and Proposition \ref{estimate.for.difference.in.cylindrical.region}, we can show that in the cylindrical region the norm of $P_0 H_{\mathcal{C}}^{\alpha\beta\gamma}$ dominates over the norm of $\hat{H}_{\mathcal{C}}^{\alpha\beta\gamma}$. More precisely, we have the following result: 

\begin{proposition}
\label{neutral.mode.dominates}
Fix $\theta>0$ small enough so that the conclusion of Proposition \ref{estimate.for.difference.in.tip.region} holds. Suppose that $-\tau_*$ is sufficiently large, and that the triplet $(\alpha,\beta,\gamma)$ is chosen as in Proposition \ref{choice.of.parameters}. Then 
\begin{align*}
&(-\tau_*) \sup_{\tau \leq \tau_*} \int_{\tau-1}^\tau \int_{\mathbb{R}} e^{-\frac{\xi^2}{4}} \, (\hat{H}_{\mathcal{C},\xi}^{\alpha\beta\gamma}(\xi,\tau')^2 + \hat{H}_{\mathcal{C}}^{\alpha\beta\gamma}(\xi,\tau')^2) \, d\xi \, d\tau' \\ 
&\leq C(\theta) \sup_{\tau \leq \tau_*} \int_{\tau-1}^\tau a^{\alpha\beta\gamma}(\tau')^2 \, d\tau'. 
\end{align*}
\end{proposition}

The proof of Proposition \ref{neutral.mode.dominates} will be given in Section \ref{overlap.region}. \\

Using Proposition \ref{neutral.mode.dominates}, we are able to derive an ODE for the function $a^{\alpha\beta\gamma}(\tau)$:

\begin{proposition}
\label{ode.for.a}
Fix $\theta>0$ small enough so that the conclusion of Proposition \ref{estimate.for.difference.in.tip.region} holds. Let $\delta > 0$ be given. Suppose that $-\tau_*$ is sufficiently large (depending on $\delta$), and the triplet $(\alpha,\beta,\gamma)$ is chosen as in Proposition \ref{choice.of.parameters}. Let $Q^{\alpha\beta\gamma}(\tau) := \frac{d}{d\tau} a^{\alpha\beta\gamma}(\tau) - 2 \, (-\tau)^{-1} \, a^{\alpha\beta\gamma}(\tau)$. Then 
\[\sup_{\tau \leq \tau_*} (-\tau) \int_{\tau-1}^\tau |Q^{\alpha\beta\gamma}(\tau')| \, d\tau' \leq \delta \, \sup_{\tau \leq \tau_*} \bigg ( \int_{\tau-1}^\tau a^{\alpha\beta\gamma}(\tau')^2 \, d\tau' \bigg )^{\frac{1}{2}}.\] 
\end{proposition}

The proof of Proposition \ref{ode.for.a} will be given in Section \ref{analysis.of.neutral.mode}. \\

We now finish the proof of Theorem \ref{uniqueness.theorem}. Using the ODE $\frac{d}{d\tau} a^{\alpha\beta\gamma}(\tau) = 2 \, (-\tau)^{-1} \, a^{\alpha\beta\gamma}(\tau) + Q^{\alpha\beta\gamma}(\tau)$ together with the fact that $a^{\alpha\beta\gamma}(\tau_*) = 0$, we obtain 
\[(-\tau)^2 \, a^{\alpha\beta\gamma}(\tau) = -\int_\tau^{\tau_*} (-\tau')^2 \, Q^{\alpha\beta\gamma}(\tau') \, d\tau'\] 
for all $\tau \leq \tau_*$. This implies 
\begin{align*} 
(-\tau) \, |a^{\alpha\beta\gamma}(\tau)| 
&\leq \int_\tau^{\tau_*} (-\tau') \, |Q^{\alpha\beta\gamma}(\tau')| \, d\tau' \\ 
&\leq \sum_{j=0}^{[\tau_*-\tau]} \int_{\tau_*-j-1}^{\tau_*-j} (-\tau') \, |Q^{\alpha\beta\gamma}(\tau')| \, d\tau' \\ 
&\leq (-\tau) \, \max_{0 \leq j \leq [\tau_*-\tau]} \int_{\tau_*-j-1}^{\tau_*-j} (-\tau') \, |Q^{\alpha\beta\gamma}(\tau')| \, d\tau' 
\end{align*}
for all $\tau \leq \tau_*$. We now divide by $-\tau$, and take the supremum over all $\tau \leq \tau_*$. This implies 
\[\sup_{\tau \leq \tau_*} |a^{\alpha\beta\gamma}(\tau)| \leq \sup_{\tau \leq \tau_*} \int_{\tau-1}^\tau (-\tau') \, |Q^{\alpha\beta\gamma}(\tau')| \, d\tau'.\]
On the other hand, Proposition \ref{ode.for.a} gives the following estimate for $Q^{\alpha\beta\gamma}$: 
\[\sup_{\tau \leq \tau_*} (-\tau) \int_{\tau-1}^\tau |Q^{\alpha\beta\gamma}(\tau')| \, d\tau' \leq \delta \, \sup_{\tau \leq \tau_*} |a^{\alpha\beta\gamma}(\tau)|.\] 
Hence, if we choose $\delta$ sufficiently small, and $-\tau_*$ sufficiently large (depending on $\delta$), then $\sup_{\tau \leq \tau_*} |a^{\alpha\beta\gamma}(\tau)| = 0$. Thus, $a^{\alpha\beta\gamma}(\tau) = 0$ for all $\tau \leq \tau_*$. Proposition \ref{neutral.mode.dominates} then implies $\hat{H}_{\mathcal{C}}^{\alpha\beta\gamma}(\xi,\tau) = 0$ for all $\tau \leq \tau_*$. Putting these facts together, we obtain $H_{\mathcal{C}}^{\alpha\beta\gamma}(\xi,\tau) = 0$ for all $\tau \leq \tau_*$. From this, we deduce that $W_+^{\beta\gamma}(\rho,\tau) = 0$ for $\rho \in [\theta,2\theta]$ and $\tau \leq \tau_*$. Proposition \ref{estimate.for.difference.in.tip.region} yields $W_+^{\beta\gamma}(\rho,\tau) = 0$ for $\rho \in [0,2\theta]$ and $\tau \leq \tau_*$. Thus, we conclude that $F_1(z,t) = F_2^{\alpha\beta\gamma}(z,t)$ for all $t \leq t_* = -e^{-\tau_*}$. In other words, the two ancient solutions coincide for $t \leq t_*$. \\

\section{Energy estimates in the tip region and proof of Proposition \ref{estimate.for.difference.in.tip.region}}

\label{difference.tip.region}

In this section, we give the proof of Proposition \ref{estimate.for.difference.in.tip.region}. Let $\omega_T$ denote a nonnegative smooth cutoff function satisfying $\omega_T(\rho)=1$ for $\rho \leq \theta$ and $\omega_T(\rho) = 0$ for $\rho \geq 2\theta$. We define 
\[W_{T+}^{\beta\gamma}(\rho,\tau) := \omega_T(\rho) \, W_+(\rho,\tau).\]  
To simplify the notation, we will write $W_+$ and $W_{T+}$ instead of $W_+^{\beta\gamma}$ and $W_{T+}^{\beta\gamma}$. \\

\begin{proposition}
\label{pde.for.W}
The function $W_+(\rho,\tau)$ satisfies the equation 
\begin{align*} 
V_{1+}^{-2} \, \Big ( \frac{\partial W_+}{\partial \tau} + \frac{\rho}{2} \, \frac{\partial W_+}{\partial \rho} \Big ) 
&= \frac{\partial^2 W_+}{\partial \rho^2} + \frac{\partial}{\partial \rho} \Big ( \rho^{-1} \, (V_{1+}^{-2} - 1) \, W_+ \Big ) \\ 
&- 2\rho^{-2} \, W_+ + V_{1+}^{-2} \, \mathcal{B}_+ \, W_+, 
\end{align*}
where 
\begin{align*} 
\mathcal{B}_+ 
&:= \rho^{-2} \, \big ( 1 - V_{1+} \, (V_{2+}^{\beta\gamma})^{-1} \big ) \\ 
&+ \rho^{-1} \, \Big ( 2 \, V_{1+}^{-1} \, \frac{\partial V_{1+}}{\partial \rho} - (V_{2+}^{\beta\gamma})^{-2} \, (V_{1+}+V_{2+}^{\beta\gamma}) \, \frac{\partial V_{2+}^{\beta\gamma}}{\partial \rho} \Big ) \\ 
&+ (V_{2+}^{\beta\gamma})^{-2} \, (V_{1+}+V_{2+}^{\beta\gamma}) \, \Big ( \frac{\partial V_{2+}^{\beta\gamma}}{\partial \tau} + \frac{\rho}{2} \, \frac{\partial V_{2+}^{\beta\gamma}}{\partial \rho} \Big ).
\end{align*} 
\end{proposition}

\textbf{Proof.} 
The functions $U_{1+}(r,t)$, $U_{1-}(r,t)$, $U_{2+}^{\beta\gamma}(r,t)$, and $U_{2-}^{\beta\gamma}(r,t)$ all satisfy the same PDE (see e.g. \cite{Brendle2}): 
\[U^{-1} \, \frac{\partial U}{\partial t} = \frac{\partial^2 U}{\partial r^2} - \frac{1}{2} \, U^{-1} \, \Big ( \frac{\partial U}{\partial r} \Big )^2 + r^{-2} \, (U^{-1}-1) \, \Big ( r \, \frac{\partial U}{\partial r}+2U \Big ).\] 
Consequently, the functions $V_{1+}(\rho,\tau)$, $V_{1-}(\rho,\tau)$, $V_{2+}^{\beta\gamma}(\rho,\tau)$, and $V_{2-}^{\beta\gamma}(\rho,\tau)$ satisfy the following PDE: 
\[V^{-2} \, \Big ( \frac{\partial V}{\partial \tau} + \frac{\rho}{2} \, \frac{\partial V}{\partial \rho} \Big ) = \frac{\partial^2 V}{\partial \rho^2} + \rho^{-2} \, (V^{-2}-1) \, \Big ( \rho \, \frac{\partial V}{\partial \rho}+V \Big ).\] 
The assertion now follows from a straightforward calculation. \\

\begin{proposition}
\label{divergence.identity}
The function $W_{T+}(\rho,\tau)$ satisfies 
\begin{align*} 
&\frac{1}{2} \, \frac{\partial}{\partial \tau} \big ( V_{1+}^{-2} \, W_{T+}^2 \, e^{\mu_+} \big ) - \frac{\partial}{\partial \rho} \Big [ \Big ( \frac{\partial W_{T+}}{\partial \rho} + \rho^{-1} \, (V_{1+}^{-2} - 1) \, W_{T+} \Big ) \, W_{T+} \, e^{\mu_+} \Big ] \\ 
&+ \frac{\partial}{\partial \rho} \big ( W_+^2 \, \omega_T' \, \omega_T \, e^{\mu_+} \big ) \\ 
&\leq -\frac{1}{2} \, \Big ( \frac{\partial W_{T+}}{\partial \rho} + \frac{\partial \mu_+}{\partial \rho} \, W_{T+} \Big )^2 \, e^{\mu_+} - 2\rho^{-2} \, W_{T+}^2 \, e^{\mu_+} \\ 
&+ V_{1+}^{-2} \, \Big ( \frac{1}{2} \, \frac{\partial \mu_+}{\partial \tau} - V_{1+}^{-1} \, \frac{\partial V_{1+}}{\partial \tau} + \frac{\rho}{2} \, \frac{\partial \mu_+}{\partial \rho} + \mathcal{B}_+ \Big ) \, W_{T+}^2 \, e^{\mu_+} \\ 
&+ \frac{1}{2} \, \Big ( \frac{\partial \mu_+}{\partial \rho} - \rho^{-1} \, (V_{1+}^{-2} - 1) - \frac{\rho}{2} \, V_{1+}^{-2} \Big )^2 \, W_{T+}^2 \, e^{\mu_+} \\ 
&+ \Big ( \frac{\partial \mu_+}{\partial \rho} - \rho^{-1} \, (V_{1+}^{-2} - 1) + \frac{\rho}{2} \, V_{1+}^{-2} \Big ) \, W_+^2 \, \omega_T' \, \omega_T \, e^{\mu_+} + (\omega_T')^2 \, W_+^2 \, e^{\mu_+}. 
\end{align*} 
\end{proposition}

\textbf{Proof.} 
Using Proposition \ref{pde.for.W}, we obtain  
\begin{align*} 
&V_{1+}^{-2} \, \Big ( \frac{\partial W_{T+}}{\partial \tau} + \frac{\rho}{2} \, \frac{\partial W_{T+}}{\partial \rho} \Big ) \\ 
&= \frac{\partial^2 W_{T+}}{\partial \rho^2} + \frac{\partial}{\partial \rho} \Big ( \rho^{-1} \, (V_{1+}^{-2} - 1) \, W_{T+} \Big ) - 2\rho^{-2} \, W_{T+} + V_{1+}^{-2} \, \mathcal{B}_+ \, W_{T+} \\ 
&+ \Big ( -2 \, \omega_T' \, \frac{\partial W_+}{\partial \rho} - \omega_T'' \, W_+ - \omega_T' \, \rho^{-1} \, (V_{1+}^{-2} - 1) \, W_+ + \frac{\rho}{2} \, \omega_T' \, V_{1+}^{-2} \, W_+ \Big ).
\end{align*}
We next bring in the weight $\mu_+(\rho,\tau)$. A straightforward calculation gives 
\begin{align*} 
&\frac{1}{2} \, \frac{\partial}{\partial \tau} \big ( V_{1+}^{-2} \, W_{T+}^2 \, e^{\mu_+} \big ) - \frac{\partial}{\partial \rho} \Big [ \Big ( \frac{\partial W_{T+}}{\partial \rho} + \rho^{-1} \, (V_{1+}^{-2} - 1) \, W_{T+} \Big ) \, W_{T+} \, e^{\mu_+} \Big ] \\ 
&+ \frac{\partial}{\partial \rho} \big ( W_+^2 \, \omega_T' \, \omega_T \, e^{\mu_+} \big ) \\ 
&= -\Big ( \frac{\partial W_{T+}}{\partial \rho} + \frac{\partial \mu_+}{\partial \rho} \, W_{T+} \Big )^2 \, e^{\mu_+} - 2\rho^{-2} \, W_{T+}^2 \, e^{\mu_+} \\ 
&+ V_{1+}^{-2} \, \Big ( \frac{1}{2} \, \frac{\partial \mu_+}{\partial \tau} - V_{1+}^{-1} \, \frac{\partial V_{1+}}{\partial \tau} + \frac{\rho}{2} \, \frac{\partial \mu_+}{\partial \rho} + \mathcal{B}_+ \Big ) \, W_{T+}^2 \, e^{\mu_+} \\ 
&+ \Big ( \frac{\partial \mu_+}{\partial \rho} - \rho^{-1} \, (V_{1+}^{-2} - 1) - \frac{\rho}{2} \, V_{1+}^{-2} \Big ) \Big ( \frac{\partial W_{T+}}{\partial \rho} + \frac{\partial \mu_+}{\partial \rho} \, W_{T+} \Big ) \, W_{T+} \, e^{\mu_+} \\ 
&+ \Big ( \frac{\partial \mu_+}{\partial \rho} - \rho^{-1} \, (V_{1+}^{-2} - 1) + \frac{\rho}{2} \, V_{1+}^{-2} \Big ) \, W_+^2 \, \omega_T' \, \omega_T \, e^{\mu_+} + (\omega_T')^2 \, W_+^2 \, e^{\mu_+}. 
\end{align*} 
The assertion follows now from Young's inequality. \\

\begin{corollary}
\label{divergence.identity.2}
Fix a small number $\eta > 0$. Then we can find a small number $\theta \in (0,\eta)$ and a small number $\varepsilon \in (0,\eta)$ (both depending on $\eta$) with the following property. If $-\tau_*$ sufficiently large (depending on $\eta$ and $\theta$) and the triplet $(\alpha,\beta,\gamma)$ is $\varepsilon$-admissible with respect to time $t_* = -e^{-\tau_*}$, then we have 
\begin{align*} 
&\frac{1}{2} \, \frac{\partial}{\partial \tau} \big ( V_{1+}^{-2} \, W_{T+}^2 \, e^{\mu_+} \big ) - \frac{\partial}{\partial \rho} \Big [ \Big ( \frac{\partial W_{T+}}{\partial \rho} + \rho^{-1} \, (V_{1+}^{-2} - 1) \, W_{T+} \Big ) \, W_{T+} \, e^{\mu_+} \Big ] \\ 
&+ \frac{\partial}{\partial \rho} \big ( W_+^2 \, \omega_T' \, \omega_T \, e^{\mu_+} \big ) \\ 
&\leq -\frac{1}{2} \, \Big ( \frac{\partial W_{T+}}{\partial \rho} + \frac{\partial \mu_+}{\partial \rho} \, W_{T+} \Big )^2 \, e^{\mu_+} - 2\rho^{-2} \, W_{T+}^2 \, e^{\mu_+} \\ 
&+ \eta \, \rho^{-2} \, V_{1+}^{-4} \, W_{T+}^2 \, e^{\mu_+} + \eta \, \rho^{-2} \, V_{1+}^{-2} \, W_+^2 \, e^{\mu_+} \, 1_{\{\theta \leq \rho \leq 2\theta\}}
\end{align*} 
for $\rho \leq 2\theta$ and $\tau \leq \tau_*$.
\end{corollary}

\textbf{Proof.}
By Proposition \ref{comparison.of.V.with.bryant.soliton.profile}, Proposition \ref{higher.derivatives.of.V.in.collar.region}, and Proposition \ref{comparison.of.modified.V.with.bryant.soliton.profile}, we can choose $\theta \in (0,\eta)$ (depending on $\eta$) sufficiently small and $-\tau_*$ sufficiently large (depending on $\eta$ and $\theta$) such that 
\[|\mathcal{B}_+| \leq \eta \, \rho^{-2} \, V_{1+}^{-2}\] 
for $\rho \leq 2\theta$ and $\tau \leq \tau_*$. By Corollary \ref{derivative.of.V.wrt.tau}, Lemma \ref{derivative.of.mu.wrt.rho}, and Lemma \ref{derivative.of.mu.wrt.tau}, we can choose $\theta \in (0,\eta)$ sufficiently small (depending on $\eta$) and $-\tau_*$ sufficiently large (depending on $\eta$ and $\theta$) such that 
\begin{align*} 
&\Big | \frac{1}{2} \, \frac{\partial \mu_+}{\partial \tau} - V_{1+}^{-1} \, \frac{\partial V_{1+}}{\partial \tau} + \frac{\rho}{2} \, \frac{\partial \mu_+}{\partial \rho} \Big | \leq \eta \, \rho^{-2} \, V_{1+}^{-2}, \\ 
&\Big | \frac{\partial \mu_+}{\partial \rho} - \rho^{-1} \, (V_{1+}^{-2} - 1) - \frac{\rho}{2} \, V_{1+}^{-2} \Big | \leq \eta \, \rho^{-1} \, V_{1+}^{-2}, \\ 
&\Big | \frac{\partial \mu_+}{\partial \rho} - \rho^{-1} \, (V_{1+}^{-2} - 1) + \frac{\rho}{2} \, V_{1+}^{-2} \Big | \leq \eta \, \rho^{-1} \, V_{1+}^{-2} 
\end{align*} 
for $\rho \leq 2\theta$ and $\tau \leq \tau_*$. Hence, the assertion follows from Proposition \ref{divergence.identity}. \\

In the next step, we finalize our choice of $\theta$.

\begin{proposition}
\label{integral.estimate}
We can find small numbers $\theta>0$, $\lambda>0$, and $\varepsilon>0$ with the following property. If $-\tau_*$ is sufficiently large (depending on $\theta$) and the triplet $(\alpha,\beta,\gamma)$ is $\varepsilon$-admissible with respect to time $t_* = -e^{-\tau_*}$, then 
\begin{align*} 
\frac{1}{2} \, \frac{d}{d\tau} \bigg ( \int_0^{2\theta} V_{1+}^{-2} \, W_{T+}^2 \, e^{\mu_+} \, d\rho \bigg ) 
&\leq -\lambda \, (-\tau) \int_0^{2\theta} V_{1+}^{-2} \, W_{T+}^2 \, e^{\mu_+} \, d\rho \\ 
&+ \int_\theta^{2\theta} \rho^{-2} \, V_{1+}^{-2} \, W_+^2 \, e^{\mu_+} \, d\rho 
\end{align*}
for $\tau \leq \tau_*$.
\end{proposition}

\textbf{Proof.} 
Let us fix a small number $\eta > 0$. In the following, we choose $\theta$ and $\varepsilon$ sufficiently small (depending on $\eta$), and we choose $-\tau_*$ sufficiently large (depending on $\eta$ and $\theta$). Using Corollary \ref{divergence.identity.2}, we obtain 
\begin{align*} 
&\frac{1}{2} \, \frac{d}{d\tau} \bigg ( \int_0^{2\theta} V_{1+}^{-2} \, W_{T+}^2 \, e^{\mu_+} \, d\rho \bigg ) \\ 
&\leq -\frac{1}{2} \int_0^{2\theta} \Big ( \frac{\partial W_{T+}}{\partial \rho} + \frac{\partial \mu_+}{\partial \rho} \, W_{T+} \Big )^2 \, e^{\mu_+} \, d\rho - 2 \int_0^{2\theta} \rho^{-2} \, W_{T+}^2 \, e^{\mu_+} \, d\rho \\ 
&+ \eta \int_0^{2\theta} \rho^{-2} \, V_{1+}^{-4} \, W_{T+}^2 \, e^{\mu_+} \, d\rho + \eta \int_\theta^{2\theta} \rho^{-2} \, V_{1+}^{-2} \, W_+^2 \, e^{\mu_+} \, d\rho 
\end{align*} 
for $\tau \leq \tau_*$. Applying Proposition \ref{Poincare.inequality} to the function $f := e^{\mu_+} \, W_{T+}$ gives 
\begin{align*} 
0 &\leq 8 \int_0^{2\theta} \Big ( \frac{\partial W_{T+}}{\partial \rho} + \frac{\partial \mu_+}{\partial \rho} \, W_{T+} \Big )^2  \, e^{\mu_+} \, d\rho \\ 
&+ K_* \int_0^{2\theta} \rho^{-2} \, W_{T+}^2 \, e^{\mu_+} \, d\rho - \int_0^{2\theta} \Big ( \frac{\partial \mu_+}{\partial \rho} \Big )^2 \, W_{T+}^2 \, e^{\mu_+} \, d\rho 
\end{align*}
for $\tau \leq \tau_*$. Using Lemma \ref{derivative.of.mu.wrt.rho}, we obtain $(\frac{\partial \mu_+}{\partial \rho})^2 \geq \frac{1}{4} \, \rho^{-2} \, (V_{1+}^{-2}-1)^2$ for $\rho \leq 2\theta$, hence 
\begin{align*} 
0 &\leq 128\eta \int_0^{2\theta} \Big ( \frac{\partial W_{T+}}{\partial \rho} + \frac{\partial \mu_+}{\partial \rho} \, W_{T+} \Big )^2  \, e^{\mu_+} \, d\rho \\ 
&+ 16\eta K_* \int_0^{2\theta} \rho^{-2} \, W_{T+}^2 \, e^{\mu_+} \, d\rho - 4\eta \int_0^{2\theta} \rho^{-2} \, (V_{1+}^{-2}-1)^2 \, W_{T+}^2 \, e^{\mu_+} \, d\rho 
\end{align*}
for $\tau \leq \tau_*$. Adding the two inequalities, we conclude that 
\begin{align*} 
&\frac{1}{2} \, \frac{d}{d\tau} \bigg ( \int_0^{2\theta} V_{1+}^{-2} \, W_{T+}^2 \, e^{\mu_+} \, d\rho \bigg ) \\ 
&\leq -\Big ( \frac{1}{2}-128\eta \Big ) \int_0^{2\theta} \Big ( \frac{\partial W_{T+}}{\partial \rho} + \frac{\partial \mu_+}{\partial \rho} \, W_{T+} \Big )^2 \, e^{\mu_+} \, d\rho \\ 
&- (2-4\eta-16\eta K_*) \int_0^{2\theta} \rho^{-2} \, W_{T+}^2 \, e^{\mu_+} \, d\rho \\ 
&- \eta \int_0^{2\theta} \rho^{-2} \, [4 \, (V_{1+}^{-2}-1)^2 + 4 - V_{1+}^{-4}] \, W_{T+}^2 \, e^{\mu_+} \, d\rho \\ 
&+ \eta \int_\theta^{2\theta} \rho^{-2} \, V_{1+}^{-2} \, W_+^2 \, e^{\mu_+} \, d\rho 
\end{align*} 
for $\tau \leq \tau_*$. We assume that $\eta>0$ is chosen small enough so that $\frac{1}{2}-128\eta > 0$ and $2-4\eta-16\eta K_* > 0$. (Here, it is crucial that the constant $K_*$ in the weighted Poincar\'e inequality does not depend on $\theta$.) This ensures that the first two terms on the right hand side have a favorable sign. To estimate the third term on the right hand side, we observe that $\rho^{-2} \, [4 \, (V_{1+}^{-2}-1)^2 + 4 - V_{1+}^{-4}] \geq \rho^{-2} \, V_{1+}^{-4}$. In view of Proposition \ref{comparison.of.V.with.bryant.soliton.profile}, the term $\rho^{-2} \, V_{1+}^{-4}$ is bounded from below by a small positive multiple of $(-\tau) \, V_{1+}^{-2}$. This completes the proof of Proposition \ref{integral.estimate}. \\

We now complete the proof of Proposition \ref{estimate.for.difference.in.tip.region}. Let $\theta$, $\lambda$, and $\varepsilon$ be chosen as in Proposition \ref{integral.estimate}. Let 
\[I(\tau) := \int_{\tau-1}^\tau \int_0^{2\theta} V_{1+}^{-2} \, W_{T+}^2 \, e^{\mu_+}\] 
and 
\[J(\tau) := \int_{\tau-1}^\tau \int_\theta^{2\theta} V_{1+}^{-2} \, W_+^2 \, e^{\mu_+}.\] 
If we choose $-\tau_*$ sufficiently large, then Proposition \ref{integral.estimate} gives 
\[\frac{1}{2} \, I'(\tau) + \lambda \, (-\tau) \, I(\tau) \leq \theta^{-2} \, J(\tau),\] 
hence 
\[\frac{d}{d\tau} (e^{-\lambda \tau^2} \, I(\tau)) \leq 2\theta^{-2} \, e^{-\lambda \tau^2} \, J(\tau)\] 
for $\tau \leq \tau_*$. Clearly, $\lim_{\tau \to -\infty} e^{-\lambda \tau^2} \, I(\tau) = 0$. Consequently, 
\begin{align*} 
e^{-\lambda \tau^2} \, I(\tau) 
&\leq 2\theta^{-2} \int_{-\infty}^\tau e^{-\lambda {\tau'}^2} \, J(\tau') \, d\tau' \\ 
&\leq 2\theta^{-2} \, \Big ( \sup_{\tau' \leq \tau} (-\tau')^{-1} \, J(\tau') \Big ) \, \int_{-\infty}^\tau e^{-\lambda {\tau'}^2} \, (-\tau') \, d\tau' \\ 
&\leq \theta^{-2} \lambda^{-1} \, e^{-\lambda \tau^2} \, \sup_{\tau' \leq \tau} (-\tau')^{-1} \, J(\tau') 
\end{align*} 
for $\tau \leq \tau_*$. This finally gives 
\begin{align*} 
(-\tau)^{-\frac{1}{2}} \, I(\tau) 
&\leq \theta^{-2} \lambda^{-1} \, (-\tau)^{-\frac{1}{2}} \sup_{\tau' \leq \tau} (-\tau')^{-1} \, J(\tau') \\ 
&\leq \theta^{-2} \lambda^{-1} \, (-\tau)^{-1} \sup_{\tau' \leq \tau} (-\tau')^{-\frac{1}{2}} \, J(\tau')  
\end{align*}
for $\tau \leq \tau_*$. Taking the supremum over $\tau \leq \tau_*$ gives 
\[\sup_{\tau \leq \tau_*} (-\tau)^{-\frac{1}{2}} \, I(\tau) \leq \theta^{-2} \lambda^{-1} \, (-\tau_*)^{-1} \, \sup_{\tau \leq \tau_*} (-\tau)^{-\frac{1}{2}} \, J(\tau).\] 
From this, the conclusion of Proposition \ref{estimate.for.difference.in.tip.region} follows immediately. \\

\section{Energy estimates in the cylindrical region and proof of Proposition \ref{estimate.for.difference.in.cylindrical.region}}

\label{difference.cylindrical.region}

In this section, we give the proof of Proposition \ref{estimate.for.difference.in.cylindrical.region}. Throughout this section, we assume that $\theta$ is chosen as in Proposition \ref{estimate.for.difference.in.tip.region}. To simplify the notation, we will write $H$, $H_{\mathcal{C}}$, $\hat{H}_{\mathcal{C}}$, and $a$ instead of $H^{\alpha\beta\gamma}$, $H_{\mathcal{C}}^{\alpha\beta\gamma}$, $\hat{H}_{\mathcal{C}}^{\alpha\beta\gamma}$, and $a^{\alpha\beta\gamma}$. 

Our goal is to study the evolution equation satisfied by the function $H$. The linearized operator 
\[\mathcal{L} f := f_{\xi\xi} - \frac{1}{2} \, \xi \, f_\xi + f\] 
is the same as in \cite{Angenent-Daskalopoulos-Sesum2}, and hence the linear theory from \cite{Angenent-Daskalopoulos-Sesum2} carries over to the Ricci flow case as well. In order for this article to be self-contained, we will state the results from \cite{Angenent-Daskalopoulos-Sesum2} that we will use later, but for the proofs of the same we refer the reder to \cite{Angenent-Daskalopoulos-Sesum2}. 

As in \cite{Angenent-Daskalopoulos-Sesum2}, we consider the Hilbert space $\mathcal{H} = L^2(\mathbb{R},e^{-\frac{\xi^2}{4}} \, d\xi)$. The norm on $\mathcal{H}$ is given by 
\[\|f\|_{\mathcal{H}}^2 := \int_{\mathbb{R}} e^{-\frac{\xi^2}{4}} \, f(\xi)^2 \, d\xi.\] 
Moreover, we denote by $\mathcal{D} \subset \mathcal{H}$ the Hilbert space of all functions $f$ such that $f \in \mathcal{H}$ and $f' \in \mathcal{H}$. The norm on $\mathcal{D}$ is given by 
\[\|f\|_{\mathcal{D}}^2 := \int_{\mathbb{R}} e^{-\frac{\xi^2}{4}} \, (f'(\xi)^2+f(\xi)^2) \, d\xi.\] 
Let $\mathcal{D}^*$ denote the dual space of $\mathcal{D}$. Clearly, the dual space $\mathcal{H}^*$ is a subspace of $\mathcal{D}^*$. After identifying $\mathcal{H}^*$ with $\mathcal{H}$ in the standard way, we can view $\mathcal{H}$ as a subspace of $\mathcal{D}^*$. The restriction of $\|\cdot\|_{\mathcal{D}^*}$ to $\mathcal{H}$ is given by 
\[\|f\|_{\mathcal{D}^*} := \sup \bigg \{ \int_{\mathbb{R}} e^{-\frac{\xi^2}{4}} \, f(\xi) \, g(\xi) \, d\xi: \|g\|_{\mathcal{D}} \leq 1 \bigg \}\] 
for $f \in \mathcal{H}$. For later reference, we collect some basic facts from \cite{Angenent-Daskalopoulos-Sesum2}. \\

\begin{proposition}
\label{boundedness.of.operators}
The following statements hold:
\begin{itemize}
\item[(i)] The operators $f \mapsto \xi \, f$, $f \mapsto f'$, $f \mapsto -f' + \frac{1}{2} \, \xi \, f$ are bounded from $\mathcal{D}$ to $\mathcal{H}$. 
\item[(ii)] The operators $f \mapsto \xi \, f$, $f \mapsto f'$, $f \mapsto -f' + \frac{1}{2} \, \xi \, f$ are bounded from $\mathcal{H}$ to $\mathcal{D}^*$. 
\item[(iii)] The operators $f \mapsto \xi^2 \, f$, $f \mapsto \xi \, f'$, $f \mapsto f''$ are bounded from $\mathcal{D}$ to $\mathcal{D}^*$.
\item[(iv)] The operator $f \mapsto \int_0^\xi f$ is bounded from $\mathcal{H}$ to $\mathcal{D}$.
\end{itemize}
\end{proposition}

\textbf{Proof.} 
Statements (i), (ii), and (iii) were proved in \cite{Angenent-Daskalopoulos-Sesum2}. To prove statement (iv), let us consider a function $f \in \mathcal{H}$, and let $g(\xi) := \int_0^\xi f(\xi') \, d\xi'$. Then $g(\xi)^2 \leq \xi \int_0^\xi f(\xi')^2 \, d\xi'$ for $\xi \geq 0$. Using Fubini's theorem, we obtain 
\begin{align*} 
\int_0^\infty e^{-\frac{\xi^2}{4}} \, g(\xi)^2 \, d\xi 
&\leq \int_0^\infty e^{-\frac{\xi^2}{4}} \, \xi \, \bigg ( \int_0^\xi f(\xi')^2 \, d\xi' \bigg ) \, d\xi \\ 
&= \int_0^\infty \bigg ( \int_{\xi'}^\infty e^{-\frac{\xi^2}{4}} \, \xi \, d\xi \bigg ) \, f(\xi')^2 \, d\xi' \\ 
&= 2 \int_0^\infty e^{-\frac{\xi'^2}{4}} \, f(\xi')^2 \, d\xi'. 
\end{align*}
An analogous argument gives $\int_{-\infty}^0 e^{-\frac{\xi^2}{4}} \, g(\xi)^2 \, d\xi \leq 2 \int_{-\infty}^0 e^{-\frac{\xi'^2}{4}} \, f(\xi')^2 \, d\xi'$. Therefore, $\|g\|_{\mathcal{H}} \leq C \, \|f\|_{\mathcal{H}}$. Since $g'=f$, it follows that $\|g\|_{\mathcal{D}} \leq C \, \|f\|_{\mathcal{H}}$, as claimed. \\

For a time-dependent function $f$, we introduce the following norms:
\begin{align*} 
&\|f\|_{\mathcal{H},\infty,\tau_*}^2 := \sup_{\tau \leq \tau_*} \int_{\tau-1}^\tau \|f(\cdot,\tau')\|_{\mathcal{H}}^2 \, d\tau', \\ 
&\|f\|_{\mathcal{D},\infty,\tau_*}^2 := \sup_{\tau \leq \tau_*} \int_{\tau-1}^\tau \|f(\cdot,\tau')\|_{\mathcal{D}}^2 \, d\tau', \\ 
&\|f\|_{\mathcal{D}^*,\infty,\tau_*}^2 := \sup_{\tau \leq \tau_*} \int_{\tau-1}^\tau \|f(\cdot,\tau')\|_{\mathcal{D}^*}^2 \, d\tau'.
\end{align*}
The following energy estimate was proved in \cite{Angenent-Daskalopoulos-Sesum2}: \\


\begin{proposition}
\label{linear.estimate}
Let $g: (-\infty,\tau_*] \to \mathcal{D}^*$ be a bounded function. Let $f: (-\infty,\tau_*] \to \mathcal{D}$ be a bounded function which satisfies the linear equation 
\[\frac{\partial}{\partial \tau} f(\tau) - \mathcal{L} f(\tau) = g(\tau).\] 
Then the function $\hat{f} := P_+ f + P_- f$ satisfies the estimate 
\[\sup_{\tau \leq \tau_*} \|\hat{f}(\tau)\|_{\mathcal{H}} + \Lambda^{-1} \, \|\hat{f}\|_{\mathcal{D},\infty,\tau_*} \leq \|P_+ f(\tau_*)\|_{\mathcal{H}} + \Lambda \, \|g\|_{\mathcal{D}^*,\infty,\tau_*},\] 
where $\Lambda$ is a universal constant.
\end{proposition}

\textbf{Proof.} 
See \cite{Angenent-Daskalopoulos-Sesum2}, Lemma 6.6. \\

We now continue with the proof of Proposition \ref{estimate.for.difference.in.cylindrical.region}. The functions $G_1(\xi,\tau)$ and $G_2^{\alpha\beta\gamma}(\xi,\tau)$ satisfy the equation 
\begin{align*} 
G_\tau(\xi,\tau) 
&= G_{\xi\xi}(\xi,\tau) - \frac{1}{2} \, \xi \, G_\xi(\xi,\tau) \\
&+ \frac{1}{2} \, (\sqrt{2}+G(\xi,\tau)) - (\sqrt{2}+G(\xi,\tau))^{-1} \\ 
&-(\sqrt{2}+G(\xi,\tau))^{-1} \, G_\xi(\xi,\tau)^2 \\ 
&+ 2 \, G_\xi(\xi,\tau) \, \bigg [ \frac{G_\xi(0,\tau)}{\sqrt{2}+G(0,\tau)} - \int_0^\xi \frac{G_\xi(\xi',\tau)^2}{(\sqrt{2}+G(\xi',\tau))^2} \, d\xi' \bigg ].
\end{align*} 
Consequently, the difference $H(\xi,\tau) = G_1(\xi,\tau) - G_2^{\alpha\beta\gamma}(\xi,\tau)$ satisfies 
\[H_\tau(\xi,\tau) = H_{\xi\xi}(\xi,\tau) - \frac{1}{2} \, \xi \, H_\xi(\xi,\tau) + H(\xi,\tau) + \sum_{k=1}^6 E_k(\xi,\tau),\] 
where 
\begin{align*} 
E_1(\xi,\tau) &= \Big [ (\sqrt{2}+G_1(\xi,\tau))^{-1} (\sqrt{2}+G_2^{\alpha\beta\gamma}(\xi,\tau))^{-1} - \frac{1}{2} \Big ] \, H(\xi,\tau) \\ 
E_2(\xi,\tau) &= (\sqrt{2}+G_1(\xi,\tau))^{-1} (\sqrt{2}+G_2^{\alpha\beta\gamma}(\xi,\tau))^{-1} \, G_{1\xi}(\xi,\tau)^2 \, H(\xi,\tau), \\ 
E_3(\xi,\tau) &= -(\sqrt{2}+G_2^{\alpha\beta\gamma}(\xi,\tau))^{-1} \, (G_{1\xi}(\xi,\tau) + G_{2\xi}^{\alpha\beta\gamma}(\xi,\tau)) \, H_\xi(\xi,\tau) \\ 
E_4(\xi,\tau) &= 2 \, \bigg [ \frac{G_{1\xi}(0,\tau)}{\sqrt{2}+G_1(0,\tau)} - \int_0^\xi \frac{G_{1\xi}(\xi',\tau)^2}{(\sqrt{2}+G_1(\xi',\tau))^2} \, d\xi' \bigg ] \, H_\xi(\xi,\tau), \\ 
E_5(\xi,\tau) &= 2 \, G_{2\xi}^{\alpha\beta\gamma}(\xi,\tau) \, \frac{H_\xi(0,\tau)}{\sqrt{2}+G_1(0,\tau)} \\ &- 2 \, G_{2\xi}^{\alpha\beta\gamma}(\xi,\tau) \, \frac{G_{2\xi}^{\alpha\beta\gamma}(0,\tau) \, H(0,\tau)}{(\sqrt{2}+G_1(0)) (\sqrt{2}+G_2^{\alpha\beta\gamma}(0,\tau))}, \\ 
E_6(\xi,\tau) &= 2 \, G_{2\xi}^{\alpha\beta\gamma}(\xi,\tau) \, \bigg [ -\int_0^\xi \frac{(G_{1\xi}(\xi',\tau)+G_{2\xi}^{\alpha\beta\gamma}(\xi',\tau)) \, H_\xi(\xi',\tau)}{(\sqrt{2}+G_2^{\alpha\beta\gamma}(\xi',\tau))^2} \, d\xi' \\ 
&+ \int_0^\xi \frac{(2\sqrt{2}+G_1(\xi',\tau)+G_2^{\alpha\beta\gamma}(\xi',\tau)) \, H(\xi',\tau) \, G_{1\xi}(\xi',\tau)^2}{(\sqrt{2}+G_1(\xi',\tau))^2 (\sqrt{2}+G_2^{\alpha\beta\gamma}(\xi',\tau))^2} \, d\xi' \bigg ]. 
\end{align*}
Consequently, the function $H_{\mathcal{C}}(\xi,\tau) = \chi_{\mathcal{C}}((-\tau)^{-\frac{1}{2}} \xi) \, H(\xi,\tau)$ satisfies 
\[H_{\mathcal{C},\tau}(\xi,\tau) = H_{\mathcal{C},\xi\xi}(\xi,\tau) - \frac{1}{2} \, \xi \, H_{\mathcal{C},\xi}(\xi,\tau) + H_{\mathcal{C}}(\xi,\tau) + \sum_{k=1}^{10} E_{\mathcal{C},k}(\xi,\tau),\] 
where 
\begin{align*} 
E_{\mathcal{C},1}(\xi,\tau) &= \Big [ (\sqrt{2}+G_1(\xi,\tau))^{-1} (\sqrt{2}+G_2^{\alpha\beta\gamma}(\xi,\tau))^{-1} - \frac{1}{2} \Big ] \, H_{\mathcal{C}}(\xi,\tau), \\ 
E_{\mathcal{C},2}(\xi,\tau) &= (\sqrt{2}+G_1(\xi,\tau))^{-1} (\sqrt{2}+G_2^{\alpha\beta\gamma}(\xi,\tau))^{-1} \, G_{1\xi}(\xi,\tau)^2 \, H_{\mathcal{C}}(\xi,\tau), \\ 
E_{\mathcal{C},3}(\xi,\tau) &= -(\sqrt{2}+G_2^{\alpha\beta\gamma}(\xi,\tau))^{-1} \, (G_{1\xi}(\xi,\tau) + G_{2\xi}^{\alpha\beta\gamma}(\xi,\tau)) \, H_{\mathcal{C},\xi}(\xi,\tau), \\ 
E_{\mathcal{C},4}(\xi,\tau) &= 2 \, \bigg [ \frac{G_{1\xi}(0,\tau)}{\sqrt{2}+G_1(0,\tau)} - \int_0^\xi \frac{G_{1\xi}(\xi',\tau)^2}{(\sqrt{2}+G_1(\xi',\tau))^2} \, d\xi' \bigg ] \, H_{\mathcal{C},\xi}(\xi,\tau), \\ 
E_{\mathcal{C},5}(\xi,\tau) &= 2 \, \chi_{\mathcal{C}}((-\tau)^{-\frac{1}{2}} \xi) \, G_{2\xi}^{\alpha\beta\gamma}(\xi,\tau) \, \frac{H_\xi(0,\tau)}{\sqrt{2}+G_1(0,\tau)} \\ &- 2 \, \chi_{\mathcal{C}}((-\tau)^{-\frac{1}{2}} \xi) \, G_{2\xi}^{\alpha\beta\gamma}(\xi,\tau) \, \frac{G_{2\xi}^{\alpha\beta\gamma}(0,\tau) \, H(0,\tau)}{(\sqrt{2}+G_1(0)) (\sqrt{2}+G_2^{\alpha\beta\gamma}(0,\tau))}, \\ 
E_{\mathcal{C},6}(\xi,\tau) &= 2 \, \chi_{\mathcal{C}}((-\tau)^{-\frac{1}{2}} \xi) \, G_{2\xi}^{\alpha\beta\gamma}(\xi,\tau) \\ 
&\cdot \bigg [ -\int_0^\xi \frac{(G_{1\xi}(\xi',\tau)+G_{2\xi}^{\alpha\beta\gamma}(\xi',\tau)) \, H_\xi(\xi',\tau)}{(\sqrt{2}+G_2^{\alpha\beta\gamma}(\xi',\tau))^2} \, d\xi' \\ 
&+ \int_0^\xi \frac{(2\sqrt{2}+G_1(\xi',\tau)+G_2^{\alpha\beta\gamma}(\xi',\tau)) \, H(\xi',\tau) \, G_{1\xi}(\xi',\tau)^2}{(\sqrt{2}+G_1(\xi',\tau))^2 (\sqrt{2}+G_2^{\alpha\beta\gamma}(\xi',\tau))^2} \, d\xi' \bigg ], \\ 
E_{\mathcal{C},7}(\xi,\tau) &= (\sqrt{2}+G_2^{\alpha\beta\gamma}(\xi,\tau))^{-1} \, (G_{1\xi}(\xi,\tau) + G_{2\xi}^{\alpha\beta\gamma}(\xi,\tau)) \\ 
&\cdot (-\tau)^{-\frac{1}{2}} \, \chi_{\mathcal{C}}'((-\tau)^{-\frac{1}{2}} \xi) \, H(\xi,\tau), \\ 
E_{\mathcal{C},8}(\xi,\tau) &= -2 \, \bigg [ \frac{G_{1\xi}(0,\tau)}{\sqrt{2}+G_1(0,\tau)} - \int_0^\xi \frac{G_{1\xi}(\xi',\tau)^2}{(\sqrt{2}+G_1(\xi',\tau))^2} \, d\xi' \bigg ] \\ 
&\cdot (-\tau)^{-\frac{1}{2}} \, \chi_{\mathcal{C}}'((-\tau)^{-\frac{1}{2}} \xi) \, H(\xi,\tau), \\ 
E_{\mathcal{C},9}(\xi,\tau) &= (-\tau)^{-1} \, \chi_{\mathcal{C}}''((-\tau)^{-\frac{1}{2}} \xi) \, H(\xi,\tau) \\ 
&+ \frac{1}{2} \, (-\tau)^{-\frac{3}{2}} \xi \, \chi_{\mathcal{C}}'((-\tau)^{-\frac{1}{2}} \xi) \, H(\xi,\tau), \\ 
E_{\mathcal{C},10}(\xi,\tau) &= -2 \, (-\tau)^{-\frac{1}{2}} \, \frac{\partial}{\partial \xi} \big [ \chi_{\mathcal{C}}'((-\tau)^{-\frac{1}{2}} \xi) \, H(\xi,\tau) \big ] \\ 
&+ \frac{1}{2} \, (-\tau)^{-\frac{1}{2}} \xi \, \chi_{\mathcal{C}}'((-\tau)^{-\frac{1}{2}} \xi) \, H(\xi,\tau).
\end{align*}
In the following, we will estimate the terms $\sum_{k=1}^6 \|E_{\mathcal{C},k}\|_{\mathcal{H},\infty,\tau_*}$ and $\sum_{k=7}^{10} \|E_{\mathcal{C},k}\|_{\mathcal{D}^*,\infty,\tau_*}$. To that end, we need the following estimates for the functions $G_1(\xi,\tau)$ and $G_2^{\alpha\beta\gamma}(\xi,\tau)$:

\begin{proposition} 
\label{estimates.for.G_1.and.G_2}
Fix a small number $\theta > 0$ and a small number $\eta > 0$. Then there exists a small number $\varepsilon>0$ (depending on $\theta$ and $\eta$) with the following property. If the triplet $(\alpha,\beta,\gamma)$ is $\varepsilon$-admissible with respect to time $t_* = -e^{-\tau_*}$ and $-\tau_*$ is sufficiently large, then  
\begin{align*} 
&\Big | (\sqrt{2}+G_1(\xi,\tau))^2 - 2 + \frac{\xi^2-2}{(-2\tau)} \Big | \leq \eta \, \frac{\xi^2+1}{(-\tau)}, \\ 
&\Big | (\sqrt{2}+G_2^{\alpha\beta\gamma}(\xi,\tau))^2 - 2 + \frac{\xi^2-2}{(-2\tau)} \Big | \leq \eta \, \frac{\xi^2+1}{(-\tau)} 
\end{align*} 
and 
\begin{align*} 
&\Big | (\sqrt{2}+G_1(\xi,\tau)) \, G_{1\xi}(\xi,\tau) + \frac{\xi}{(-2\tau)} \Big | \leq \eta \, \frac{|\xi|+1}{(-\tau)}, \\ 
&\Big | (\sqrt{2}+G_2^{\alpha\beta\gamma}(\xi,\tau)) \, G_{2\xi}^{\alpha\beta\gamma}(\xi,\tau) + \frac{\xi}{(-2\tau)} \Big | \leq \eta \, \frac{|\xi|+1}{(-\tau)} 
\end{align*} 
for $|\xi| \leq \sqrt{4-\frac{\theta^2}{8}} \, (-\tau)^{\frac{1}{2}}$ and $\tau \leq \tau_*$. 
\end{proposition}

\textbf{Proof.} 
This follows directly from Proposition \ref{precise.estimate.for.F}, Proposition \ref{precise.estimate.for.F_z}, and Proposition \ref{estimate.for.modified.profile}. \\

In order to estimate the term $\|E_{\mathcal{C},6}\|_{\mathcal{H},\infty,\tau_*}$, we need the following pointwise estimate: 

\begin{lemma}
\label{pointwise.estimate.for.E6}
We have 
\begin{align*} 
|E_{\mathcal{C},6}(\xi,\tau)| 
&\leq C(\theta) \, (-\tau)^{-\frac{1}{2}} \, |G_{2\xi}^{\alpha\beta\gamma}(\xi,\tau)| \, \bigg | \int_0^\xi |H_{\mathcal{C}}(\xi',\tau)| \, d\xi' \bigg | \\ 
&+ C(\theta) \, (-\tau)^{-\frac{1}{2}} \, |G_{2\xi}^{\alpha\beta\gamma}(\xi,\tau)| \, (|H_{\mathcal{C}}(\xi,\tau)|+|H(0,\tau)|)
\end{align*} 
for all $\tau \leq \tau_*$.
\end{lemma}

\textbf{Proof.} 
Proposition \ref{higher.derivative.bounds.in.cylindrical.region} implies 
\[|G_{1\xi\xi}(\xi,\tau)| + |G_{2\xi\xi}^{\alpha\beta\gamma}(\xi,\tau)| \leq C(\theta) \, (-\tau)^{-\frac{1}{2}}\] 
for $|\xi| \leq \sqrt{4-\frac{\theta^2}{4}} \, (-\tau)^{\frac{1}{2}}$. Using integration by parts, we obtain 
\begin{align*} 
&\bigg | \int_0^\xi \frac{(G_{1\xi}(\xi',\tau)+G_{2\xi}^{\alpha\beta\gamma}(\xi',\tau)) \, H_\xi(\xi',\tau)}{(\sqrt{2}+G_2^{\alpha\beta\gamma}(\xi',\tau))^2} \, d\xi' \\ 
&+ \int_0^\xi \frac{(G_{1\xi\xi}(\xi',\tau)+G_{2\xi\xi}^{\alpha\beta\gamma}(\xi',\tau)) \, H(\xi',\tau)}{(\sqrt{2}+G_2^{\alpha\beta\gamma}(\xi',\tau))^2} \, d\xi' \\ 
&- 2 \int_0^\xi \frac{G_{2\xi}^{\alpha\beta\gamma}(\xi',\tau) \, (G_{1\xi}(\xi',\tau)+G_{2\xi}^{\alpha\beta\gamma}(\xi',\tau)) \, H(\xi',\tau)}{(\sqrt{2}+G_2^{\alpha\beta\gamma}(\xi',\tau))^3} \, d\xi' \bigg | \\ 
&\leq C(\theta) \, (-\tau)^{-\frac{1}{2}} \, (|H(\xi,\tau)|+|H(0,\tau)|)
\end{align*} 
for $|\xi| \leq \sqrt{4-\frac{\theta^2}{4}} \, (-\tau)^{\frac{1}{2}}$. This gives 
\begin{align*} 
&\bigg | \int_0^\xi \frac{(G_{1\xi}(\xi',\tau)+G_{2\xi}^{\alpha\beta\gamma}(\xi',\tau)) \, H_\xi(\xi',\tau)}{(\sqrt{2}+G_2^{\alpha\beta\gamma}(\xi',\tau))^2} \, d\xi' \bigg | \\ 
&\leq C(\theta) \, (-\tau)^{-\frac{1}{2}} \, \bigg | \int_0^\xi |H(\xi',\tau)| \, d\xi' \bigg | \\ 
&+ C(\theta) \, (-\tau)^{-\frac{1}{2}} \, (|H(\xi,\tau)|+|H(0,\tau)|) 
\end{align*} 
for $|\xi| \leq \sqrt{4-\frac{\theta^2}{4}} \, (-\tau)^{\frac{1}{2}}$. Consequently,  
\begin{align*} 
|E_{\mathcal{C},6}(\xi,\tau)| 
&\leq C(\theta) \, (-\tau)^{-\frac{1}{2}} \, |G_{2\xi}^{\alpha\beta\gamma}(\xi,\tau)| \, \chi_{\mathcal{C}}((-\tau)^{-\frac{1}{2}} \xi) \, \bigg | \int_0^\xi |H(\xi',\tau)| \, d\xi' \bigg | \\ 
&+ C(\theta) \, (-\tau)^{-\frac{1}{2}} \, |G_{2\xi}^{\alpha\beta\gamma}(\xi,\tau)| \, \chi_{\mathcal{C}}((-\tau)^{-\frac{1}{2}} \xi) \, (|H(\xi,\tau)|+|H(0,\tau)|)
\end{align*}
for $|\xi| \leq \sqrt{4-\frac{\theta^2}{4}} \, (-\tau)^{\frac{1}{2}}$. Since $\chi_{\mathcal{C}}$ is monotone decreasing on the interval $[0,\infty)$, we obtain $0 \leq \chi_{\mathcal{C}}((-\tau)^{-\frac{1}{2}} \xi) \leq \chi_{\mathcal{C}}((-\tau)^{-\frac{1}{2}} \xi')$ for $|\xi'| \leq |\xi|$. Putting these facts together, the assertion follows. \\

In order to estimate the term $\|E_{\mathcal{C},5}\|_{\mathcal{H},\infty,\tau_*}$, we need the following estimate for $H_\xi(0,\tau)$: 

\begin{lemma}
\label{derivative.of.H.at.0}
We have 
\[\sup_{\tau \leq \tau_*} \bigg ( \int_{\tau-1}^\tau H_\xi(0,\tau')^2 \, d\tau' \bigg )^{\frac{1}{2}} \leq C \, \|H_{\mathcal{C}}\|_{\mathcal{H},\infty,\tau_*} + C \sum_{k=1}^6 \|E_{\mathcal{C},k}\|_{\mathcal{H},\infty,\tau_*}.\]
\end{lemma}

\textbf{Proof.} 
In the region $\{|\xi| \leq 1\}$, we have $\frac{\partial}{\partial \tau} H_{\mathcal{C}} = \mathcal{L} H_{\mathcal{C}} + \sum_{k=1}^6 E_{\mathcal{C},k}$. Using the embedding of the Sobolev space $H^1([-1,1])$ into $C^0([-1,1])$ together with standard interior estimates for linear parabolic equations, we obtain 
\begin{align*} 
&\sup_{\tau \leq \tau_*} \bigg ( \int_{\tau-1}^\tau H_{\mathcal{C},\xi}(0,\tau')^2 \, d\tau' \bigg )^{\frac{1}{2}} \\ 
&\leq C \sup_{\tau \leq \tau_*} \bigg ( \int_{\tau-1}^\tau \int_{-1}^1 (H_{\mathcal{C},\xi\xi}(\xi,\tau')^2 + H_{\mathcal{C},\xi}(\xi,\tau')^2) \, d\xi \, d\tau' \bigg )^{\frac{1}{2}} \\ 
&\leq C \, \|H_{\mathcal{C}}\|_{\mathcal{H},\infty,\tau_*} + C \sum_{k=1}^6 \|E_{\mathcal{C},k}\|_{\mathcal{H},\infty,\tau_*}. 
\end{align*}
Since $H_{\mathcal{C},\xi}(0,\tau) = H_\xi(0,\tau)$, the assertion follows. \\

\begin{lemma}
\label{E1.E6}
We have 
\[\sum_{k=1}^6 \|E_{\mathcal{C},k}\|_{\mathcal{H},\infty,\tau_*} \leq C(\theta) \, (-\tau_*)^{-\frac{1}{2}} \, \|H_{\mathcal{C}}\|_{\mathcal{D},\infty,\tau_*}.\]
\end{lemma}

\textbf{Proof.} 
Using Proposition \ref{estimates.for.G_1.and.G_2}, we obtain the pointwise estimate 
\[|E_{\mathcal{C},1}(\xi,\tau)| \leq C(\theta) \, (-\tau)^{-\frac{1}{2}} \, (|\xi|+1) \, |H_{\mathcal{C}}(\xi,\tau)|\] 
for all $\tau \leq \tau_*$. Consequently, 
\[\|E_{\mathcal{C},1}\|_{\mathcal{H},\infty,\tau_*} \leq C(\theta) \, (-\tau_*)^{-\frac{1}{2}} \, \|H_{\mathcal{C}}\|_{\mathcal{D},\infty,\tau_*}\] 
by Proposition \ref{boundedness.of.operators}. In the next step, we estimate the terms $\|E_{\mathcal{C},k}\|_{\mathcal{H},\infty,\tau_*}$, where $k \in \{2,3,4,5\}$. Using the embedding of the Sobolev space $H^1([-1,1])$ into $C^0([-1,1])$, we can bound $|H(0,\tau)| \leq C \, \|H_{\mathcal{C}}(\cdot,\tau)\|_{\mathcal{D}}$. This implies 
\begin{align*} 
\sum_{k=2}^5 \|E_{\mathcal{C},k}\|_{\mathcal{H},\infty,\tau_*} 
&\leq C(\theta) \, (-\tau_*)^{-\frac{1}{2}} \, \|H_{\mathcal{C}}\|_{\mathcal{D},\infty,\tau_*} \\ 
&+ C(\theta) \, (-\tau_*)^{-\frac{1}{2}} \sup_{\tau \leq \tau_*} \bigg ( \int_{\tau-1}^\tau H_\xi(0,\tau')^2 \, d\tau' \bigg )^{\frac{1}{2}}. 
\end{align*}
Finally, using the pointwise estimate in Lemma \ref{pointwise.estimate.for.E6} together with Proposition \ref{boundedness.of.operators}, we obtain 
\[\|E_{\mathcal{C},6}\|_{\mathcal{H},\infty,\tau_*} \leq C(\theta) \, (-\tau_*)^{-\frac{1}{2}} \, \|H_{\mathcal{C}}\|_{\mathcal{D},\infty,\tau_*}.\] 
Putting these facts together, we conclude that 
\begin{align*} 
\sum_{k=1}^6 \|E_{\mathcal{C},k}\|_{\mathcal{H},\infty,\tau_*} 
&\leq C(\theta) \, (-\tau_*)^{-\frac{1}{2}} \, \|H_{\mathcal{C}}\|_{\mathcal{D},\infty,\tau_*} \\ 
&+ C(\theta) \, (-\tau_*)^{-\frac{1}{2}} \sup_{\tau \leq \tau_*} \bigg ( \int_{\tau-1}^\tau H_\xi(0,\tau')^2 \, d\tau' \bigg )^{\frac{1}{2}} \\ 
&\leq C(\theta) \, (-\tau_*)^{-\frac{1}{2}} \, \|H_{\mathcal{C}}\|_{\mathcal{D},\infty,\tau_*} \\ 
&+ C(\theta) \, (-\tau_*)^{-\frac{1}{2}} \, \sum_{k=1}^6 \|E_{\mathcal{C},k}\|_{\mathcal{H},\infty,\tau_*}, 
\end{align*}
where in the last step we have used Lemma \ref{derivative.of.H.at.0}. If $-\tau_*$ is sufficiently large, the last term on the right hand side can be absorbed into the left hand side. From this, the assertion follows. \\

\begin{lemma}
\label{E7.E10}
We have 
\begin{align*} 
&\sum_{k=7}^9 \|E_{\mathcal{C},k}\|_{\mathcal{H},\infty,\tau_*} + \|E_{\mathcal{C},10}\|_{\mathcal{D}^*,\infty,\tau_*} \\ 
&\leq C(\theta) \, (-\tau_*)^{-\frac{1}{2}} \, \Big \| H \, 1_{\{\sqrt{4-\frac{\theta^2}{2}} \, (-\tau)^{\frac{1}{2}} \leq |\xi| \leq \sqrt{4-\frac{\theta^2}{4}} \, (-\tau)^{\frac{1}{2}}\}} \Big \|_{\mathcal{H},\infty,\tau_*}. 
\end{align*}
\end{lemma}

\textbf{Proof.} 
Using Proposition \ref{boundedness.of.operators}, we obtain 
\[\|E_{\mathcal{C},10}\|_{\mathcal{D}^*,\infty,\tau_*} \leq C \, (-\tau_*)^{-\frac{1}{2}} \, \|\chi_{\mathcal{C}}'((-\tau)^{-\frac{1}{2}} \xi) \, H\|_{\mathcal{H},\infty,\tau_*}.\] 
This gives the desired estimate for $E_{\mathcal{C},10}$. The estimates for $E_{\mathcal{C},7}$, $E_{\mathcal{C},8}$, and $E_{\mathcal{C},9}$ follow directly from the respective definitions. This completes the proof of Lemma \ref{E7.E10}. \\

We now complete the proof of Proposition \ref{estimate.for.difference.in.cylindrical.region}. To that end, we apply Proposition \ref{linear.estimate} to the function $H_{\mathcal{C}}$. Since $P_+ H_{\mathcal{C}}(\tau_*) = 0$, we obtain 
\[\sup_{\tau \leq \tau_*} \|\hat{H}_{\mathcal{C}}(\tau)\|_{\mathcal{H}} + \Lambda^{-1} \, \|\hat{H}_{\mathcal{C}}\|_{\mathcal{D},\infty,\tau_*} \leq \Lambda \sum_{k=1}^{10} \|E_{\mathcal{C},k}(\xi,\tau)\|_{\mathcal{D}^*,\infty,\tau_*}.\] 
We use Lemma \ref{E1.E6} and Lemma \ref{E7.E10} to estimate the terms on the right hand side. This implies 
\begin{align*} 
&\sup_{\tau \leq \tau_*} \|\hat{H}_{\mathcal{C}}(\tau)\|_{\mathcal{H}} + \Lambda^{-1} \, \|\hat{H}_{\mathcal{C}}\|_{\mathcal{D},\infty,\tau_*} \\ 
&\leq C(\theta) \, (-\tau_*)^{-\frac{1}{2}} \, \|\hat{H}_{\mathcal{C}}\|_{\mathcal{D},\infty,\tau_*} + C(\theta) \, (-\tau_*)^{-\frac{1}{2}} \, \|P_0 H_{\mathcal{C}}\|_{\mathcal{D},\infty,\tau_*} \\ 
&+ C(\theta) \, (-\tau_*)^{-\frac{1}{2}} \, \Big \| H \, 1_{\{\sqrt{4-\frac{\theta^2}{2}} \, (-\tau)^{\frac{1}{2}} \leq |\xi| \leq \sqrt{4-\frac{\theta^2}{4}} \, (-\tau)^{\frac{1}{2}}\}} \Big \|_{\mathcal{H},\infty,\tau_*}. 
\end{align*}
If $-\tau_*$ is sufficiently large (depending on $\theta$), then the first term on the right hand side can be absorbed into the left hand side. This completes the proof of Proposition \ref{estimate.for.difference.in.cylindrical.region}. \\

\section{Analysis of the overlap region and proof of Proposition \ref{neutral.mode.dominates}}

\label{overlap.region}

In this section, we give the proof of Proposition \ref{neutral.mode.dominates}. We remind the reader that $\theta$ is chosen as in Proposition \ref{estimate.for.difference.in.tip.region}. As in the previous section, we write $H$, $H_{\mathcal{C}}$, $\hat{H}_{\mathcal{C}}$, and $a$ instead of $H^{\alpha\beta\gamma}$, $H_{\mathcal{C}}^{\alpha\beta\gamma}$, $\hat{H}_{\mathcal{C}}^{\alpha\beta\gamma}$, and $a^{\alpha\beta\gamma}$. We begin with an elementary lemma:

\begin{lemma}
\label{weighted.Poincare.inequality.cylindrical.region}
Assume that $4 \leq L_1 < L_2 < L_3$. Then 
\begin{align*} 
L_2^2 \int_{\{L_2 \leq \xi \leq L_3\}} e^{-\frac{\xi^2}{4}} \, f(\xi)^2 \, d\xi 
&\leq C \int_{\{L_1 \leq \xi \leq L_3\}} e^{-\frac{\xi^2}{4}} \, f'(\xi)^2 \, d\xi \\ 
&+ C \, (L_2-L_1)^{-2} \int_{\{L_1 \leq \xi \leq L_2\}} e^{-\frac{\xi^2}{4}} \, f(\xi)^2 \, d\xi, 
\end{align*}
where $C$ is a numerical constant that is independent of $L_1$, $L_2$, $L_3$, and $f$.
\end{lemma} 

\textbf{Proof.} 
Note that 
\begin{align*} 
\frac{d}{d\xi} (e^{-\frac{\xi^2}{4}} \, \xi \, f(\xi)^2) 
&= e^{-\frac{\xi^2}{4}} \, \Big ( f(\xi)^2 - \frac{\xi^2}{2} \, f(\xi)^2 + 2\xi \, f(\xi) \, f'(\xi) \Big ) \\ 
&\leq e^{-\frac{\xi^2}{4}} \, \Big ( f(\xi)^2 - \frac{\xi^2}{4} \, f(\xi)^2 + 4 \, f'(\xi)^2 \Big ). 
\end{align*} 
Integrating over $\xi \in [0,L_3]$ gives 
\begin{align*} 
&e^{-\frac{L_3^2}{4}} \, L_3 \, f(L_3)^2 + \int_{\{0 \leq \xi \leq L_3\}} e^{-\frac{\xi^2}{4}} \, \Big ( \frac{\xi^2}{4} - 1 \Big ) \, f(\xi)^2 \, d\xi \\ 
&\leq 4 \int_{\{0 \leq \xi \leq L_3\}} e^{-\frac{\xi^2}{4}} \, f'(\xi)^2 \, d\xi. 
\end{align*} 
Hence, if $f$ vanishes for $0 \leq \xi \leq 4$, then we obtain
\[\Big ( \frac{L_2^2}{4}-1 \Big ) \int_{\{L_2 \leq \xi \leq L_3\}} e^{-\frac{\xi^2}{4}} \, f(\xi)^2 \, d\xi \leq 4 \int_{\{4 \leq \xi \leq L_3\}} e^{-\frac{\xi^2}{4}} \, f'(\xi)^2 \, d\xi.\] 
Finally, we multiply the given function $f$ by a smooth cutoff function which is equal to $0$ on the interval $(-\infty,L_1]$, and which is equal to $1$ on the interval $[L_2,\infty)$. This completes the proof of Lemma \ref{weighted.Poincare.inequality.cylindrical.region}. \\

The following lemma relates the function $H(\xi,\tau)$ to the function $W_+(\rho,\tau)$:

\begin{lemma}
\label{pointwise.estimate.in.overlap.region}
If we choose $-\tau_*$ sufficiently large (depending on $\theta$), then 
\[\big | H_\xi(\xi,\tau) + W_+(\sqrt{2}+G_1(\xi,\tau),\tau) \big | \leq C(\theta) \, |H(\xi,\tau)|\] 
provided that $\sqrt{4-400 \, \theta^2} \, (-\tau)^{\frac{1}{2}} \leq \xi \leq \sqrt{4-\frac{\theta^2}{100}} \, (-\tau)^{\frac{1}{2}}$ and $\tau \leq \tau_*$. 
\end{lemma}

\textbf{Proof.} 
Suppose that $\sqrt{4-400 \, \theta^2} \, (-\tau)^{\frac{1}{2}} \leq \xi \leq \sqrt{4-\frac{\theta^2}{100}} \, (-\tau)^{\frac{1}{2}}$ and $\tau \leq \tau_*$. Let 
\begin{align*} 
\rho_1 &:= e^{\frac{\tau}{2}} \, F_1(e^{-\frac{\tau}{2}} \xi,-e^{-\tau}) = \sqrt{2}+G_1(\xi,\tau), \\ 
\rho_2 &:= e^{\frac{\tau}{2}} \, F_2^{\alpha\beta\gamma}(e^{-\frac{\tau}{2}} \xi,-e^{-\tau}) = \sqrt{2}+G_2^{\alpha\beta\gamma}(\xi,\tau). 
\end{align*} 
By Proposition \ref{precise.estimate.for.F} and Proposition \ref{estimate.for.modified.profile}, $\rho_1 \in [\frac{\theta}{20},20\theta]$ and $\rho_2 \in [\frac{\theta}{20},20\theta]$. Moreover, $\frac{\partial}{\partial \xi} G_1(\xi,\tau) = -V_{1+}(\rho_1,\tau)$, $\frac{\partial}{\partial \xi} G_2^{\alpha\beta\gamma}(\xi,\tau) = -V_{2+}^{\beta\gamma}(\rho_2,\tau)$, and $\rho_1 - \rho_2 = H(\xi,\tau)$. This implies 
\begin{align*} 
\frac{\partial}{\partial \xi} H(\xi,\tau) 
&= \frac{\partial}{\partial \xi} G_1(\xi,\tau) - \frac{\partial}{\partial \xi} G_2^{\alpha\beta\gamma}(\xi,\tau) \\ 
&= -V_{1+}(\rho_1,\tau) + V_{2+}^{\beta\gamma}(\rho_2,\tau) \\ 
&= -W_+(\rho_1,\tau) - V_{2+}^{\beta\gamma}(\rho_1,\tau) + V_{2+}^{\beta\gamma}(\rho_2,\tau). 
\end{align*} 
Using Proposition \ref{bound.for.V.in.transition.region}, we obtain $\big | \frac{\partial}{\partial \rho} V_{2+}^{\beta\gamma}(\rho,\tau) \big | \leq C(\theta)$ for every $\rho \in [\frac{\theta}{20},20\theta]$ and every $\tau \leq \tau_*$. This gives 
\begin{align*} 
\Big | \frac{\partial}{\partial \xi} H(\xi,\tau) + W_+(\rho_1,\tau) \Big | 
&\leq \big | V_{2+}^{\beta\gamma}(\rho_1,\tau) - V_{2+}^{\beta\gamma}(\rho_2,\tau) \big | \\ 
&\leq C(\theta) \, |\rho_1-\rho_2| \\ 
&\leq C(\theta) \, |H(\xi,\tau)|, 
\end{align*}
as claimed. \\

\begin{lemma} 
\label{overlap.1}
We have 
\begin{align*} 
&(-\tau) \int_{\{\sqrt{4-\frac{\theta^2}{2}} \, (-\tau)^{\frac{1}{2}} \leq \xi \leq \sqrt{4-\frac{\theta^2}{4}} \, (-\tau)^{\frac{1}{2}}\}} e^{-\frac{\xi^2}{4}} \, H(\xi,\tau)^2 \, d\xi \\ 
&\leq C(\theta) \, (-\tau)^{-\frac{1}{2}} \int_{\frac{\theta}{4}}^\theta V_{1+}(\rho,\tau)^{-2} \, W_+(\rho,\tau)^2 \, e^{\mu_+(\rho,\tau)} \, d\rho \\ 
&+ C(\theta) \int_{\{\sqrt{4-\theta^2} \, (-\tau)^{\frac{1}{2}} \leq \xi \leq \sqrt{4-\frac{\theta^2}{2}} \, (-\tau)^{\frac{1}{2}}\}} e^{-\frac{\xi^2}{4}} \, H(\xi,\tau)^2 \, d\xi,
\end{align*} 
provided that $\tau \leq \tau_*$ and $-\tau_*$ is sufficiently large.
\end{lemma}

\textbf{Proof.} 
We apply Lemma \ref{weighted.Poincare.inequality.cylindrical.region} with $L_1 = \sqrt{4-\theta^2} \, (-\tau)^{\frac{1}{2}}$, $L_2 = \sqrt{4-\frac{\theta^2}{2}} \, (-\tau)^{\frac{1}{2}}$, $L_3 = \sqrt{4-\frac{\theta^2}{4}} \, (-\tau)^{\frac{1}{2}}$, and $f(\xi) = H(\xi,\tau)$. This implies 
\begin{align*} 
&(-\tau) \int_{\{\sqrt{4-\frac{\theta^2}{2}} \, (-\tau)^{\frac{1}{2}} \leq \xi \leq \sqrt{4-\frac{\theta^2}{4}} \, (-\tau)^{\frac{1}{2}}\}} e^{-\frac{\xi^2}{4}} \, H(\xi,\tau)^2 \, d\xi \\ 
&\leq C(\theta) \int_{\{\sqrt{4-\theta^2} \, (-\tau)^{\frac{1}{2}} \leq \xi \leq \sqrt{4-\frac{\theta^2}{4}} \, (-\tau)^{\frac{1}{2}}\}} e^{-\frac{\xi^2}{4}} \, H_\xi(\xi,\tau)^2 \, d\xi \\ 
&+ C(\theta) \int_{\{\sqrt{4-\theta^2} \, (-\tau)^{\frac{1}{2}} \leq \xi \leq \sqrt{4-\frac{\theta^2}{2}} \, (-\tau)^{\frac{1}{2}}\}} e^{-\frac{\xi^2}{4}} \, H(\xi,\tau)^2 \, d\xi. 
\end{align*}
By Lemma \ref{pointwise.estimate.in.overlap.region}, $H_\xi(\xi,\tau)^2 \leq 4 \, W_+(\sqrt{2}+G_1(\xi,\tau),\tau)^2 + C(\theta) \, H(\xi,\tau)^2$ for $\sqrt{4-\theta^2} \, (-\tau)^{\frac{1}{2}} \leq \xi \leq \sqrt{4-\frac{\theta^2}{4}} \, (-\tau)^{\frac{1}{2}}$. This gives 
\begin{align*} 
&(-\tau) \int_{\{\sqrt{4-\frac{\theta^2}{2}} \, (-\tau)^{\frac{1}{2}} \leq \xi \leq \sqrt{4-\frac{\theta^2}{4}} \, (-\tau)^{\frac{1}{2}}\}} e^{-\frac{\xi^2}{4}} \, H(\xi,\tau)^2 \, d\xi \\ 
&\leq C(\theta) \int_{\{\sqrt{4-\theta^2} \, (-\tau)^{\frac{1}{2}} \leq \xi \leq \sqrt{4-\frac{\theta^2}{4}} \, (-\tau)^{\frac{1}{2}}\}} e^{-\frac{\xi^2}{4}} \, W_+(\sqrt{2}+G_1(\xi,\tau),\tau)^2 \, d\xi \\ 
&+ C(\theta) \int_{\{\sqrt{4-\theta^2} \, (-\tau)^{\frac{1}{2}} \leq \xi \leq \sqrt{4-\frac{\theta^2}{4}} \, (-\tau)^{\frac{1}{2}}\}} e^{-\frac{\xi^2}{4}} \, H(\xi,\tau)^2 \, d\xi.
\end{align*} 
By Proposition \ref{bound.for.V.in.transition.region}, we have $V_{1+}(\rho,\tau) \leq C(\theta) \, (-\tau)^{-\frac{1}{2}}$ for $\rho \in [\frac{\theta}{4},\theta]$. Moreover, Lemma \ref{weight} gives $\mu_+(\rho,\tau) = -\frac{\xi_{1+}(\rho,\tau)^2}{4}$ for $\rho \in [\frac{\theta}{4},\theta]$. Consequently, 
\begin{align*} 
&\int_{\{\sqrt{4-\theta^2} \, (-\tau)^{\frac{1}{2}} \leq \xi \leq \sqrt{4-\frac{\theta^2}{4}} \, (-\tau)^{\frac{1}{2}}\}} e^{-\frac{\xi^2}{4}} \, W_+(\sqrt{2}+G_1(\xi,\tau),\tau)^2 \, d\xi \\ 
&\leq \int_{\frac{\theta}{4}}^\theta V_{1+}(\rho,\tau)^{-1} \, W_+(\rho,\tau)^2 \, e^{-\frac{\xi_{1+}(\rho,\tau)^2}{4}} \, d\rho \\ 
&\leq C(\theta) \, (-\tau)^{-\frac{1}{2}} \int_{\frac{\theta}{4}}^\theta V_{1+}(\rho,\tau)^{-2} \, W_+(\rho,\tau)^2 \, e^{\mu_+(\rho,\tau)} \, d\rho. 
\end{align*} 
Putting these facts together, we conclude that 
\begin{align*} 
&(-\tau) \int_{\{\sqrt{4-\frac{\theta^2}{2}} \, (-\tau)^{\frac{1}{2}} \leq \xi \leq \sqrt{4-\frac{\theta^2}{4}} \, (-\tau)^{\frac{1}{2}}\}} e^{-\frac{\xi^2}{4}} \, H(\xi,\tau)^2 \, d\xi \\ 
&\leq C(\theta) \, (-\tau)^{-\frac{1}{2}} \int_{\frac{\theta}{4}}^\theta V_{1+}(\rho,\tau)^{-2} \, W_+(\rho,\tau)^2 \, e^{\mu_+(\rho,\tau)} \, d\rho \\ 
&+ C(\theta) \int_{\{\sqrt{4-\theta^2} \, (-\tau)^{\frac{1}{2}} \leq \xi \leq \sqrt{4-\frac{\theta^2}{4}} \, (-\tau)^{\frac{1}{2}}\}} e^{-\frac{\xi^2}{4}} \, H(\xi,\tau)^2 \, d\xi.
\end{align*} 
From this, the assertion follows easily. \\

\begin{lemma} 
\label{overlap.2}
We have 
\begin{align*} 
&(-\tau)^{-\frac{1}{2}} \int_\theta^{2\theta} V_{1+}(\rho,\tau)^{-2} \, W_+(\rho,\tau)^2 \, e^{\mu_+(\rho,\tau)} \, d\rho \\ 
&\leq C(\theta) \int_{\{\sqrt{4-16\theta^2} \, (-\tau)^{\frac{1}{2}} \leq \xi \leq \sqrt{4-\theta^2} \, (-\tau)^{\frac{1}{2}}\}} e^{-\frac{\xi^2}{4}} \, (H_\xi(\xi,\tau)^2 + H(\xi,\tau)^2) \, d\xi 
\end{align*} 
provided that $\tau \leq \tau_*$ and $-\tau_*$ is sufficiently large.
\end{lemma}

\textbf{Proof.} 
By Proposition \ref{bound.for.V.in.transition.region}, we have $V_{1+}(\rho,\tau) \geq \frac{1}{C(\theta)} \, (-\tau)^{-\frac{1}{2}}$ for $\rho \in [\theta,2\theta]$. Moreover, Lemma \ref{weight} gives $\mu_+(\rho,\tau) = -\frac{\xi_{1+}(\rho,\tau)^2}{4}$ for $\rho \in [\theta,2\theta]$. Consequently, 
\begin{align*} 
&\int_{\{\sqrt{4-16\theta^2} \, (-\tau)^{\frac{1}{2}} \leq \xi \leq \sqrt{4-\theta^2} \, (-\tau)^{\frac{1}{2}}\}} e^{-\frac{\xi^2}{4}} \, W_+(\sqrt{2}+G_1(\xi,\tau),\tau)^2 \, d\xi \\ 
&\geq \int_\theta^{2\theta} V_{1+}(\rho,\tau)^{-1} \, W_+(\rho,\tau)^2 \, e^{-\frac{\xi_{1+}(\rho,\tau)^2}{4}} \, d\rho \\ 
&\geq \frac{1}{C(\theta)} \, (-\tau)^{-\frac{1}{2}} \int_\theta^{2\theta} V_{1+}(\rho,\tau)^{-2} \, W_+(\rho,\tau)^2 \, e^{\mu_+(\rho,\tau)} \, d\rho. 
\end{align*} 
By Lemma \ref{pointwise.estimate.in.overlap.region}, $W_+(\sqrt{2}+G_1(\xi,\tau),\tau)^2 \leq 4 \, H_\xi(\xi,\tau)^2 + C(\theta) \, H(\xi,\tau)^2$ for $\sqrt{4-16\theta^2} \, (-\tau)^{\frac{1}{2}} \leq \xi \leq \sqrt{4-\theta^2} \, (-\tau)^{\frac{1}{2}}$. From this, the assertion follows. \\

\begin{proposition}
\label{control.overlap.region}
We have 
\begin{align*} 
&\sup_{\tau \leq \tau_*} (-\tau) \int_{\tau-1}^\tau \int_{\{\sqrt{4-\frac{\theta^2}{2}} \, (-\tau')^{\frac{1}{2}} \leq |\xi| \leq \sqrt{4-\frac{\theta^2}{4}} \, (-\tau')^{\frac{1}{2}}\}} e^{-\frac{\xi^2}{4}} \, H(\xi,\tau')^2 \, d\xi \, d\tau' \\ 
&\leq C(\theta) \sup_{\tau \leq \tau_*} \int_{\tau-1}^\tau \int_{\mathbb{R}} e^{-\frac{\xi^2}{4}} \, (H_{\mathcal{C},\xi}(\xi,\tau')^2 + H_{\mathcal{C}}(\xi,\tau')^2) \, d\xi \, d\tau'. 
\end{align*} 
\end{proposition}

\textbf{Proof.} 
Using Lemma \ref{overlap.1}, Proposition \ref{estimate.for.difference.in.tip.region}, and Lemma \ref{overlap.2}, we obtain
\begin{align*} 
&\sup_{\tau \leq \tau_*} (-\tau) \int_{\tau-1}^\tau \int_{\{\sqrt{4-\frac{\theta^2}{2}} \, (-\tau')^{\frac{1}{2}} \leq \xi \leq \sqrt{4-\frac{\theta^2}{4}} \, (-\tau')^{\frac{1}{2}}\}} e^{-\frac{\xi^2}{4}} \, H(\xi,\tau')^2 \, d\xi \, d\tau' \\ 
&\leq C(\theta) \sup_{\tau \leq \tau_*} (-\tau)^{-\frac{1}{2}} \int_{\tau-1}^\tau \int_{\frac{\theta}{4}}^\theta V_{1+}(\rho,\tau')^{-2} \, W_+(\rho,\tau')^2 \, e^{\mu_+(\rho,\tau')} \, d\rho \, d\tau' \\ 
&+ C(\theta) \sup_{\tau \leq \tau_*} \int_{\tau-1}^\tau \int_{\{\sqrt{4-\theta^2} \, (-\tau')^{\frac{1}{2}} \leq \xi \leq \sqrt{4-\frac{\theta^2}{2}} \, (-\tau')^{\frac{1}{2}}\}} e^{-\frac{\xi^2}{4}} \, H(\xi,\tau')^2 \, d\xi \, d\tau' \\ 
&\leq C(\theta) \sup_{\tau \leq \tau_*} (-\tau)^{-\frac{1}{2}} \int_{\tau-1}^\tau \int_\theta^{2\theta} V_{1+}(\rho,\tau')^{-2} \, W_+(\rho,\tau')^2 \, e^{\mu_+(\rho,\tau')} \, d\rho \, d\tau' \\ 
&+ C(\theta) \sup_{\tau \leq \tau_*} \int_{\tau-1}^\tau \int_{\{\sqrt{4-\theta^2} \, (-\tau')^{\frac{1}{2}} \leq \xi \leq \sqrt{4-\frac{\theta^2}{2}} \, (-\tau')^{\frac{1}{2}}\}} e^{-\frac{\xi^2}{4}} \, H(\xi,\tau')^2 \, d\xi \, d\tau' \\ 
&\leq C(\theta) \sup_{\tau \leq \tau_*} \int_{\tau-1}^\tau \int_{\mathbb{R}} e^{-\frac{\xi^2}{4}} \, (H_{\mathcal{C},\xi}(\xi,\tau')^2 + H_{\mathcal{C}}(\xi,\tau')^2) \, d\xi \, d\tau'. 
\end{align*} 
An analogous argument works when $\xi$ is negative. This completes the proof of Proposition \ref{control.overlap.region}. \\

After these preparations, we now finish the proof of Proposition \ref{neutral.mode.dominates}. Using Proposition \ref{control.overlap.region}, we obtain
\begin{align*} 
&\sup_{\tau \leq \tau_*} (-\tau) \int_{\tau-1}^\tau \int_{\{\sqrt{4-\frac{\theta^2}{2}} \, (-\tau')^{\frac{1}{2}} \leq |\xi| \leq \sqrt{4-\frac{\theta^2}{4}} \, (-\tau')^{\frac{1}{2}}\}} e^{-\frac{\xi^2}{4}} \, H(\xi,\tau')^2 \, d\xi \, d\tau' \\ 
&\leq C(\theta) \sup_{\tau \leq \tau_*} \int_{\tau-1}^\tau a(\tau')^2 \, d\tau' \\ 
&+ C(\theta) \sup_{\tau \leq \tau_*} \int_{\tau-1}^\tau \int_{\mathbb{R}} e^{-\frac{\xi^2}{4}} \, (\hat{H}_{\mathcal{C},\xi}(\xi,\tau')^2 + \hat{H}_{\mathcal{C}}(\xi,\tau')^2) \, d\xi \, d\tau'. 
\end{align*} 
Combining this estimate with Proposition \ref{estimate.for.difference.in.cylindrical.region} gives 
\begin{align*}
&(-\tau_*) \sup_{\tau \leq \tau_*} \int_{\tau-1}^\tau \int_{\mathbb{R}} e^{-\frac{\xi^2}{4}} \, (\hat{H}_{\mathcal{C},\xi}(\xi,\tau')^2 + \hat{H}_{\mathcal{C}}(\xi,\tau')^2) \, d\xi \, d\tau' \\ 
&\leq C(\theta) \sup_{\tau \leq \tau_*} \int_{\tau-1}^\tau a(\tau')^2 \, d\tau' \\ 
&+ C(\theta) \, (-\tau_*)^{-1} \, \sup_{\tau \leq \tau_*} \int_{\tau-1}^\tau \int_{\mathbb{R}} e^{-\frac{\xi^2}{4}} \, (\hat{H}_{\mathcal{C},\xi}(\xi,\tau')^2 + \hat{H}_{\mathcal{C}}(\xi,\tau')^2) \, d\xi \, d\tau'. 
\end{align*} 
If $-\tau_*$ is chosen sufficiently large (depending on $\theta$), then the last term on the right hand side can be absorbed into the left hand side. This completes the proof of Proposition \ref{neutral.mode.dominates}. \\

\section{Analysis of the neutral mode and proof of Proposition \ref{ode.for.a}}

\label{analysis.of.neutral.mode}

In this final section, we give the proof of Proposition \ref{ode.for.a}. As before, we write $H$, $H_{\mathcal{C}}$, $\hat{H}_{\mathcal{C}}$, and $a$ instead of $H^{\alpha\beta\gamma}$, $H_{\mathcal{C}}^{\alpha\beta\gamma}$, $\hat{H}_{\mathcal{C}}^{\alpha\beta\gamma}$, and $a^{\alpha\beta\gamma}$. 

\begin{lemma}
\label{control.overlap.region.2}
We have 
\begin{align*} 
&\sup_{\tau \leq \tau_*} (-\tau) \int_{\tau-1}^\tau \int_{\{\sqrt{4-\frac{\theta^2}{2}} \, (-\tau')^{\frac{1}{2}} \leq |\xi| \leq \sqrt{4-\frac{\theta^2}{4}} \, (-\tau')^{\frac{1}{2}}\}} e^{-\frac{\xi^2}{4}} \, H(\xi,\tau')^2 \, d\xi \, d\tau' \\ 
&\leq C(\theta) \sup_{\tau \leq \tau_*} \int_{\tau-1}^\tau a(\tau')^2 \, d\tau'. 
\end{align*} 
\end{lemma}

\textbf{Proof.} 
This follows by combining Proposition \ref{neutral.mode.dominates} and Proposition \ref{control.overlap.region}. \\

We next establish an improved version of Lemma \ref{derivative.of.H.at.0}: 

\begin{lemma}
\label{derivative.of.H.at.0.improved.version}
We have 
\[(-\tau_*) \sup_{\tau \leq \tau_*} \int_{\tau-1}^\tau H_\xi(0,\tau')^2 \, d\tau' \leq C(\theta) \sup_{\tau \leq \tau_*} \int_{\tau-1}^\tau a(\tau')^2 \, d\tau'.\] 
\end{lemma}

\textbf{Proof.}
The function $H_{\mathcal{C}}$ satisfies the evolution equation $\frac{\partial}{\partial \tau} H_{\mathcal{C}} = \mathcal{L} H_{\mathcal{C}} + \sum_{k=1}^{10} E_{\mathcal{C},k}$. Consequently, the function $\hat{H}_{\mathcal{C}}$ satisfies the evolution equation $\frac{\partial}{\partial \tau} \hat{H}_{\mathcal{C}} = \mathcal{L} \hat{H}_{\mathcal{C}} + \sum_{k=1}^{10} E_{\mathcal{C},k} - \sum_{k=1}^{10} P_0 E_{\mathcal{C},k}$. In particular, in the region $\{|\xi| \leq 1\}$, we have $\frac{\partial}{\partial \tau} \hat{H}_{\mathcal{C}} = \mathcal{L} \hat{H}_{\mathcal{C}} + \sum_{k=1}^6 E_{\mathcal{C},k} - \sum_{k=1}^{10} P_0 E_{\mathcal{C},k}$. Using the embedding of the Sobolev space $H^1([-1,1])$ into $C^0([-1,1])$ together with standard interior estimates for linear parabolic equations, we obtain 
\begin{align*}
&\sup_{\tau \leq \tau_*} \bigg ( \int_{\tau-1}^\tau \hat{H}_{\mathcal{C},\xi}(0,\tau')^2 \, d\tau' \bigg )^{\frac{1}{2}} \\ 
&\leq C \sup_{\tau \leq \tau_*} \bigg ( \int_{\tau-1}^\tau \int_{-1}^1 (\hat{H}_{\mathcal{C},\xi\xi}(\xi,\tau')^2 + \hat{H}_{\mathcal{C},\xi}(\xi,\tau')^2) \, d\xi \, d\tau' \bigg )^{\frac{1}{2}} \\ 
&\leq C \, \|\hat{H}_{\mathcal{C}}\|_{\mathcal{H},\infty,\tau_*} + C \sum_{k=1}^6 \|E_{\mathcal{C},k}\|_{\mathcal{H},\infty,\tau_*} + C \sum_{k=1}^{10} \|P_0 E_{\mathcal{C},k}\|_{\mathcal{H},\infty,\tau_*} \\ 
&\leq C \, \|\hat{H}_{\mathcal{C}}\|_{\mathcal{H},\infty,\tau_*} + C \sum_{k=1}^6 \|E_{\mathcal{C},k}\|_{\mathcal{H},\infty,\tau_*} + C \sum_{k=7}^{10} \|E_{\mathcal{C},k}\|_{\mathcal{D}^*,\infty,\tau_*}.
\end{align*} 
Note that $\hat{H}_{\mathcal{C},\xi}(0,\tau) = H_{\mathcal{C},\xi}(0,\tau) = H_\xi(0,\tau)$ for each $\tau$. In the next step, we use Lemma \ref{E1.E6} and Lemma \ref{E7.E10} to estimate the terms on the right hand side. This gives 
\begin{align*}
&\sup_{\tau \leq \tau_*} \bigg ( \int_{\tau-1}^\tau H_\xi(0,\tau')^2 \, d\tau' \bigg )^{\frac{1}{2}} \\ 
&\leq C \, \|\hat{H}_{\mathcal{C}}\|_{\mathcal{H},\infty,\tau_*} + C(\theta) \, (-\tau_*)^{-\frac{1}{2}} \, \|H_{\mathcal{C}}\|_{\mathcal{D},\infty,\tau_*} \\ 
&+ C(\theta) \, (-\tau_*)^{-\frac{1}{2}} \, \Big \| H \, 1_{\{\sqrt{4-\frac{\theta^2}{2}} \, (-\tau)^{\frac{1}{2}} \leq |\xi| \leq \sqrt{4-\frac{\theta^2}{4}} \, (-\tau)^{\frac{1}{2}}\}} \Big \|_{\mathcal{H},\infty,\tau_*} \\ 
&\leq C(\theta) \, (-\tau_*)^{-\frac{1}{2}} \, \sup_{\tau \leq \tau_*} \bigg ( \int_{\tau-1}^\tau a(\tau')^2 \, d\tau' \bigg )^{\frac{1}{2}}, 
\end{align*} 
where in the last step we have used Proposition \ref{neutral.mode.dominates} and Lemma \ref{control.overlap.region.2}. This completes the proof of Lemma \ref{derivative.of.H.at.0.improved.version}. \\

After these preparations, we now study the evolution of the function $a(\tau)$. Using the evolution equation $\frac{\partial}{\partial \tau} H_{\mathcal{C}} = \mathcal{L} H_{\mathcal{C}} + \sum_{k=1}^{10} E_{\mathcal{C},k}$, we obtain 
\[\frac{d}{d\tau} a(\tau) = \sum_{k=1}^{10} I_k(\tau),\] 
where $I_k(\tau)$ is defined by 
\[I_k(\tau) = \frac{1}{16\sqrt{2\pi}} \int_{\mathbb{R}} e^{-\frac{\xi^2}{4}} \, (\xi^2-2) \, E_{\mathcal{C},k}(\xi,\tau) \, d\xi\] 
for $\tau \leq \tau_*$. In the remainder of this section, we estimate the terms $I_k(\tau)$. \\

\begin{lemma} 
\label{I1}
Let $\delta > 0$ be given. If $-\tau_*$ is sufficiently large (depending on $\delta$), then  
\[\sup_{\tau \leq \tau_*} (-\tau) \int_{\tau-1}^\tau |I_1(\tau') - (-\tau')^{-1} \, a(\tau')| \, d\tau' \leq \delta \, \sup_{\tau \leq \tau_*} \bigg ( \int_{\tau-1}^\tau a(\tau')^2 \, d\tau' \bigg )^{\frac{1}{2}}.\] 
\end{lemma}

\textbf{Proof.}
We define a function $\hat{I}_1(\tau)$ by  
\begin{align*} 
\hat{I}_1(\tau) 
&= \frac{1}{16\sqrt{2\pi}} \int_{\{|\xi| \leq \sqrt{4-\frac{\theta^2}{4}} \, (-\tau)^{\frac{1}{2}}\}} e^{-\frac{\xi^2}{4}} \, (\xi^2-2) \, \hat{H}_{\mathcal{C}}(\xi,\tau) \\ 
&\hspace{25mm} \cdot \Big [ (\sqrt{2}+G_1(\xi,\tau))^{-1} (\sqrt{2}+G_2^{\alpha\beta\gamma}(\xi,\tau))^{-1} - \frac{1}{2} \Big ] \, d\xi 
\end{align*}
for $\tau \leq \tau_*$. Using the asymptotic estimates in Proposition \ref{estimates.for.G_1.and.G_2} together with the Cauchy-Schwarz inequality, we obtain 
\[(-\tau) \int_{\tau-1}^\tau |\hat{I}_1(\tau')| \, d\tau' \leq C(\theta) \, \bigg ( \int_{\tau-1}^\tau \int_{\mathbb{R}} e^{-\frac{\xi^2}{4}} \, \hat{H}_{\mathcal{C}}(\xi,\tau')^2 \, d\xi \, d\tau' \bigg )^{\frac{1}{2}}\] 
for all $\tau \leq \tau_*$. Using Proposition \ref{neutral.mode.dominates}, we conclude that 
\[\sup_{\tau \leq \tau_*} (-\tau) \int_{\tau-1}^\tau |\hat{I}_1(\tau')| \, d\tau' \leq \delta \, \sup_{\tau \leq \tau_*} \bigg ( \int_{\tau-1}^\tau a(\tau')^2 \, d\tau' \bigg )^{\frac{1}{2}},\] 
provided that $-\tau_*$ is sufficiently large (depending on $\delta$). On the other hand, using the identity $H_{\mathcal{C}}(\xi,\tau) - \hat{H}_{\mathcal{C}}(\xi,\tau) = \sqrt{2} \, a(\tau) \, (\xi^2-2)$, we obtain 
\begin{align*} 
I_1(\tau) - \hat{I}_1(\tau)
&= \frac{a(\tau)}{16\sqrt{\pi}} \int_{\{|\xi| \leq \sqrt{4-\frac{\theta^2}{4}} \, (-\tau)^{\frac{1}{2}}\}} e^{-\frac{\xi^2}{4}} \, (\xi^2-2)^2 \\ 
&\hspace{25mm} \cdot \Big [ (\sqrt{2}+G_1(\xi,\tau))^{-1} (\sqrt{2}+G_2^{\alpha\beta\gamma}(\xi,\tau))^{-1} - \frac{1}{2} \Big ] \, d\xi 
\end{align*} 
for all $\tau \leq \tau_*$. Using the asymptotic estimates in Proposition \ref{estimates.for.G_1.and.G_2} and the formula $\int_{\mathbb{R}} e^{-\frac{\xi^2}{4}} \, (\xi^2-2)^3 \, d\xi = 128\sqrt{\pi}$, we conclude that 
\[|I_1(\tau) - \hat{I}_1(\tau) - (-\tau)^{-1} \, a(\tau)| \leq \delta \, (-\tau)^{-1} \, |a(\tau)|\] 
for all $\tau \leq \tau_*$. Putting these facts together, the assertion follows. \\

\begin{lemma} 
\label{I2}
Let $\delta > 0$ be given. If $-\tau_*$ is sufficiently large (depending on $\delta$), then  
\[\sup_{\tau \leq \tau_*} (-\tau) \int_{\tau-1}^\tau |I_2(\tau')| \, d\tau' \leq \delta \, \sup_{\tau \leq \tau_*} \bigg ( \int_{\tau-1}^\tau a(\tau')^2 \, d\tau' \bigg )^{\frac{1}{2}}.\] 
\end{lemma}

\textbf{Proof.}
Using the Cauchy-Schwarz inequality, we obtain 
\begin{align*} 
&\int_{\tau-1}^\tau |I_2(\tau')| \, d\tau' \\ 
&\leq C(\theta) \int_{\tau-1}^\tau \int_{\mathbb{R}} e^{-\frac{\xi^2}{4}} \, |\xi^2-2| \, |G_{1\xi}(\xi,\tau')|^2 \, |H_{\mathcal{C}}(\xi,\tau')| \, d\xi \, d\tau' \\ 
&\leq C(\theta) \, \bigg ( \int_{\tau-1}^\tau \int_{\{|\xi| \leq \sqrt{4-\frac{\theta^2}{4}} \, (-\tau')^{\frac{1}{2}}\}} e^{-\frac{\xi^2}{4}} \, |\xi^2-2|^2 \, |G_{1\xi}(\xi,\tau')|^4 \, d\xi \, d\tau' \bigg )^{\frac{1}{2}} \\ 
&\hspace{10mm} \cdot \bigg ( \int_{\tau-1}^\tau \int_{\mathbb{R}} e^{-\frac{\xi^2}{4}} \, H_{\mathcal{C}}(\xi,\tau')^2 \, d\xi \, d\tau' \bigg )^{\frac{1}{2}} 
\end{align*}
for all $\tau \leq \tau_*$. To bound the term on the right hand side, we use the asymptotic estimates in Proposition \ref{estimates.for.G_1.and.G_2}. This gives 
\[(-\tau) \int_{\tau-1}^\tau |I_2(\tau')| \, d\tau' \leq \delta \, \bigg ( \int_{\tau-1}^\tau \int_{\mathbb{R}} e^{-\frac{\xi^2}{4}} \, H_{\mathcal{C}}(\xi,\tau')^2 \, d\xi \, d\tau' \bigg )^{\frac{1}{2}}\] 
for all $\tau \leq \tau_*$. Therefore, the assertion follows from Proposition \ref{neutral.mode.dominates}. \\

\begin{lemma} 
\label{I3}
Let $\delta > 0$ be given. If $-\tau_*$ is sufficiently large (depending on $\delta$), then  
\[\sup_{\tau \leq \tau_*} (-\tau) \int_{\tau-1}^\tau |I_3(\tau') - (-\tau')^{-1} \, a(\tau')| \, d\tau' \leq \delta \, \sup_{\tau \leq \tau_*} \bigg ( \int_{\tau-1}^\tau a(\tau')^2 \, d\tau' \bigg )^{\frac{1}{2}}.\] 
\end{lemma}

\textbf{Proof.}
We define a function $\hat{I}_3(\tau)$ by  
\begin{align*} 
\hat{I}_3(\tau)  
&= -\frac{1}{16\sqrt{2\pi}} \int_{\{|\xi| \leq \sqrt{4-\frac{\theta^2}{4}} \, (-\tau)^{\frac{1}{2}}\}} e^{-\frac{\xi^2}{4}} \, (\xi^2-2) \, \hat{H}_{\mathcal{C},\xi}(\xi,\tau) \\ 
&\hspace{25mm} \cdot (\sqrt{2}+G_2^{\alpha\beta\gamma}(\xi,\tau))^{-1} \, (G_{1\xi}(\xi,\tau) + G_{2\xi}^{\alpha\beta\gamma}(\xi,\tau)) \, d\xi 
\end{align*} 
for $\tau \leq \tau_*$. Using the asymptotic estimates in Proposition \ref{estimates.for.G_1.and.G_2}, we obtain 
\[(-\tau) \int_{\tau-1}^\tau |\hat{I}_3(\tau')| \, d\tau' \leq C(\theta) \, \bigg ( \int_{\tau-1}^\tau \int_{\mathbb{R}} e^{-\frac{\xi^2}{4}} \, \hat{H}_{\mathcal{C},\xi}(\xi,\tau')^2 \, d\xi \, d\tau' \bigg )^{\frac{1}{2}}\] 
for all $\tau \leq \tau_*$. Using Proposition \ref{neutral.mode.dominates}, we conclude that 
\[\sup_{\tau \leq \tau_*} (-\tau) \int_{\tau-1}^\tau |\hat{I}_3(\tau')| \, d\tau' \leq \delta \, \sup_{\tau \leq \tau_*} \bigg ( \int_{\tau-1}^\tau a(\tau')^2 \, d\tau' \bigg )^{\frac{1}{2}},\] 
provided that $-\tau_*$ is sufficiently large (depending on $\delta$). On the other hand, using the identity $H_{\mathcal{C},\xi}(\xi,\tau) - \hat{H}_{\mathcal{C},\xi}(\xi,\tau) = 2\sqrt{2} \, a(\tau) \, \xi$, we obtain 
\begin{align*} 
I_3(\tau) - \hat{I}_3(\tau) 
&= -\frac{a(\tau)}{8\sqrt{\pi}} \int_{\{|\xi| \leq \sqrt{4-\frac{\theta^2}{4}} \, (-\tau)^{\frac{1}{2}}\}}  e^{-\frac{\xi^2}{4}} \, (\xi^2-2) \, \xi \\ 
&\hspace{25mm} \cdot (\sqrt{2}+G_2^{\alpha\beta\gamma}(\xi,\tau))^{-1} \, (G_{1\xi}(\xi,\tau) + G_{2\xi}^{\alpha\beta\gamma}(\xi,\tau)) \, d\xi 
\end{align*} 
for all $\tau \leq \tau_*$. Using the asymptotic estimates in Proposition \ref{estimates.for.G_1.and.G_2} and the formula $\int_{\mathbb{R}} e^{-\frac{\xi^2}{4}} \, (\xi^2-2) \, \xi^2 \, d\xi = 16\sqrt{\pi}$, we conclude that 
\[|I_3(\tau) - \hat{I}_3(\tau) - (-\tau)^{-1} \, a(\tau)| \leq \delta \, (-\tau)^{-1} \, |a(\tau)|\] 
for all $\tau \leq \tau_*$. Putting these facts together, the assertion follows. \\

\begin{lemma} 
\label{I4}
Let $\delta > 0$ be given. If $-\tau_*$ is sufficiently large (depending on $\delta$), then 
\[\sup_{\tau \leq \tau_*} (-\tau) \int_{\tau-1}^\tau |I_4(\tau')| \, d\tau' \leq \delta \, \sup_{\tau \leq \tau_*} \bigg ( \int_{\tau-1}^\tau a(\tau')^2 \, d\tau' \bigg )^{\frac{1}{2}}.\] 
\end{lemma}

\textbf{Proof.}
Using Proposition \ref{estimates.for.G_1.and.G_2}, we obtain $|G_{1\xi}(0,\tau)| \leq o(1) \, (-\tau)^{-1}$. This implies 
\[|E_{\mathcal{C},4}(\xi,\tau)| \leq o(1) \, (-\tau)^{-1} \, |H_{\mathcal{C},\xi}(\xi,\tau)| + C(\theta) \, (-\tau)^{-2} \, (|\xi|^4+1) \, |H_{\mathcal{C},\xi}(\xi,\tau)|\] 
for all $\tau \leq \tau_*$. Using the Cauchy-Schwarz inequality, we conclude that 
\[(-\tau) \int_{\tau-1}^\tau |I_4(\tau')| \, d\tau' \leq \delta \, \bigg ( \int_{\tau-1}^\tau \int_{\mathbb{R}} e^{-\frac{\xi^2}{4}} \, H_{\mathcal{C},\xi}(\xi,\tau')^2 \, d\xi \, d\tau' \bigg )^{\frac{1}{2}}\] 
for all $\tau \leq \tau_*$. Hence, the assertion follows from Proposition \ref{neutral.mode.dominates}. \\

\begin{lemma} 
\label{I5}
Let $\delta > 0$ be given. If $-\tau_*$ is sufficiently large (depending on $\delta$), then 
\[\sup_{\tau \leq \tau_*} (-\tau) \int_{\tau-1}^\tau |I_5(\tau')| \, d\tau' \leq \delta \, \sup_{\tau \leq \tau_*} \bigg ( \int_{\tau-1}^\tau a(\tau')^2 \, d\tau' \bigg )^{\frac{1}{2}}.\] 
\end{lemma}

\textbf{Proof.}
Using the asymptotic estimates in Proposition \ref{estimates.for.G_1.and.G_2}, we obtain 
\begin{align*} 
&\int_{\tau-1}^\tau |I_5(\tau')| \, d\tau' \\ 
&\leq C(\theta) \, (-\tau)^{-1} \int_{\tau-1}^\tau |H_\xi(0,\tau')| \, d\tau' + C(\theta) \, (-\tau)^{-2} \int_{\tau-1}^\tau |H(0,\tau')| \, d\tau' 
\end{align*} 
for all $\tau \leq \tau_*$. The first term on the right hand side can be estimated using Lemma \ref{derivative.of.H.at.0.improved.version}. To estimate the second term on the right hand side, we use Proposition \ref{neutral.mode.dominates} together with the embedding of the Sobolev space $H^1([-1,1])$ into $C^0([-1,1])$. Putting these facts together, the assertion follows. \\

\begin{lemma} 
\label{I6}
Let $\delta > 0$ be given. If $-\tau_*$ is sufficiently large (depending on $\delta$), then 
\[\sup_{\tau \leq \tau_*} (-\tau) \int_{\tau-1}^\tau |I_6(\tau')| \, d\tau' \leq \delta \, \sup_{\tau \leq \tau_*} \bigg ( \int_{\tau-1}^\tau a(\tau')^2 \, d\tau' \bigg )^{\frac{1}{2}}.\] 
\end{lemma}

\textbf{Proof.}
For abbreviation, let $M(\xi,\tau) := \big | \int_0^\xi |H_{\mathcal{C}}(\xi',\tau)| \, d\xi' \big | + |H_{\mathcal{C}}(\xi,\tau)| + |H(0,\tau)|$. Proposition \ref{boundedness.of.operators} implies 
\[\int_{\tau-1}^\tau \int_{\mathbb{R}} e^{-\frac{\xi^2}{4}} \, M(\xi,\tau')^2 \, d\xi \, d\tau' \leq C \int_{\tau-1}^\tau \int_{\mathbb{R}} e^{-\frac{\xi^2}{4}} \, (H_{\mathcal{C},\xi}(\xi,\tau')^2+H_{\mathcal{C}}(\xi,\tau')^2) \, d\xi \, d\tau'\] 
for all $\tau \leq \tau_*$. Using Lemma \ref{pointwise.estimate.for.E6}, we obtain 
\[(-\tau)^{\frac{1}{2}} \, |E_{\mathcal{C},6}(\xi,\tau)| \leq C(\theta) \, |G_{2\xi}^{\alpha\beta\gamma}(\xi,\tau)| \, M(\xi,\tau)\] 
for all $\tau \leq \tau_*$. This gives 
\begin{align*} 
&(-\tau)^{\frac{1}{2}} \int_{\tau-1}^\tau |I_6(\tau')| \, d\tau' \\ 
&\leq C(\theta) \int_{\tau-1}^\tau \int_{\{|\xi| \leq \sqrt{4-\frac{\theta^2}{4}} \, (-\tau')^{\frac{1}{2}}\}}  e^{-\frac{\xi^2}{4}} \, |\xi^2-2| \, |G_{2\xi}^{\alpha\beta\gamma}(\xi,\tau')| \, M(\xi,\tau') \, d\xi \, d\tau' 
\end{align*} 
for all $\tau \leq \tau_*$. Using the Cauchy-Schwarz inequality, we conclude that 
\begin{align*} 
&(-\tau)^{\frac{1}{2}} \int_{\tau-1}^\tau |I_6(\tau')| \, d\tau' \\ 
&\leq C(\theta) \, \bigg ( \int_{\tau-1}^\tau \int_{\{|\xi| \leq \sqrt{4-\frac{\theta^2}{4}} \, (-\tau')^{\frac{1}{2}}\}} e^{-\frac{\xi^2}{4}} \, |\xi^2-2|^2 \, |G_{2\xi}^{\alpha\beta\gamma}(\xi,\tau')|^2 \, d\xi \, d\tau' \bigg )^{\frac{1}{2}} \\ 
&\hspace{10mm} \cdot \bigg ( \int_{\tau-1}^\tau \int_{\mathbb{R}} e^{-\frac{\xi^2}{4}} \, (H_{\mathcal{C},\xi}(\xi,\tau')^2+H_{\mathcal{C}}(\xi,\tau')^2) \, d\xi \, d\tau' \bigg )^{\frac{1}{2}} 
\end{align*} 
for all $\tau \leq \tau_*$. To bound the term on the right hand side, we use the asymptotic estimates in Proposition \ref{estimates.for.G_1.and.G_2}. This gives 
\[(-\tau) \int_{\tau-1}^\tau |I_6(\tau')| \, d\tau' \leq \delta \, \bigg ( \int_{\tau-1}^\tau \int_{\mathbb{R}} e^{-\frac{\xi^2}{4}} \, (H_{\mathcal{C},\xi}(\xi,\tau')^2+H_{\mathcal{C}}(\xi,\tau')^2) \, d\xi \, d\tau' \bigg )^{\frac{1}{2}}\] 
for all $\tau \leq \tau_*$. Therefore, the assertion follows from Proposition \ref{neutral.mode.dominates}. \\

\begin{lemma} 
\label{I7.I10}
Let $\delta > 0$ be given. If $-\tau_*$ is sufficiently large (depending on $\delta$), then 
\[\sup_{\tau \leq \tau_*} (-\tau) \int_{\tau-1}^\tau \sum_{k=7}^{10} |I_k(\tau')| \, d\tau' \leq \delta \, \sup_{\tau \leq \tau_*} \bigg ( \int_{\tau-1}^\tau a(\tau')^2 \, d\tau' \bigg )^{\frac{1}{2}}.\] 
\end{lemma}

\textbf{Proof.}
We first observe that 
\begin{align*} 
&\int_{\tau-1}^\tau \sum_{k=7}^{10} |I_k(\tau')| \, d\tau' \\ 
&\leq C(\theta) \int_{\tau-1}^\tau \int_{\{\sqrt{4-\frac{\theta^2}{2}} \, (-\tau')^{\frac{1}{2}} \leq |\xi| \leq \sqrt{4-\frac{\theta^2}{4}} \, (-\tau')^{\frac{1}{2}}\}} e^{-\frac{\xi^2}{4}} \, |\xi|^{10} \, |H(\xi,\tau')| \, d\xi \, d\tau' 
\end{align*}
for all $\tau \leq \tau_*$. Indeed, the estimates for $I_7$, $I_8$, and $I_9$ follow directly from the respective definitions. The estimate for $I_{10}$ follows by integration by parts. 

Using the Cauchy-Schwarz inequality, we obtain 
\begin{align*} 
&\int_{\tau-1}^\tau \sum_{k=7}^{10} |I_k(\tau')| \, d\tau' \\ 
&\leq C(\theta) \, \bigg ( \int_{\tau-1}^\tau \int_{\{\sqrt{4-\frac{\theta^2}{2}} \, (-\tau')^{\frac{1}{2}} \leq |\xi| \leq \sqrt{4-\frac{\theta^2}{4}} \, (-\tau')^{\frac{1}{2}}\}} e^{-\frac{\xi^2}{4}} \, |\xi|^{20} \, d\xi \bigg )^{\frac{1}{2}} \\ 
&\hspace{10mm} \cdot \bigg ( \int_{\tau-1}^\tau \int_{\{\sqrt{4-\frac{\theta^2}{2}} \, (-\tau')^{\frac{1}{2}} \leq |\xi| \leq \sqrt{4-\frac{\theta^2}{4}} \, (-\tau')^{\frac{1}{2}}\}} e^{-\frac{\xi^2}{4}} \, H(\xi,\tau')^2 \, d\xi \, d\tau' \bigg )^{\frac{1}{2}}
\end{align*} 
for all $\tau \leq \tau_*$. Consequently, 
\begin{align*} 
&(-\tau) \int_{\tau-1}^\tau \sum_{k=7}^{10} |I_k(\tau')| \, d\tau' \\ 
&\leq C(\theta) \, \bigg ( \int_{\tau-1}^\tau \int_{\{\sqrt{4-\frac{\theta^2}{2}} \, (-\tau')^{\frac{1}{2}} \leq |\xi| \leq \sqrt{4-\frac{\theta^2}{4}} \, (-\tau')^{\frac{1}{2}}\}} e^{-\frac{\xi^2}{4}} \, H(\xi,\tau')^2 \, d\xi \, d\tau' \bigg )^{\frac{1}{2}} 
\end{align*}
for all $\tau \leq \tau_*$. Hence, the assertion follows from Lemma \ref{control.overlap.region.2}. This completes the proof of Lemma \ref{I7.I10}. \\

Proposition \ref{ode.for.a} follows immediately from Lemma \ref{I1} -- Lemma \ref{I7.I10} together with the identity $\frac{d}{d\tau} a(\tau) = \sum_{k=1}^{10} I_k(\tau)$.

\end{document}